\documentclass[a4paper,12pt]{article}
\usepackage[utf8]{inputenc}
\usepackage[left=2.5cm, right=2.5cm, top=2.5cm, bottom=2.5cm]{geometry}
\usepackage{subcaption}
\usepackage{siunitx}
\usepackage{tikz}
\usetikzlibrary{arrows,decorations.markings,arrows.meta}
\usetikzlibrary{patterns.meta}

\usepackage{amssymb}
\usepackage{amsmath,amsfonts,bm,enumitem}
\usepackage{booktabs, multirow} 
\usepackage{soul}
\usepackage{graphicx}
\usepackage{changepage,threeparttable} 
\newcommand{\zetabis}{s}
\DeclareMathOperator{\supp}{supp}

 \newtheorem{theorem}{Theorem}[section]
 
 \newtheorem{proposition}[theorem]{Proposition}

 \newtheorem{remark}[theorem]{Remark}
 \newenvironment{proof}{{\it Proof}. }{\hfill\END\\[0.5ex]}

 \newcommand{\R}{\mathbb{R}}

\newcommand{\END}{\hfill$\Box$}

\newcounter{algo}
\newenvironment{algobox}[1]{%
  \refstepcounter{algo}\par\medskip\noindent
  \hrule height 0.6pt\smallskip
  \noindent\textbf{Algorithm \thealgo} (#1)\par\smallskip
  \hrule height 0.3pt\medskip\noindent}%
  {\par\medskip\hrule height 0.6pt\par\medskip}

\begin{document}
\title{An all-steps-at-once algorithm for time-domain wave propagation in heterogeneous and unbounded media}
\author{V. Domínguez\thanks{Departamento de Estadística, Informática y
    Matemáticas, Universidad Pública de Navarra, Campus de Tudela,
    31500 Tudela, Spain, and Institute for Advanced Materials and
    Mathematics (INAMAT\textsuperscript{2}), Universidad Pública de Navarra,
    31006 Pamplona, Spain. E-mail: \texttt{victor.dominguez@unavarra.es}}
  \and M. Ganesh\thanks{Department of Applied Mathematics and Statistics,
    Colorado School of Mines, Golden, CO 80401, USA.
    E-mail: \texttt{mganesh@mines.edu}}}
\date{November,  2026}

\maketitle
\begin{abstract}
  
We propose, analyze, and implement a convolution quadrature (CQ) algorithm for transient acoustic waves in an unbounded medium consisting of a heterogeneous domain embedded in free space, where the wave speed may be spatially variable and discontinuous. Multistep CQ, applied to the time-domain solution operator, reduces the evolution problem to a family of Helmholtz problems with damped complex frequencies, which can be solved independently in parallel, all time steps at once. The frequency problems are discretized by an overlapping coupling of finite elements (FEM) on a polygonal/polyhedral domain containing the inhomogeneity, and boundary elements (BEM) on a smooth curve enclosed by that domain, in the region where the medium is homogeneous, which accounts exactly for the unbounded exterior. Neither discretization needs to conform to the geometry of the heterogeneous domain, nor to each other.

For the resulting FEM--BEM coupling we prove bounds for the solution operators of the transformed problems that are explicit in the Laplace parameter and the wave speed. These bounds yield convergence of the time discretization for A-stable multistep schemes. We also prove explicit estimates on the exponential decay of the transformed solutions away from the heterogeneity, which show that the correction provided by the BEM part is small but relevant, and cheap to compute: a few GMRES iterations suffice, and the correction is needed only over a fixed range of weakly damped frequencies, beyond which it becomes negligible and can be dropped. Numerical experiments confirm the predicted convergence rates in time and space and explore iterative (multigrid) methods for the interior FEM problems.
\end{abstract}


\smallskip

 
 \section{Introduction}
Computational algorithms for simulating the propagation of time-dependent acoustic waves in heterogeneous and unbounded media are fundamental to numerous applications~\cite{bergman2018computational}. Three main strategies can be found in the literature: 
finite elements in a truncated region, with the exterior replaced by approximate non-reflecting boundary conditions (NRBCs)~\cite{EngquistMajda1977,Givoli1992,holter2026tutorial} or perfectly matched layers (PMLs)~\cite{Berenger1994,BermudezEtAl2010};
 boundary elements, which represent the unbounded exterior exactly but are restricted to piecewise homogeneous media, since a variable wave speed cannot be accommodated by boundary potentials alone~\cite{BetckeSallesSmigaj2017, GaLo:2023}; and couplings of finite and boundary elements, which combine both, typically with convolution quadrature (CQ) applied only to the boundary integral operators and a time-stepping scheme in the interior, as discussed below.
 
  In this article we consider a time-domain wave propagation model problem in $\mathbb{R}^m$, $m=2,3$, involving a transient scalar wave field that propagates through a heterogeneous domain $\Omega$, where the spatially varying wave speed may be discontinuous, and into the unbounded free-space exterior $\Omega^c=\mathbb{R}^m\setminus\overline{\Omega}$. The inhomogeneity is excited by an incident wave, and the resulting time-dependent wave field can be described, through the Laplace transform, by a family of frequency-domain scattering problems with complex wave numbers. CQ replaces this family by a finite set of such problems and reconstructs the full time history from their solutions, providing a robust and flexible framework for numerical simulation.  Each transformed problem is solved by an overlapping coupling of the finite element method (FEM) and the boundary element method (BEM), whose construction and analysis are the subject of this paper.

Let us mention \cite{LechMond:2015}, where the model problem solved at each complex frequency is exactly the one considered here; the difference lies in its discretization. There, the transmission problem is recast as a volume Lippmann--Schwinger equation (see \cite[Ch.~8]{ColKress:2013}, \cite[Ch.~10]{SaVa:2002}) and solved by trigonometric collocation, so that the unknown is supported throughout the inhomogeneity and the resulting operator is dense. Here, we retain the differential formulation and discretize it by finite elements, using boundary integral operators only to account for the unbounded exterior. Consequently, the analysis in \cite{LechMond:2015} requires the sound speed to be globally smooth, whereas our formulation accommodates discontinuities in the sound speed across material interfaces.

Problems of this kind have also been treated numerically using non-overlapping FEM--BEM coupling \cite{Banjai:2022,BaLuSa:2015,Eber:2022,HaSa:2016,MeRi:2017}. In all of these approaches, the boundary integral equation, discretized in time by CQ, is posed on the boundary of the finite element domain, which serves as both the truncation boundary and the coupling interface. The interior problem is advanced in time either by an explicit leapfrog scheme~\cite{BaLuSa:2015,Eber:2022,Banjai:2022} or by the same A-stable scheme used to generate the CQ. Among the implicit schemes, we mention the trapezoidal rule in \cite{HaSa:2016} and a Runge--Kutta method in \cite{MeRi:2017}, which is devoted to the time-dependent linear Schr\"odinger equation.

The method we propose and analyze here differs from these approaches in its coupling strategy and, in particular, in its geometric flexibility. The finite element domain is chosen as a simple polygonal/polyhedral domain containing the heterogeneous domain, while the boundary element problem is posed on an auxiliary curve/surface  lying inside the finite element domain and surrounding the heterogeneity, in the region where the medium is homogeneous, so that the two discretizations overlap. Consequently, neither the finite element domain nor the boundary element curve/surface needs to conform to, or depend on, the geometry of the heterogeneity. This separation of the discretization geometries provides considerable flexibility in handling complex or irregular heterogeneous domains. Moreover, the auxiliary boundary element curve/surface can be chosen to be smooth, or even analytic, independently of the geometry of the heterogeneity, which is particularly advantageous for high-order BEM.

In contrast, in a non-overlapping coupling, the boundary integral operators are posed on the polygonal/polyhedral boundary of the finite element domain and therefore inherit its corners. The finite element mesh may of course be chosen to resolve the material interfaces, as we do in our numerical experiments, to preserve the accuracy of the finite elements across jumps in the speed of sound; however, this is a property of the mesh and is not a requirement of the coupling strategy; see Section~\ref{sec:numerics} for the alternative. This work can be viewed as a time-domain (TD) counterpart of the frequency-domain framework in~\cite{DoGa:2022,DoGaSa:2020}.

The advantages of our overlapping TD-FEM--BEM approach are several. First, the all-steps-at-once implementation is naturally parallel. This is an algorithmic choice rather than a choice of discretization: the same discrete solution can also be computed by marching in time, in which case one recovers a scheme of the type used in \cite{HaSa:2016}; see Remark~\ref{rem:marching}. In the all-steps-at-once formulation, however, the decoupled frequency-domain problems are independent and uniformly damped, and can therefore be solved fully in parallel.

Second, the exterior is accounted for exactly. Truncating the computational domain by means of NRBCs~\cite{holter2026tutorial} introduces an additional approximation whose accuracy must be controlled independently, whereas the boundary element component provides an exact transparent condition on $\Sigma$, the boundary of the FEM domain.

The price to be paid is a non-local operator. However, this price is seldom significant: unlike in classical acoustic scattering, the complex frequencies generated by the Laplace transform introduce strong damping, so that the transmitted field is negligible on $\Sigma$ for all but a small fraction of the frequencies. The coupling is therefore activated only within a fixed frequency window, whose size we quantify in this article. Even within this window, the coupled FEM--BEM system is solved by GMRES, which converges in very few iterations, typically one or two.

Finally, the interior problems associated with the transformed formulation are well suited to iterative solvers. They are coercive for every Laplace parameter $s$, with positive real part, with a coercivity constant proportional to the relative damping $\Re (s)/|s|$. 
The complex shift therefore improves the smoothing properties of standard relaxation schemes and mitigates the resonance phenomena that would otherwise arise at purely real frequencies; this is the mechanism exploited by shifted-Laplacian preconditioners for the Helmholtz equation \cite{ErOoVu:2006}, where, in contrast to our setting, the shift is introduced artificially. This makes multigrid a viable option, although, as our experiments show, not over the entire range of frequencies produced by the CQ discretization. This observation is particularly important in three dimensions, where sparse direct solvers quickly become prohibitively expensive.

On the theoretical side, this article provides bounds for the solution operators of the transformed problems that are explicit in the Laplace parameter and in the wave speed (Proposition~\ref{prop:estimates}). These bounds place the problem within the framework of \cite{BanjaiSayas2023} and yield convergence of the time discretization for BDF1, BDF2, and the trapezoidal rule, with the order of the underlying multistep method and with constants that grow at most polynomially with the final time (Theorem~\ref{th:main}). We also prove exponential decay estimates for the transformed solutions away from the inhomogeneity (Proposition~\ref{prop:estimates:02}). These estimates quantify the frequency window in which the boundary element correction is needed (Theorem~\ref{theo:the:theo:01}) and show that the coupled operator is a contraction outside this window (Theorem~\ref{theo:the:theo:02}).

The analysis in this article is restricted to the time discretization. The dependence on the spatial discretization of the BEM and FEM components is more delicate, but we believe that such an analysis is attainable in a future work. The damping introduced by the Laplace transform, with $\Re (s)\ge\sigma_0>0$ uniformly along the quadrature contour, keeps all transformed problems coercive and avoids the resonances associated with purely real frequencies. The main difficulty is therefore to track the dependence of the relevant constants on $s$, in the spirit of the estimates in Section~\ref{sec:estimates} and the semidiscrete analyses in \cite{BaLuSa:2015,HaSa:2016,Sayas:2016}.

Let us point out, however, that this damping does not by itself prevent the pollution effect \cite{MelenkSauter2011}, since a sizable portion of the quadrature frequencies has small relative damping $\Re (s)/|s|$. In our experiments, however, the frequencies that carry non-negligible data are resolved by the mesh, while those that are not carry data only at the level of round-off error; see Section~\ref{sec:numerics}. 

This paper is structured as follows. In Section~\ref{sec:FULLCQ}, we introduce the physical model problem under consideration and the mathematical model used for its numerical approximation. We then describe the linear multistep CQ method employed for the time discretization, together with its all-steps-at-once implementation, and present the convergence results in time, relying on the well-established theory for CQ methods \cite{BanjaiSayas2023}; see also \cite{LalienaSayas2009,Lubich1994,Sayas:2016} for multistep CQ based on BDF-type schemes. In Section~\ref{sec:coupling}, we explain how the transformed problems are addressed using our overlapping FEM--BEM framework and analyze the resulting coupled equation. A key point is the decay of the solution in the exterior of the inhomogeneity, which determines the form of the transformed problem and the role of the BEM component. This issue is examined in detail in Section~\ref{sec:estimates}, where we also establish bounds for the continuity constants with respect to the transformed variable---precisely the type of estimates required in the preceding sections. Section~\ref{sec:numerics} contains an extensive set of numerical experiments illustrating the efficiency of the algorithm and providing further support for the theoretical results. Finally, Section~\ref{sec:conclusions} summarizes the main conclusions and outlines directions for future work.

\section{A model problem and time-discretization algorithm}
\label{sec:FULLCQ}
\subsection{Transient acoustic wave problem}
We consider a model problem for transient acoustic wave propagation in
$\mathbb{R}^m$, $m=2,3$, consisting of a heterogeneous domain $\Omega$
embedded in an unbounded homogeneous exterior region. The wave speed
$c\in L^\infty(\mathbb{R}^m)$ is positive and may vary spatially and be
discontinuous within $\Omega$. In the exterior free space
$\Omega^c=\mathbb{R}^m\setminus\overline{\Omega}$, the wave speed is
normalized to the constant $c_{\rm ext}=1$.

The evolution of the {\em total} wave field in $\mathbb{R}^m$ is described by
the scalar field $\vartheta=\vartheta(\mathbf{x},t)$, which satisfies
the classical time-domain wave equation:
\[
  c^{-2}(\mathbf{x})\,\vartheta_{tt}(\mathbf{x},t) = \Delta \vartheta(\mathbf{x},t),
  \qquad \mathbf{x}\in\mathbb{R}^m,\ t>0.
\]
The total field is induced by a smooth incident field $\vartheta^{\mathrm{inc}}$, which originates in the exterior region $\Omega^c$, 
and with support that does not intersect $\Omega$ for $t\le 0$. The incident field satisfies the wave equation
\[
  \vartheta^{\mathrm{inc}}_{tt}(\mathbf{x},t) = \Delta\vartheta^{\mathrm{inc}}(\mathbf{x},t),
  \qquad \mathbf{x}\in\mathbb{R}^m.
\]

The interaction of this wave with the inhomogeneity generates a scattered and propagating field
$u = \vartheta - \vartheta^{\mathrm{inc}}$, which satisfies homogeneous initial
data and the evolution equation
\begin{equation}\label{eq:theproblem}
 \left|
 \begin{array}{l}
  u\in H^2([0,T], L^2(\mathbb{R}^m)),  \\
 \dfrac{1}{c^2} \ddot{u}(t)  =     \Delta u(t) +f(t) ,\quad t\in [0,T],\\
  u(0) = 0, \quad \dot{u} (0) = 0,
 \end{array}
 \right.
\end{equation}
where
\[
  f(\mathbf{x},t) := c_{{\rm R}}({\mathbf{x}})\,\vartheta^{\mathrm{inc}}_{tt}(\mathbf{x},t),\quad
  c_{{\rm R}}(\mathbf{x}) := 1 - c^{-2}(\mathbf{x})
\]
is supported inside $\Omega$ since $c(\mathbf{x})=1$ in the free space (i.e.  for $\mathbf{x}\in\Omega^{\rm c}$). The function $c_{{\rm R}}$ is the contrast function and the regularity in space of the solution $u$ depends on the regularity of $c_{{\rm R}}$  as a function in $\mathbb{R}^m$.

\subsection{Convolution Quadrature Method}

Let us introduce, for now only in a formal way, the Laplace transform:
\[
 F(s) := (\mathcal{L} f)(s) := \int_0^{\infty} f(t)\, e^{-s t}\, \mathrm{d}t,
 \quad s \in \mathbb{C}_+ := \{s \in \mathbb{C} : \Re s > 0\}.
\]
Then $U = \mathcal{L}u$ solves the following problem:
\begin{subequations}\label{eq:theproblem:02}
\begin{equation}\label{eq:theproblem:02a}
 \left|
 \begin{array}{l}
 U(s) \in H^1(\mathbb{R}^m),\\[4pt]
 \dfrac{s^2}{c^2}\, U(s) - \Delta U(s) = F(s).
 \end{array}
 \right.
\end{equation}
The analysis of problem~\eqref{eq:theproblem:02}, together with the derivation of suitable bounds and its numerical resolution, will occupy most of this work. For now, we simply state that the problem is well posed and gives rise to a solution operator $\mathrm{K}(s)$ such that
\begin{equation}\label{eq:theproblem:02b}
U = \mathrm{K}(s)F,
\end{equation}
\end{subequations}
where $\mathrm{K}(s): L^2(\Omega) \to H^1(\mathbb{R}^m)$; its range is in fact contained in $H^2(\mathbb{R}^m)$. Moreover, $\mathrm{K}(s)$ is analytic for $s \in \mathbb{C}_+$.

Fix $t>0$ and set $\widetilde f(\tau):=f(\tau)\,\chi_{[0,t]}(\tau)$. By causality of the wave equation
(the solution of \eqref{eq:theproblem} at time $t$ only depends on $f|_{[0,t]}$),
$u(t)$ is also the value at $t$ of the solution with datum $\widetilde f$, so that,
formally,
\begin{equation}\label{eq:u}
 \begin{aligned}
 u(t)\ &=\ \frac{1}{2\pi{\rm i}}\int_{\mu+{\rm i}\mathbb{R}} e^{\zeta t}\,\mathrm{K}(\zeta)\widetilde F(\zeta)\,{\rm d}\zeta
 \ =\ \frac{1}{2\pi{\rm i}}\int_{\mu+{\rm i}\mathbb{R}} \mathrm{K}(\zeta)
   \Big[\int_0^\infty \widetilde f(\tau)\,e^{\zeta(t-\tau)}\,{\rm d}\tau\Big]{\rm d}\zeta\\
 &=\ \frac{1}{2\pi{\rm i}}\int_{\mu+{\rm i}\mathbb{R}} \mathrm{K}(\zeta)
   \underbrace{\Big[\int_0^{t} f(\tau)\,e^{\zeta(t-\tau)}\,{\rm d}\tau\Big]}_{y(t;\zeta)}\,{\rm d}\zeta ,
 \qquad \mu>0 .
 \end{aligned}
\end{equation}
(The argument can be made rigorous in $L^2(\mathbb{R}^m)$: by Proposition~\ref{prop:estimates} below, and since $\widetilde F(\zeta)=O(|\zeta|^{-1})$, the integral exists in the classical sense.)
Here,  $y(t;\zeta)$, with $\zeta\in\mathbb{C}_+$, is the solution of
\begin{equation}\label{eq:ode}
 y' - \zeta y =f, \quad y (0)=0.
\end{equation}
Multistep CQ methods are obtained by discretizing \eqref{eq:ode}: define $y^\kappa(t;\zeta)$ as the solution of the multistep recursion 
\[
\frac{1}{\kappa} \sum_{j = 0}^r \alpha_{j} y^\kappa (t-j\kappa;\zeta) =
 \sum_{j = 0}^r \beta_{j} f (t-j\kappa)+ \zeta \sum_{j = 0}^r \beta_{j} y^\kappa(t-j\kappa;\zeta), \qquad y^\kappa(t;\zeta) =0,\ \text{for}\ t\le 0.
\] 
Applying the Laplace transform, with $
Y^\kappa(s;\zeta):= {\cal L}(y^\kappa(\cdot; \zeta))(s)$ we have
\[
 Y^\kappa({\zetabis;\zeta})=\frac{1}{\dfrac{\delta(e^{-\kappa \zetabis})}{\kappa}-\zeta} F(\zetabis)
\]
where
\[
 \delta(z) = \frac{\alpha_{r}z^r+\alpha_{r-1}z^{r-1}+\ldots+\alpha_0}{\beta_{r}z^r+\beta_{r-1}z^{r-1}+\ldots+\beta_0}
\]
is the generating function of the linear multistep method, also known in the CQ literature as the discrete differentiation symbol.

The CQ approximation is then defined by replacing $y$ with $y^\kappa$ in \eqref{eq:u}: with $\sigma=\Re s>0$,
\begin{equation}\label{eq:u2}
\begin{aligned}
u^\kappa (t) &\ =\  \frac{1}{2\pi \mathrm{i}}\int_{\mu +\mathrm{i}\mathbb{R}} \mathrm{K}(\zeta)y^\kappa(t;\zeta)\,{\rm d}\zeta =
\frac{1}{2\pi \mathrm{i}}\int_{\mu +\mathrm{i}\mathbb{R}} \mathrm{K}(\zeta) {\cal L}^{-1}( Y^\kappa(\cdot;\zeta) )(t)\,{\rm d}\zeta\\
&\ =\  \frac{1}{2\pi \mathrm{i}}\int_{\mu+\mathrm{i}\mathbb{R}}\mathrm{K}(\zeta) \left[\frac{1}{2\pi \mathrm{i}}\int_{ {\sigma} +\mathrm{i}\mathbb{R}} e^{t\,\zetabis}\frac{1}{\dfrac{\delta(e^{-\kappa  \zetabis})}{\kappa}-\zeta} F( {\zetabis})\,{\rm d} {\zetabis} \right]\,{\rm d}\zeta \\
&\ =\  \frac{1}{2\pi \mathrm{i}} \int_{ {\sigma} +\mathrm{i}\mathbb{R}} e^{t\,\zetabis}
\left[  \frac{1}{2\pi \mathrm{i}}\int_{ \mu +\mathrm{i}\mathbb{R}} \mathrm{K}(\zeta) \frac{1}{\dfrac{\delta(e^{-\kappa \zetabis})}{\kappa}-\zeta} \,{\rm d}\zeta  \right]F(\zetabis)\,{\rm d}\zetabis\\
&\ =\  \frac{1}{2\pi \mathrm{i}} \int_{ {\sigma} +\mathrm{i}\mathbb{R}}e^{t\,\zetabis}
 \mathrm{K}\left(\dfrac{\delta(e^{-\kappa \zetabis})}{\kappa} \right) F(\zetabis)\,{\rm d}\zetabis.
\end{aligned}
\end{equation}
Observe that the Cauchy integral formula has been applied above in the last step. This requires, first, that $\mathrm{K}\big(\delta(e^{-\kappa s})/\kappa\big)$ be analytic for $\Re s>0$, i.e.\ $\Re\delta(z)>0$ for $|z|<1$, which imposes a severe restriction on the multistep method (A--stability, see Remark~\ref{rem:newmark}); and, second, that the pole $\zeta=\delta(e^{-\kappa s})/\kappa$ lie to the right of the contour $\mu+{\rm i}\mathbb{R}$ for every $s\in\sigma+{\rm i}\mathbb{R}$, that is, $\mu<\inf_{s\in\sigma+{\rm i}\mathbb{R}}\Re\delta(e^{-\kappa s})/\kappa$. For the three methods considered below this infimum equals $\delta(e^{-\kappa\sigma})/\kappa<\sigma$, so that $\mu$ has to be taken strictly smaller than the abscissa $\sigma$ of the outer integral.

The same assumption makes $z\mapsto\mathrm{K}(\delta(z)/\kappa)$ analytic in the unit disc, so that it admits the power series expansion
\begin{equation}\label{eq:analyticity:K}
 \mathrm{K}\left(\frac1{\kappa}\delta(z)\right) = \sum_{j=0}^\infty \mathrm{W}_j(\kappa) z^j,
\end{equation}
which converges absolutely on $|z|=e^{-\kappa\sigma}<1$. Inserting it in \eqref{eq:u2} and exchanging sum and integral,
\begin{equation}\label{eq:u_kappa}
  u^\kappa (t) \ =\ \sum_{j=0}^\infty  \mathrm{W}_j(\kappa)\left[  \frac{1}{2\pi \mathrm{i}} \int_{ {\sigma} +\mathrm{i}\mathbb{R}}e^{s(t-j\kappa)} F(s)\,{\rm d}s\right]=\sum_{j=0}^\infty  \mathrm{W}_j(\kappa) f(t-j\kappa),
\end{equation}
where $f$ is extended by zero to $t<0$, so that by causality the sum is finite.

\begin{remark}\label{rem:newmark}
Convolution Quadrature methods have been an active topic of research since the seminal papers of Lubich \cite{Lubich1988CQI,Lubich1988CQII}, originally developed for the solution of both parabolic and hyperbolic problems. The method hinges on the expansion \eqref{eq:analyticity:K} and the calculation, or approximation, of the weights $\{\mathrm{W}_j(\kappa)\}_{j=0}^\infty$. For such an expansion to exist the linear multistep method has to be A--stable, as mentioned above, and by the second Dahlquist barrier no A--stable linear multistep method has order higher than two. Higher order in time can be achieved with Runge--Kutta based CQ, a case that we will not consider here.

Let us also point out that convolution quadrature is compatible with composition, $\mathrm{K}^{-1}(\partial_t^{\kappa})\mathrm{K}(\partial_t^{\kappa})=I$, so that $u^\kappa$ can equivalently be characterized as the solution of
\[
\mathrm{K}^{-1}(\delta(e^{-\kappa s})/\kappa)\,U^\kappa(s)=F(s)
\]
i.e.\ as the result of applying the multistep method directly to the wave equation. For the trapezoidal rule, since
\[
 \mathrm{K}^{-1}\!\left(\frac{\delta(z)}{\kappa}\right)
  = \frac{4}{1+2z+z^2}\left[\frac{1-2z+z^2}{\kappa^2c^2}-\frac{1+2z+z^2}{4}\,\Delta\right],
\]
this is the recursion
\begin{multline*}
 \frac{1}{c^2}\,\frac{u^\kappa(t)-2u^\kappa(t-\kappa)+u^\kappa(t-2\kappa)}{\kappa^2}
 -\frac14\Delta\big(u^\kappa(t)+2u^\kappa(t-\kappa)+u^\kappa(t-2\kappa)\big)\\
 =\frac14\big(f(t)+2f(t-\kappa)+f(t-2\kappa)\big),
\end{multline*}
that is, the Newmark method with $\gamma=1/2$, $\beta=1/4$, unconditionally stable and of
order two. The all--steps--at--once implementation of Subsection~\ref{subsec:cq:allsteps} and
this marching form compute the same discrete solution; the choice between them is one of
algorithm, not of discretization.
\end{remark}

\subsection{Convergence  of linear multistep CQ methods}

We restrict our attention to the BDF1, BDF2 and trapezoidal--rule CQ methods; see Remark~\ref{rem:main} below for other linear multistep methods. Following C.~Lubich \cite{Lubich1988CQI}, and as in \cite{BanjaiSayas2023,Sayas:2016}, the CQ solution is written as
\[
 \mathrm{K}(\partial_t^\kappa)f(t) := \frac{1}{2\pi{\rm i}}\int_{\sigma +{\rm i}\mathbb{R}} e^{s t} \mathrm{K}(s^\kappa ) F(s)\,{\rm d}s = U^\kappa(t),
 \qquad
 s^\kappa := \frac{\delta(e^{-s \kappa})}{\kappa},
\]
as an approximation of the exact solution
\[
  \mathrm{K}(\partial_t)f(t) := \frac{1}{2\pi{\rm i}}\int_{\sigma +{\rm i}\mathbb{R}} e^{s t} \mathrm{K}(s) F(s)\,{\rm d}s = U(t).
\]
The use of $\partial_t$ as the \emph{argument} of $\mathrm{K}$ emphasizes the role of the Laplace transform (since ${\cal L}(\partial_t)=s$) in both the definition of the method and the convergence analysis, even though the Laplace transform is never actually applied to the right--hand side of the wave equation.

We will consider the linear multistep CQ method for three representative cases: for BDF1 (implicit Euler), with
\[
 \delta_{\rm E}(z) = 1-z,
\]
of order $p=1$; for BDF2,
\[
 \delta_{\rm BDF2}(z) = \tfrac12 z^2 - 2z + \tfrac32;
\]
and for the trapezoidal rule,
\[
 \delta_{\rm TR}(z) = \frac{2(1-z)}{1+z},
\]
of order $p=2$. 

The proof of the following theorem is given in  Section~\ref{proofs}.
\begin{theorem}\label{th:main}
Let $\delta$ be the generating function of BDF1 ($p=1$), BDF2 or the trapezoidal rule ($p=2$).
Let $r\in[1,2]$, let $\widetilde{p}\in\mathbb{N}$ with $\widetilde{p}>r+p+1$, and let
$f\in C^{\widetilde{p}}([0,T];L^2(\Omega))$ satisfy $f(0)=\cdots=f^{(\widetilde{p}-1)}(0)=0$.
Then there exist $\kappa_0>0$ and $C>0$, depending only on $p$, $r$ and $\delta$, such that
for all $\kappa\in(0,\kappa_0)$, $\sigma\in(0,1)$ and $t\in(0,T]$
\[
 \|\mathrm{K}(\partial_t^\kappa)f(t)-\mathrm{K}(\partial_t)f(t)\|_{H^r(\mathbb{R}^m)}
 \le C\, c_{\rm max}^2\, c_{\rm min}^{-2(r-1)}\, \sigma^{p-1-\widetilde{p}}\, e^{\sigma t}\,\kappa^p
 \int_0^t e^{-\sigma\tau}\,\|f^{(\widetilde{p})}(\tau)\|_{L^2(\Omega)}\,{\rm d}\tau .
\]
\end{theorem}

\begin{remark}\label{rem:main}
Convergence in $H^r(\mathbb{R}^m)$ for $r\in[0,1)$ needs no separate statement, since the case
$r=1$ already provides it. For $r\in[1,2)$ it suffices to take $\widetilde{p}=p+3$, and the
constant in Theorem~\ref{th:main} is then proportional to $\sigma^{-4}$; for $r=2$ one needs
$\widetilde{p}=p+4$, and the constant is proportional to $\sigma^{-5}$. The dependence on
$\sigma$ is the same for the three schemes. The cases $r=1+\varepsilon$ for $m=2$ and
$r=3/2+\varepsilon$ for $m=3$ yield, via the Sobolev embedding theorem, estimates in the
uniform norm with $\widetilde{p}=p+3$. Choosing $\sigma=1/t$ for $t>1$ (and, say,
$\sigma=1/2$ for $t\le1$), the estimate takes the form
\[
 \|\mathrm{K}(\partial_t^\kappa)f(t)-\mathrm{K}(\partial_t)f(t)\|_{H^r(\mathbb{R}^m)}
 \le C\,\max\{1,t\}^{\widetilde{p}+1-p}\,\kappa^p\int_0^t\|f^{(\widetilde{p})}(\tau)\|_{L^2(\Omega)}\,{\rm d}\tau,
\]
so that the constant grows at most like $t^4$ ($r<2$) or $t^5$ ($r=2$) with the final time,
with no exponential factor.

Theorem~\ref{th:main} is stated for the three multistep methods used in this paper, but nothing in its proof is specific to them: the derivation of Section~\ref{sec:FULLCQ} only uses that $\Re\delta(z)>0$ for $|z|<1$, and the convergence results of \cite{BanjaiSayas2023} on which the proof rests apply to any A--stable linear multistep method whose generating function satisfies the hypotheses stated there, in which case the theorem holds with $p$, the order of the method. Among the A--stable methods of order two, BDF2 and the trapezoidal rule are the two extreme choices, maximal damping of the high frequencies and none at all, and they lead to markedly different distributions of the quadrature frequencies, as discussed in Section~\ref{sec:numerics}. Intermediate schemes, such as the truncated trapezoidal rules mentioned in \cite{Banjai:2022}, which maintain a polynomial generating function and thus exhibit bounded quadrature frequencies similar to BDF2, can also be considered. Their error constants are close to those of the trapezoidal rule, but still satisfy the criteria for applying the convergence result from Theorem~\ref{th:main}.
\end{remark}

\subsection{All--steps--at--once implementation}
\label{subsec:cq:allsteps}
The discretization in time thus depends on computing  sufficiently accurate approximations  of
 \begin{equation}\label{eq:analyticity:K:2}
  \begin{aligned}
 \mathrm{W}_j(\kappa)\ = \ & \frac{1}{j!}
 \left.\frac{{\rm d}^j}{{\rm d}z^j} \mathrm{K}\left(\frac1{\kappa}\delta(z)\right) \right|_{z=0}
 = \frac{1}{2\pi{\rm i}}
  \oint_{|z|=r} z^{-j-1} \mathrm{K}\left(\frac1{\kappa}\delta(z)\right)\,{\rm d}z\\
  = \ &
\left[\frac{1}{r^j}\int_{0}^1 \exp(-2\pi {\rm i}j   x)   \mathrm{K}\left(\frac1{\kappa}\delta(r\exp(2\pi{\rm i}   x) )\right)\,{\rm d}x\right]
  \end{aligned}
 \end{equation}
 (by the Cauchy integral formula on the circle $|z|=r$, $0<r<1$, which lies inside the disc of analyticity of $z\mapsto\mathrm{K}(\delta(z)/\kappa)$).

Take a positive integer $n_{\rm Q}$ and set
\[
 e_{\ell,n_{\rm Q}}:= \exp\left(\frac{2\pi{\rm i}\ell}{n_{\rm Q}+1}  \right), \quad
 z_{\ell,n_{\rm Q}}^{\kappa} = \frac{1}{\kappa}\delta\left(r e_{\ell,n_{\rm Q}}
 \right),
\]
and define the following approximations for $\mathrm{W}_j(\kappa)$:
\begin{equation}
\label{eq:wjncqkappa}
 \mathrm{W}^{n_{\rm Q}}_j(\kappa) g := \frac{1}{n_{\rm Q}+1 }  r^{-j} \sum_{\ell = 0}^{n_{\rm Q} }
 e_{-j\ell,n_{\rm Q}}
 \mathrm{K}(z_{\ell,n_{\rm Q}}^\kappa)  g .
\end{equation}
That is, the weights are approximated by the trapezoidal rule on a uniform grid with $n_{\rm Q}+1$ points, which leads to the following approximation of $u^\kappa$ (cf. \eqref{eq:u_kappa})
 \[
\begin{aligned}
u_{n_{\rm Q}}^\kappa (t) &\ :=\ \sum_{j =0}^{\lceil t/\kappa \rceil }  \mathrm{W}^{n_{\rm Q}}_j(\kappa)  f(t-j\kappa).
 \end{aligned}
\]
Let $n_{\rm CQ}$ denote from now on the number of time steps, so that the final time is $T = n_{\rm CQ}\,\kappa$. We next apply \eqref{eq:wjncqkappa} to approximate $\mathrm{W}^{n_{\rm Q}}_{n-j}(\kappa)$, even for $j>n$ (notice that $\mathrm{W}_{n}=0$  for $n<0$).
We then define
\[
 u(n\kappa)   \approx   u^{\kappa,n_{\rm Q}} (n\kappa ) \ :=\
\sum_{j = 0}^{n_{\rm CQ} } \frac{1}{ n_{\rm Q} +1} r^{j-n}
\bigg[
\sum_{\ell = 0}^{n_{\rm Q} }
 e_{(j-n)\ell,n_{\rm Q}}
 \mathrm{K}(z_{\ell,n_{\rm Q}}^\kappa)  f(j\kappa)\bigg].
\]
Interchanging the order of summation,
\begin{eqnarray*}
u^{\kappa,n_{\rm Q}} (n\kappa) \  =\ r^{-n} \sum_{\ell = 0}^{n_{\rm Q} }
 e_{-\ell n,n_{\rm Q} }  \mathrm{K}(z_{\ell,n_{\rm Q}}^\kappa)   \bigg[\frac{1}{ n_{\rm Q} +1}
\sum_{j = 0}^{n_{\rm  CQ} } r^{j} e_{j\ell,n_{\rm Q}}
 f(j\kappa)
\bigg]
\end{eqnarray*}
Although the approximation
\eqref{eq:analyticity:K:2}--\eqref{eq:wjncqkappa} was already considered in the seminal papers of C.~Lubich
\cite{Lubich1988CQI,Lubich1988CQII}, this implementation--commonly referred to in the literature as the all-steps-at-once approach--was introduced in \cite{BaSa:2008}. 

The approximation \eqref{eq:wjncqkappa} requires $n_{\rm Q}+1$ evaluations of $\mathrm{K}$, half of them by symmetry (see \eqref{eq:sym:01} below), and does not affect the accuracy of the method for the choice of $r$ indicated below \cite{BaSa:2008}. Its interest is that it allows the whole computation to be organized around the Fast Fourier Transform, which is the so-called all--steps--at--once implementation:

\begin{algobox}{all--steps--at--once evaluation of $u^\kappa$ with $n_{\rm Q}=n_{\rm CQ}$}\label{alg:allsteps}%
\textbf{Input:} $f(j\kappa)$, $j=0,\dots,n_{\rm Q}$; radius $0<r<1$; generating function $\delta$.\\
\textbf{Output:} $\widetilde u^\kappa_r(n\kappa)$, $n=1,\dots,n_{\rm Q}$, approximating $u^\kappa(n\kappa)$ in \eqref{eq:u_kappa}.
\begin{subequations}\label{eq:allsteps}
\begin{enumerate}
 \item \textbf{Scale the data:}
 \begin{equation}\label{eq:frkappa}
  {\bf f}^{\kappa}_{r} := (r^{j} f(j\kappa))_{j=0}^{n_{\rm Q}} \in \left[L^2(\Omega)\right]^{n_{\rm Q}+1}.
 \end{equation}
 \item \textbf{Transform (inverse FFT):}\label{step:2}
 \begin{equation}\label{eq:def:F}
  F^\kappa_{n,r} := \frac{1}{n_{\rm Q}+1}\sum_{j=0}^{n_{\rm Q}} e_{jn,n_{\rm Q}}\, f^\kappa_{j,r},
  \qquad n=0,\dots,n_{\rm Q}.
 \end{equation}
 \item \textbf{Solve the transformed problems, in parallel:}\label{step:3}
 \begin{equation}\label{eq:Ynk}
  Y^\kappa_{\ell,r} := \mathrm{K}(z_{\ell,n_{\rm Q}}^\kappa)\, F^\kappa_{\ell,r},
  \qquad \ell=0,\dots,n_{\rm Q},
 \end{equation}
 by Algorithm~\ref{alg:coupling} in  Section~\ref{sec:coupling}; due to \eqref{eq:sym:01} only $\lfloor n_{\rm Q}/2\rfloor+1$ of them are needed.
 \item \textbf{Transform back (FFT) and rescale:}\label{step:4}
 \begin{equation}\label{eq:urk}
  \widetilde{u}^\kappa_{r}(n\kappa) := r^{-n}\sum_{\ell=0}^{n_{\rm Q}} Y^\kappa_{\ell,r}\, e_{-\ell n,n_{\rm Q}},
  \qquad n=1,\dots,n_{\rm Q}.
 \end{equation}
\end{enumerate}
\end{subequations}
\end{algobox}

Steps~\ref{step:2} and \ref{step:4} are FFTs of length $n_{\rm Q}+1$; the cost is dominated by Step~\ref{step:3}, which is obviously parallelizable, one of the strong points of the algorithm.

 Note  that
\[
  z_{n_{\rm Q}+1-\ell,n_{\rm Q}}^\kappa = \overline{z_{\ell,n_{\rm Q}}^\kappa}, \quad \ell = 1,\ldots, n_{\rm Q},
\]
since $e_{n_{\rm Q}+1-\ell,n_{\rm Q}} = \overline{e_{\ell,n_{\rm Q}}}$ and $\delta$ has real coefficients. This implies that, for real-valued data~$f$, as in our case,
\begin{equation}\label{eq:sym:01}
  Y^\kappa_{\ell,r} = \overline{Y^\kappa_{n_{\rm Q}+1-\ell,r}}, \quad \ell = 1,\ldots, n_{\rm Q},
\end{equation}
so that the computations in Step~\ref{step:3}, \eqref{eq:Ynk}, the most expensive part, can be halved.

Two parameters remain to be chosen. A typical choice is $n_{\rm Q}=n_{\rm CQ}$, which we adopt in the experiments of Section~\ref{sec:numerics}; we refer to \cite{BetckeSallesSmigaj2017} for an analysis of the dependence on $n_{\rm Q}$ and for the advantages of overresolving, $n_{\rm Q}>n_{\rm CQ}$. As for the radius, the aliasing error of \eqref{eq:wjncqkappa} is $O(r^{\,n_{\rm Q}+1})$ while the scaling by $r^{-n}$ in \eqref{eq:urk} amplifies round--off errors by up to $r^{-n_{\rm Q}}$; balancing the two gives the standard choice $r={\tt eps}^{1/(2n_{\rm Q})}$, with ${\tt eps}$ the machine precision, for which both errors are $O({\tt eps}^{1/2})$. Better trade--offs are possible by taking $n_{\rm Q}>n_{\rm CQ}$: with $n_{\rm Q}=2n_{\rm CQ}$ and $r={\tt eps}^{1/(3n_{\rm CQ}+1)}$ the error becomes $O({\tt eps}^{2/3})$, and so on; see \cite[Sect.~3.4]{BanjaiSayas2023} and \cite{BetckeSallesSmigaj2017}. Taking $r=1$, which would avoid the round--off amplification, is not an option here: for purely imaginary $s$ problem \eqref{eq:theproblem:02} becomes the heterogeneous Helmholtz equation, which has no solution in $H^1(\mathbb{R}^m)$, and indeed the bound \eqref{eq:01:prop:estimates} degenerates as $\Re s\to0$; the quadrature frequencies must therefore be kept at a positive distance from the imaginary axis.

\section{A FEM--BEM spatial-discretization algorithm}
\label{sec:coupling}

\subsection{Overlapping coupling algorithm for evaluating $\mathrm{K}(s)$}

The method relies on the calculation of \eqref{eq:Ynk}. So we are led to consider the (numerical) solution of the problem
\begin{equation}\label{eq:theproblem:BEMFEM:01}
 U = \mathrm{K}(s) F,\quad \text{where}\quad
 \left|
 \begin{array}{l}
 U \in H^1(\mathbb{R}^m),\\
 (s/c)^2  U  -\Delta U   =F.
 \end{array}
 \right.
\end{equation}
Here we are just assuming that $F\in L^2(\mathbb{R}^m)$, but in our application $\supp F\subset \overline{\Omega}$, a fact that we emphasize by writing $F\in L^2(\Omega)$. We also recall the notation
$\sigma = \Re s>0$. Although the problem is naturally formulated for finding $U$ in $H^1(\R^m)$, we will see later that the solution is smoother, and in fact $U\in {H}^2(\mathbb{R}^m)$.

\begin{figure}
 \[
 \begin{tikzpicture}[scale=1.2]

  \draw[thick,fill=gray!10] (-3,-2) rectangle (3,2);
  \node at (3.25,-0.8) {$\Sigma$};

  \draw[thick, rotate=45,fill=gray!0] (0,0) ellipse (1.25cm and 1.75cm);
  \node at ( 1.85, 0.25) {$\Gamma$};

  \draw[thick, fill=gray!30] (0:.8) -- (72:.8) -- (144:.8) -- (216:.8) -- (288:.8) -- cycle;
  \node at (0,0) {$\Omega$};

 \node at ( 0.8,-0.8 ) {$\Omega_1$};
\node at (2.1,-1.5) {$\Omega_2$};
\end{tikzpicture}
\]
\caption{\label{fig:01}Domains for the space solution algorithm: $\Omega$ (dark gray), with $\overline{\Omega}$ being, or containing, the support of $c_{\rm R}=1-c^{-2}$, the contrast function; $\Omega_1$ (white), the interior of the smooth curve $\Gamma$ (so that $\Omega_1^{\rm c}$ is the exterior); and $\Omega_2$, the interior of the polygonal curve $\Sigma$. The overlapping domain  $\Omega_{12}=\Omega_1^{\rm c}\cap \Omega_2$ (light gray) is also depicted in the figure.}
\end{figure}
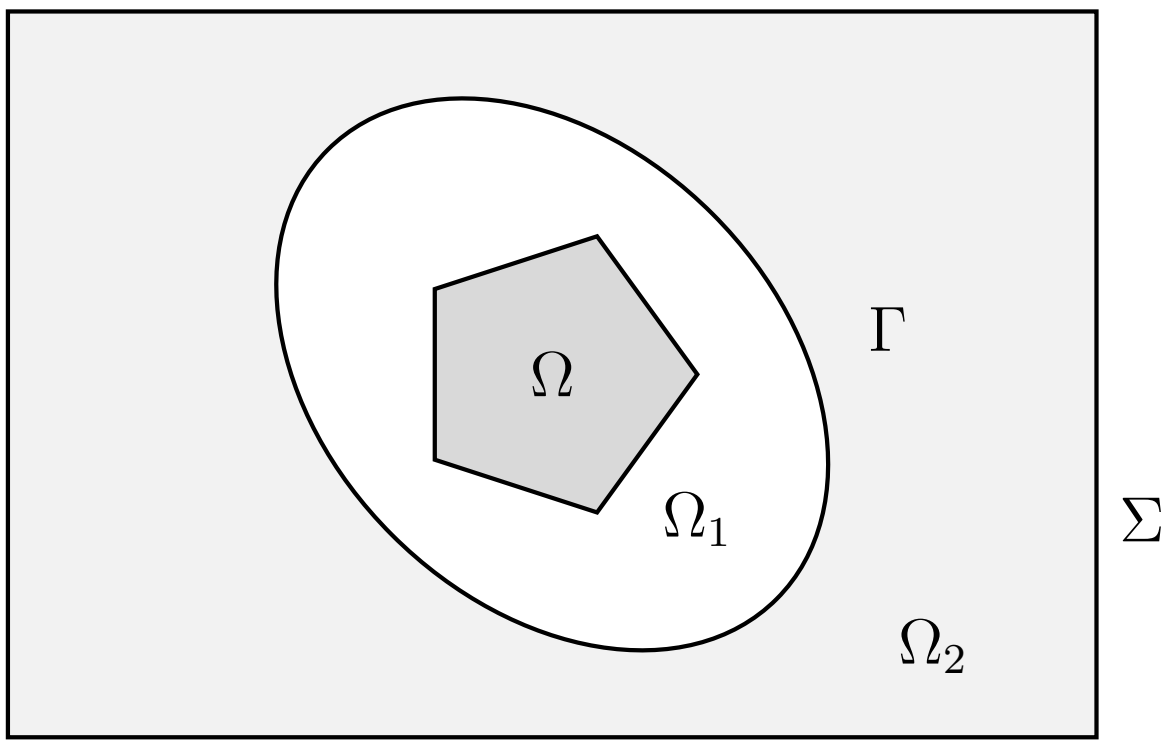

The numerical solution of this problem is obtained by means of an overlapping FEM--BEM approach. We first introduce the continuous framework on which the discrete procedures will be applied. Let $\Omega_1$ be a smooth domain and $\Omega_2$ a polygonal one, 
such that $\overline{\Omega}\subset\Omega_1$ and $\overline{\Omega}_1\subset\Omega_2$. Their boundaries are denoted by $\Gamma$ (for $\Omega_1$) and $\Sigma$ (for $\Omega_2$), respectively (see Figure~\ref{fig:01} for a schematic representation of the domains involved).

Let $\gamma_\Sigma$ and  $\gamma_\Gamma$ be trace  operators  that evaluate functions on the artificial boundaries $\Sigma$ and $\Gamma$, respectively.  Let $W_F^{\rm tr}$ be the solution of
\begin{equation}\label{eq:2:2}
  \left|\begin{array}{ll}
            (s/c)^2   W_F^{\rm tr}- \Delta W_F^{\rm tr}    = F,\quad \text{in $\Omega_2$},\\
            \gamma_\Sigma W_F^{\rm tr} = 0.
        \end{array}
        \right.
 \end{equation}
We introduce two  unknowns
\[
  W^{\rm tr} := U|_{\Omega_2}-W_F^{\rm tr} \in H^1(\Omega_2), \quad
  W^{\rm c}  := U|_{\Omega_1^{\rm c}} \in H^1(\Omega_1^{\rm c}).
\]
$W^{\rm c}$ is analytic in $\overline{\Omega}^{\rm c}_1$, since $U$ solves the Helmholtz equation $s^2U-\Delta U=0$ in the open set $\Omega^{\rm c}\supset\overline{\Omega}_1^{\rm c}$, and $W^{\rm tr}\in H^2(\Omega_2)$ whenever $\Omega_2$ is convex, see Remark~\ref{remark:sec:4}.

We take
 \[
  g_\Sigma= \gamma_\Sigma U   = \gamma_{\Sigma }W^{\rm tr} ,\quad
  g_\Gamma= \gamma_\Gamma U= \gamma_\Gamma W^{\rm c}
 \]
as new unknowns (for the first identity, note that $\gamma_\Sigma W_F^{\rm tr}=0$). Then
\begin{enumerate}

  \item Since, with $\Omega_1^{\rm c}=\mathbb{R}^m\setminus \overline{\Omega}_1$,  $c|_{\Omega_1^{\rm c}}\equiv 1$, $W^{\rm c}$ solves
  \begin{equation}\label{eq:2.3}
  \left|
  \begin{array}{l}
  W^{\rm c}\in H^1(\Omega_1^{\rm c}),\ \gamma_{\Gamma} W^{\rm c} = g_\Gamma,  \\
  s^2 W^{\rm c} - \Delta W^{\rm c} = 0,\quad \text{in }\Omega_1^{\rm c}.
  \end{array}
  \right.
  \end{equation}
  In other words, $W^{\rm c}$ can be computed using BEM schemes 
  provided the Dirichlet datum $g_\Gamma=\gamma_\Gamma W^{\rm c}=\gamma_\Gamma U$ is known, since problem~\eqref{eq:2.3} is uniquely solvable (this holds even for $\sigma = \Re s = 0$, which is never our case, provided the condition $W^{\rm c}\in H^1(\Omega_1^{\rm c})$ is replaced by a radiation condition at~$\infty$).

  \item For $W^{\rm tr}$, we instead have
  \begin{equation}\label{eq:2.4}
  \left|
  \begin{array}{ll}
    W^{\rm tr}\in H^1(\Omega_2),\ \gamma_\Sigma W^{\rm tr} = g_\Sigma,\\
    (s/c)^2\, W^{\rm tr} - \Delta W^{\rm tr} = 0,\quad \text{in }\Omega_2,
  \end{array}
  \right.
  \end{equation}
  which is well posed since $\sigma = \Re s > 0$. This problem, as well as problem~\eqref{eq:2:2}, is better suited to approximation by FEM schemes, because $\Omega_2$ is bounded and may contain inhomogeneities arising from the wave-speed field~$c$.
\end{enumerate}

Hence we can introduce
\[
\mathrm{K}_{\Omega_1^{\rm c}\Gamma}(s): H^{1/2}(\Gamma)\to H^1(\Omega_1^{\rm c}), \quad
\mathrm{K}^{\rm tr}_{\Omega_2\Sigma}(s): H^{1/2}(\Sigma)\to H^1(\Omega_2),  
\]
as the solution operators for problems \eqref{eq:2.3} and \eqref{eq:2.4},  respectively. That is,
\[
\mathrm{K}_{\Omega_1^{\rm c}\Gamma}(s)g_\Gamma := W^{\rm c},\quad
\mathrm{K}^{\rm tr}_{\Omega_2\Sigma}(s) g_\Sigma:= W^{\rm tr}
\]
so that our unknown $U$ can be expressed as
\[
  U|_{\Omega_2}= \mathrm{K}^{\rm tr}_{\Omega_2\Sigma}(s)g_\Sigma+W_{F}^{\rm tr},\quad U|_{\Omega_1^{\rm c}} = \mathrm{K}_{\Omega_1^{\rm c}\Gamma}(s)g_\Gamma,
\]
provided that the Dirichlet data $g_\Sigma$ and $g_\Gamma$ used to feed the operators $\mathrm{K}^{\rm tr}_{\Omega_2\Sigma}(s)$ and $\mathrm{K}_{\Omega_1^{\rm c}\Gamma}(s)$ are known.

The next step is to determine these unknowns $(g_\Sigma,g_\Gamma)$. We notice that it suffices to impose that both representations agree on the overlap,
\[
(W^{\rm tr}+W_{F}^{\rm tr})\big|_{\Omega_{12}}=
W^{\rm c}\big|_{\Omega_{12} }, \qquad \Omega_{12}:=\Omega_{1}^{\rm c}\cap \Omega_2,
\]
since then the two fields glue together to an element of $H^1(\mathbb{R}^m)$ which solves \eqref{eq:theproblem:BEMFEM:01} (recall that $\supp F\subset\overline{\Omega}$), and hence coincides with $U$. In the notation of the operators introduced above, this amounts to requiring
\[
V:=\big(\mathrm{K}^{\rm tr}_{\Omega_2\Sigma}(s) g_\Sigma + W_{F}^{\rm tr}-  \mathrm{K}_{\Omega_1^{\rm c}\Gamma}(s) g_\Gamma\big)\big|_{\Omega_{1}^{\rm c}\cap \Omega_2} = 0.
\]
Now, $V$ solves the homogeneous equation
 \begin{equation}\label{eq:2.5}
          -\Delta V     + s^2    V      = 0,\quad \text{in $\Omega_{12}$}
 \end{equation}
(the support of $F$ does not meet $\Omega_{12}$, and $c\equiv 1$ there) and, since $\gamma_\Sigma \mathrm{K}^{\rm tr}_{\Omega_2\Sigma}(s)$ and $\gamma_\Gamma \mathrm{K}_{\Omega_1^{\rm c}\Gamma}(s)$ are the identity on $H^{1/2}(\Sigma)$ and $H^{1/2}(\Gamma)$ respectively, and $\gamma_\Sigma W_F^{\rm tr}=0$, its traces are
\[
 \gamma_\Sigma V = g_\Sigma - \gamma_\Sigma \mathrm{K}_{\Omega_1^{\rm c}\Gamma}(s)g_\Gamma,
 \qquad
 \gamma_\Gamma V = \gamma_\Gamma \mathrm{K}^{\rm tr}_{\Omega_2\Sigma}(s)g_\Sigma + \gamma_\Gamma W_F^{\rm tr}- g_\Gamma.
\]
Since the Dirichlet problem for \eqref{eq:2.5} in $\Omega_{12}$ is uniquely solvable for $\sigma=\Re s>0$, we have $V= 0$ if and only if both traces vanish, that is, if and only if $(g_\Sigma,g_\Gamma)$ solves
 \begin{equation}\label{eq:19}
  \left|
  \begin{array}{ccccl}
g_\Sigma &-& \gamma_\Sigma \mathrm{K}_{\Omega_1^{\rm c}\Gamma}(s) g_\Gamma &=&0,\\
- \gamma_\Gamma \mathrm{K}^{\rm tr}_{\Omega_2\Sigma}(s)  g_\Sigma &+&  g_\Gamma &=& \gamma_{\Gamma} W_F^{\rm tr}.
  \end{array}
  \right.
 \end{equation}
Hence  $(g_\Sigma,g_\Gamma)$ is the only solution of \eqref{eq:19}: any solution of \eqref{eq:19} reconstructs, through the representation above, the unique solution $U$ of \eqref{eq:theproblem:BEMFEM:01}, and therefore $g_\Sigma = \gamma_\Sigma U$ and $g_\Gamma = \gamma_\Gamma U$.

The solution of \eqref{eq:19} is obtained by a Schur complement approach:
 \[
  g_\Gamma - \gamma_\Gamma \mathrm{K}^{\rm tr}_{\Omega_2\Sigma}(s) \gamma_{\Sigma}\mathrm{K}_{\Omega_1^{\rm c}\Gamma}(s) g_\Gamma = \gamma_\Gamma W_F^{\rm tr}, \quad g_\Sigma = \gamma_\Sigma \mathrm{K}_{\Omega_1^{\rm c}\Gamma}(s) g_\Gamma.
 \]

 Hence the resulting spatial FEM--BEM algorithm is as follows:

\begin{algobox}{Overlapping FEM--BEM evaluation of $\mathrm{K}(s)F$}\label{alg:coupling}
\textbf{Input:} $s\in\mathbb{C}_+$, $F\in L^2(\Omega)$; domains $\Omega\subset\Omega_1\subset\Omega_2$ with boundaries $\Gamma$ and $\Sigma$.\\
\textbf{Output:} $U=\mathrm{K}(s)F$ on $\Omega_2$ and on $\Omega_1^{\rm c}$.
\begin{enumerate}
 \item \textbf{Interior solve (FEM):} find $W_F^{\rm tr}$ with
 \begin{equation}\label{eq:alg02:01}
  \left|\begin{array}{ll}
         (s/c)^2 W_F^{\rm tr}-\Delta W_F^{\rm tr} = F, & \text{in } \Omega_2,\\
         \gamma_\Sigma W_F^{\rm tr} = 0,
        \end{array}
  \right. 
 \end{equation}
 and set
 \begin{equation}\label{eq:fgamma}
  f_\Gamma := \gamma_{\Gamma} W_F^{\rm tr}.
 \end{equation}
 \item \textbf{Coupling (GMRES on $\Gamma$):} solve the second kind equation
 \begin{equation}\label{eq:the:equation}
  g_\Gamma - \gamma_\Gamma \mathrm{K}^{\rm tr}_{\Omega_2\Sigma}(s)\,\gamma_\Sigma\mathrm{K}_{\Omega_1^{\rm c}\Gamma}(s)\, g_\Gamma = f_\Gamma ,
 \end{equation}
 each application of the operator requiring one exterior solve (BEM) and one interior solve (FEM), and set $g_\Sigma := \gamma_\Sigma \mathrm{K}_{\Omega_1^{\rm c}\Gamma}(s)\, g_\Gamma$.
 \item \textbf{Assemble:} with $W^{\rm tr} := \mathrm{K}^{\rm tr}_{\Omega_2\Sigma}(s)\,g_\Sigma$ and $W^{\rm c} := \mathrm{K}_{\Omega_1^{\rm c}\Gamma}(s)\,g_\Gamma$,
 \begin{equation}\label{eq:alg02:02}
  \mathrm{K}(s)F\big|_{\Omega_2} = W^{\rm tr}+W_F^{\rm tr},\qquad
  \mathrm{K}(s)F\big|_{\Omega_1^{\rm c}} = W^{\rm c}.
 \end{equation}
\end{enumerate}
\end{algobox}

Note that Algorithm~\ref{alg:coupling} provides a double representation of $\mathrm{K}(s)F$ on the overlap $\Omega_{12}$, by $W^{\rm tr}+W_F^{\rm tr}$ and by $W^{\rm c}$, which coincide by construction. The coupled equation \eqref{eq:the:equation} is analyzed in the next subsection.

\begin{remark}[Marching form]\label{rem:marching}
The equivalence of Remark~\ref{rem:newmark} extends to the coupled scheme: applying the  multistep method directly to \eqref{eq:alg02:01}--\eqref{eq:alg02:02} yields a time--stepping algorithm in which the differential part becomes a Newmark--type recursion on $\Omega_2$, with one real symmetric positive definite solve per step, while the boundary element part, being nonlocal in time, becomes a convolution quadrature whose history has to be stored and advanced, and the two are linked at each step through the coupled equation \eqref{eq:the:equation}. This is the structure of \cite{HaSa:2016} and, for the Schrödinger equation with Runge--Kutta methods, of \cite{MeRi:2017}, with the overlapping coupling in place of  the non--overlapping one. The history can be handled by the recursive, FFT--based algorithms of \cite{Banjai:2010,HaLuSc:1985}.

We favor the all--steps--at--once variant for two reasons. In the scattering problems considered here the right--hand side is known in advance on the whole interval $[0,T]$, so that early availability of the solution brings no benefit; and the decoupled formulation yields independent, uniformly damped frequency problems, most of which do not even require the coupling (Section~\ref{sec:estimates}), and which can be solved fully in parallel. 
The price is that no approximation of the solution $u$ of the original time-domain problem~\eqref{eq:theproblem} is available until all the transformed problems have been solved, and that the corresponding volume solutions must be stored.
\end{remark}

\subsection{Analysis of the coupled equation}

The algorithm of the previous subsection rests on two facts: that the boundary element correction $U-W^{\rm tr}_F$ is small for strongly damped frequencies, so that the coupling step can be skipped and, when it cannot, the second kind equation \eqref{eq:the:equation} is well conditioned. The two theorems below make these facts quantitative, with constants explicit in the Laplace parameter $s$, the wave speed and the distances between $\Omega$, $\Gamma$ and $\Sigma$; the proofs are given in Section~\ref{proofs}. Before stating them let us introduce, for any pair of sets $A,B\subset\mathbb{R}^m$,
\[
 d(A,B) =\inf\{|{\bf x}-{\bf y}|\ : \ {\bf x}\in A,\ {\bf y}\in B\},
\]
the distance between $A$ and $B$.

 \begin{theorem}\label{theo:the:theo:01}
Let $\sigma_0>0$ and let $\varepsilon>0$ be sufficiently small. Set 
\[
 d:= (1-\varepsilon)\,d(\Omega,\Sigma),\quad
 c_{\rm max} := \|c\|_{L^\infty(\mathbb{R}^m)}\ (\ge 1),\quad
 c_{0,\rm R} := \|c_{\rm R}\|_{L^\infty(\Omega_2)} = \|c^{-2}-1\|_{L^\infty(\Omega_2)}.
\]
Then there exists $C>0$, depending only on $\sigma_0$, $\varepsilon$, $\Omega$, $\Gamma$ and $\Sigma$ (or $\Omega_1$, $\Omega_2$), such that for all $s$ with $\sigma =\Re s\ge\sigma_0$, all positive wave--speed fields $c$ with $c|_{\mathbb{R}^m\setminus\Omega}\equiv 1$ and all $F\in L^2(\Omega)$ the following holds:
\begin{equation}\label{eq:01:theo:the:theo:01}
\begin{aligned}
|s|\|W^{\rm tr}_F-U\|_{L^2(\Omega_2)}+& \|\nabla(W^{\rm tr}_F-U)\|_{L^2(\Omega_2)}\\
& \ \le \ C c_{\rm max}^2  e^{-\sigma  d}  \sigma^{-2} \,|s|^{m/2+1} (1+c_{\rm max}^2 c_{0,\rm R}\sigma^{-1}|s|)\|F\|_{L^2(\Omega_2)}.
\end{aligned}
\end{equation}
Besides, if $d_1=(1-\varepsilon) d(\Gamma,\Sigma) $,   $d_2=(1-\varepsilon)  d(\Gamma,\Omega) $ and $f_\Gamma=\gamma_{\Gamma}W^{\rm tr}_F$ is as in \eqref{eq:fgamma},
\begin{equation}\label{eq:02:theo:the:theo:01}
\| f_\Gamma\|_{H^{1/2}(\Gamma)}\le
C  c_{\rm max}^2  e^{-\sigma   \min\{d_1,d_2\}}   \,|s|^{(m-1)/2} \sigma^{-1}\|F\|_{L^2(\Omega_2)}.
\end{equation}
 \end{theorem}

\begin{theorem}\label{theo:the:theo:02}
Under the assumptions of Theorem \ref{theo:the:theo:01}, there exists $C>0$ such that
\begin{subequations}
\label{eq:01:theo:the:theo:02}
\begin{align}
 \|\gamma_\Gamma \mathrm{K}^{\rm tr}_{\Omega_2\Sigma}(s)\|_{H^{1/2}(\Sigma)\to H^{1/2}(\Gamma)}\ &\le\  C e^{-\sigma \min\{d_1,d_2\}}|s|^{m/2+1}\sigma^{-1}  \left(1+  c_{\rm max}^2 c_{0,\rm R}\,\sigma^{-1}|s| \right)\label{eq:01a:theo:the:theo:02},\\
 \|\gamma_\Sigma \mathrm{K}_{\Omega_1^{\rm c}\Gamma}(s)\|_{H^{1/2}(\Gamma)\to H^{1/2}(\Sigma)}\ &\le\
 C e^{-\sigma d_1}|s|^{m/2+1}\sigma^{-1}. \label{eq:01b:theo:the:theo:02}
\end{align}
\end{subequations}
\end{theorem}

\begin{remark} \label{remark:considerations}Some considerations before proceeding:

 \begin{enumerate}
 \item Observe, in view of \eqref{eq:01:theo:the:theo:01} in Theorem \ref{theo:the:theo:01},
 that $W_F^{\rm tr}$  may by itself be a sufficiently good approximation of $U$ in $\Omega_2$, which renders the FEM--BEM coupling unnecessary (in this case, the approximation of $U$ is taken to be zero in $\Omega_{2}^{\rm c}$).

 \item    Let $d:=\min\{d_1,d_2\}$ and, in this remark only,
\[
 \mathrm{K}_{\Gamma}(s):= \gamma_\Gamma \mathrm{K}^{\rm tr}_{\Omega_2\Sigma}(s)\,\gamma_\Sigma \mathrm{K}_{\Omega_1^{\rm c}\Gamma}(s):H^{1/2}(\Gamma)\to H^{1/2}(\Gamma),
\]
so that \eqref{eq:the:equation} reads $(\mathrm{I}-\mathrm{K}_\Gamma(s))g_\Gamma=f_\Gamma$. Theorem~\ref{theo:the:theo:02} gives
\begin{equation}\label{eq:contraction:factor}
 \| \mathrm{K}_{\Gamma}(s)\|_{H^{1/2}(\Gamma)\to H^{1/2}(\Gamma)}\ \le\ \rho(s):=C\,e^{-2\sigma d}\,|s|^{m+2}\sigma^{-2}\big(1+c_{\rm max}^2c_{0,\rm R}\sigma^{-1}|s|\big).
\end{equation}
The separation between $\Omega$, $\Gamma$ and $\Sigma$ enters only through the exponential factor, whose role is to compensate the polynomial growth in $|s|$ of the potentials. For strongly damped frequencies, namely when $2\sigma d>(m+3)\log|s|+\log C$, we have $\rho(s)<1$, $\mathrm{K}_\Gamma(s)$ is a contraction, and the Neumann series for $(\mathrm{I}-\mathrm{K}_\Gamma(s))^{-1}$ converges geometrically with rate $\rho(s)$. Since $2\min\{d_1,d_2\}\le d_1+d_2\le d(\Omega,\Sigma)$, this exponential rate is never faster than that of the boundary correction $U-W_F^{\rm tr}$ in \eqref{eq:01:theo:the:theo:01}; contraction, however, only requires $\rho(s)<1$, whereas discarding the correction requires it to fall below the prescribed tolerance, so that \eqref{eq:contraction:factor} still covers an intermediate range of strongly damped frequencies in which the BEM correction is small but not negligible.

In practice we solve \eqref{eq:the:equation} with GMRES. Its iterates are optimal on the Krylov spaces spanned by the partial sums of the Neumann series, at the same cost per iteration (one BEM and one FEM solve), so that its residual is never larger than that of the fixed point iteration. GMRES does not depend on this estimate, though: for every $s\in\mathbb{C}_+$ the operator $\mathrm{K}_\Gamma(s)$ is compact and $\mathrm{I}-\mathrm{K}_\Gamma(s)$ is injective (a solution of the homogeneous equation yields, via \eqref{eq:alg02:02} with $F=0$, a solution of the homogeneous problem \eqref{eq:theproblem:BEMFEM:01}, which must vanish), 
so that \eqref{eq:the:equation} is uniquely solvable and the spectrum of $\mathrm{I}-\mathrm{K}_\Gamma(s)$ accumulates only at $1$; since $\mathrm{K}_\Gamma(s)$ is a smoothing operator, its eigenvalues decay very fast, which is the setting in which GMRES converges superlinearly.
 All of this is stated at the continuous level; after discretization it holds up to mesh--dependent norm equivalences.

The polynomial factors in \eqref{eq:contraction:factor} prevent $\rho(s)<1$ in the weakly damped window $\Re s\approx\sigma_0$, $|s|\gg\sigma$, which is where the coupling is mostly activated in the experiments of Section~\ref{sec:numerics}. In those experiments GMRES nevertheless converges in very few iterations---typically one or two, never more than four---uniformly in the frequency, which suggests that these factors are an artifact of the proof and that the exponential factor alone governs the convergence.

 \end{enumerate}

\end{remark}
\section{Operator estimates and proofs of main results}
\label{sec:estimates}
The aim of this section is to derive estimates for the solution operator $\mathrm{K}(s)$ and for the operators of the coupled formulation, explicit in the Laplace parameter $s$. The former are the ingredient of the convergence analysis of Theorem~\ref{th:main};
the latter, together with the decay estimates of Propositions~\ref{prop:estimates:0} and \ref{prop:estimates:02}, yield the proofs of Theorems~\ref{theo:the:theo:01} and \ref{theo:the:theo:02} and quantify when the BEM component is truly necessary in our algorithm.

\subsection{Preliminaries}
In this subsection only, let $\Omega$ be an arbitrary bounded Lipschitz domain with boundary $\Sigma$, and consider the boundary value problem:
\begin{equation}\label{eq:01:sec:03}
\left|
 \begin{array}{l}
 W \in H^1(\mathbb{R}^m \setminus \Sigma) = H^1(\Omega \cup \Omega^{\rm c}),\\
 s^2 W - \Delta W = 0,\quad \text{in } \mathbb{R}^m \setminus \Sigma,\\
 \gamma_\Sigma W = g_\Sigma. 
 \end{array}
 \right.
\end{equation}
Then 
\begin{align} 
 \|\nabla W\|_{L^2(\mathbb{R}^m \setminus \Sigma)} + |s| \|W\|_{L^2(\mathbb{R}^m)} &\le C_{\Sigma} (1+ |s|)^{1/2} |s| \sigma^{-1} \|g_\Sigma\|_{H^{1/2}(\Sigma)} \label{eq:02:sec:03} 
\end{align} 
see, for example, \cite{Sayas:2016} or \cite[Lemma~4.13 and Theorem~4.15]{BanjaiSayas2023}. The proof relies on a closely related estimate for $s>0$ originally proven in \cite{BambergerHaDuong1986}. We will use estimate \eqref{eq:02:sec:03} later in this section.

The goal of this section is twofold. First, we aim to derive analogs of estimates \eqref{eq:02:sec:03} for the operator $\mathrm{K}(s): L^2(\mathbb{R}^m) \to H^r(\mathbb{R}^m)$ with $s \in \mathbb{C}_+$, defined by

\begin{equation}\label{eq:06:sec:03}
 U = \mathrm{K}(s)F, \quad \left|
 \begin{array}{l}
 U \in H^1(\mathbb{R}^m),\\
 \dfrac{s^2}{c^2} U - \Delta U = F.
 \end{array}
 \right.
\end{equation}
Here, $c \in L^\infty(\mathbb{R}^m)$ is a positive function such that $c \equiv 1$ in $\Omega^{\rm c}$.
For any open domain $D \subset \mathbb{R}^m$, we denote the inner product in $L^2(D)$ by
\begin{eqnarray}\label{eq:08:sec:03}
(U,V)_D &:=& \int_D U\, \overline{V}\, 
\end{eqnarray}
and the associated sesquilinear form by
\begin{eqnarray}\label{eq:07:sec:03}
a_{s,c,D}(U,V) &:=& s^2 (c^{-2}U, V)_D + (\nabla U, \nabla V)_D,
\end{eqnarray}
with the convention that for vector-valued functions, such as $\nabla U$, the inner product is understood componentwise.

With this notation, problem \eqref{eq:06:sec:03} is equivalent to the weak formulation
\begin{equation}\label{eq:09:sec:03}
\left|
\begin{array}{l}
U \in H^1(\mathbb{R}^m),\\
a_{s,c,\mathbb{R}^m}(U,V) = (F,V)_{\mathbb{R}^m}, \qquad \forall V \in H^1(\mathbb{R}^m).
\end{array}
\right.
\end{equation}
Since $\Delta U = s^2c^{-2}U - F\in L^2(\mathbb{R}^m)$ and
\begin{equation}\label{eq:equivalent:norms}
 \|W\|_{H^2(\mathbb{R}^m)}^2 \le 2\left(\|W\|_{L^2(\mathbb{R}^m)}^2 + \|\Delta W\|_{L^2(\mathbb{R}^m)}^2\right),
\end{equation}
which follows from the Fourier characterization of the Sobolev norms on $\mathbb{R}^m$ \cite[Ch.~3]{mclean:2000} and the elementary inequality $(1+t)^2\le2(1+t^2)$, we conclude that in fact $U\in H^2(\mathbb{R}^m)$. For smoother $F$, $U$ gains regularity.

 In the next result, we establish essential continuity constants for the operator $\mathrm{K}(s)$ in different norms, merely assuming that $F \in L^2(\mathbb{R}^m)$.

 Recall, before proceeding, that we have followed the convention 
\[
 \sigma = \Re s>0,
 \]
 which yields the elementary bound 
\begin{equation}\label{eq:basic:ineq}
1+|s|^\alpha\le 2\,\underline{\sigma}^{-\alpha}|s|^{\alpha},\qquad \alpha\ge 0,
\end{equation}
with 
\begin{equation}\label{eq:underline:sigma}
 \underline{\sigma} := \min\{1, \sigma\},
\end{equation}
a bound that will be used repeatedly in the sequel.

\subsection{Estimates for the solution operator $\mathrm{K}(s)$ and for related operators on truncated domains}
\begin{proposition}\label{prop:estimates}
The following bounds hold:
\begin{subequations}\label{eq:prop:estimates}
\begin{eqnarray} 
\label{eq:01:prop:estimates}
\|\mathrm{K}(s)\|_{L^2(\mathbb{R}^m)\to L^2(\mathbb{R}^m)} &\le&  c_{\rm max}^{2}\sigma^{-1} |s|^{-1}\\
\label{eq:02:prop:estimates}
\|\nabla \mathrm{K}(s)\|_{L^2(\mathbb{R}^m)\to L^2(\mathbb{R}^m)} &\le& \tfrac12 c_{\rm max}\sigma^{-1}, \\
\|\Delta \mathrm{K}(s)\|_{L^2(\mathbb{R}^m)\to L^2(\mathbb{R}^m)} &\le& 1+(c_{\rm max}/c_{\rm min})^2\sigma^{-1} |s|\le 2 (c_{\rm max}/c_{\rm min})^2\sigma^{-1} |s|,
\label{eq:03:prop:estimates} \\
\| \mathrm{K}(s)\|_{L^2(\mathbb{R}^m)\to H^1(\mathbb{R}^m)} &\le&    c_{\rm max}^2 \left(1+\left|\frac{s}{2c_{\rm  max}}\right|^2\right)^{1/2}\sigma^{-1}|s|^{-1}\label{eq:04a:prop:estimates}
\\ 
&\le &\sqrt{2} c_{\rm max}^2\underline{\sigma}^{-1}\sigma^{-1}
\label{eq:04:prop:estimates}
\end{eqnarray}
\end{subequations}
with
\[
c_{\rm max} = \|c\|_{L^\infty(\mathbb{R}^m)}\ge 1,\quad c_{\rm min} = \operatorname*{ess\,inf} c \le 1.
\]
Furthermore, for $r\in  [1,2], $
 \begin{equation}
\begin{aligned}
\label{eq:05:prop:estimates}
\|\mathrm{K}(s)\|_{L^2(\mathbb{R}^m)\to H^r(\mathbb{R}^m)} \ &\le \
 c_{\rm max}^2 \left[2\sqrt{2}\,   c_{\rm min}^{-2}  \right]^{r-1}  (1+|s|^2)^{r/2}\sigma^{-1} |s|^{-1}\\
 &\le \ \sqrt{2}c_{\rm max}^2 \left[2\sqrt{2}\,   c_{\rm min}^{-2}   \right]^{r-1}  \underline{\sigma}^{-r}\sigma^{-1} |s|^{r-1}. 
 \end{aligned}
 \end{equation}
\end{proposition}
\begin{proof}
The sesquilinear form associated to $\mathrm{K}(s)$ is coercive after a rotation in the complex plane: taking $V =   s/|s|\, U$, we have
\[
s |s| \|c^{-1} U\|^2_{L^2(\mathbb{R}^m)}+\frac{\overline{s}}{|s|}\|\nabla U\|^2_{L^2(\mathbb{R}^m)} =a_{s,c,\mathbb{R}^m}(U,V) = (F,V)_{\mathbb{R}^m}
 \]
 and therefore
\[
\begin{aligned}
  \|F\|_{L^2(\mathbb{R}^m)}\|U\|_{L^2(\mathbb{R}^m)}&\ \ge\ |(F,V)_{\mathbb{R}^m}|\ge  \Re ( a_{s,c,\mathbb{R}^m}(U,V) )\\
 &\ \ge\  \sigma\left( c_{\rm max}^{-2}|s| \|U\|^2_{L^2(\mathbb{R}^m)}+ \frac{1}{|s|}\|\nabla U\|^2_{L^2(\mathbb{R}^m)}\right)\\
  &\ \ge\  2 \sigma  \left( c_{\rm max}^{-1}  \|U\|_{L^2(\mathbb{R}^m)}\|\nabla U\| _{L^2(\mathbb{R}^m)}\right).
 \end{aligned}
\]
(In the last inequality we have used $a^2+b^2\ge 2ab$). Hence
\[
 \begin{aligned}
 \|U\|_{L^2(\mathbb{R}^m)}\ &\le\
   c_{\rm max} ^2  \sigma^{-1} |s|^{-1}  \|F\|_{L^2(\mathbb{R}^m)}, \\ 
 \|\nabla U\|_{L^2(\mathbb{R}^m)}\ &\le\ 
 \tfrac12 c_{\rm max} \sigma^{-1}   \|F\|_{L^2(\mathbb{R}^m)}
 .\quad && 
   \end{aligned}
\]

Clearly,  \eqref{eq:03:prop:estimates} follows from
\[
  \|\Delta U\|_{L^2(\mathbb{R}^m)} \le  \|F\|_{L^2(\mathbb{R}^m)} + |s/c_{\rm min}|^2\|U\|_{L^2(\mathbb{R}^m)}.
\]
On the other hand, we invoke \eqref{eq:equivalent:norms} to derive
\begin{eqnarray*}
\|U\|_{H^2(\mathbb{R}^m)}^2
& \le & 2\Bigl[ \|U\|_{L^2(\mathbb{R}^m)}^2 + \|\Delta U\|_{L^2(\mathbb{R}^m)}^2 \Bigr] \\[1ex]
& \le & 4\Bigl[ (1 + |s/c_{\min}|^{4}) \|U\|_{L^2(\mathbb{R}^m)}^2
                 + \|F\|_{L^2(\mathbb{R}^m)}^2 \Bigr] \\[1ex]
& \le & 4\Bigl[ (1 + |s/c_{\min}|^{4})\, c_{\max}^4 \sigma^{-2} |s|^{-2} + 1 \Bigr]
                 \|F\|_{L^2(\mathbb{R}^m)}^2 \\[1ex]
& \le & 8\,(1 + |s/c_{\min}|^{4})\, c_{\max}^4 \sigma^{-2} |s|^{-2}\,
        \|F\|_{L^2(\mathbb{R}^m)}^2 \\[1ex]
& \le & 8\,(c_{\max}/c_{\min})^4 (1 + |s|^{4})\, \sigma^{-2} |s|^{-2}\,
        \|F\|_{L^2(\mathbb{R}^m)}^2 \\[1ex]
& \le & 16\,(c_{\max}/c_{\min})^4\, \underline{\sigma}^{-4}\sigma^{-2} |s|^{2}\, 
        \|F\|_{L^2(\mathbb{R}^m)}^2 ,
\end{eqnarray*}
which proves \eqref{eq:05:prop:estimates} for $r=2$. (Notice that \eqref{eq:basic:ineq} has been used in the  last inequality.)

Also, from \eqref{eq:01:prop:estimates}--\eqref{eq:02:prop:estimates},
\[
\begin{aligned}
  \|U\|_{H^1(\mathbb{R}^m)} \ &=\  \left[ \|U\|_{L^2(\mathbb{R}^m)}^2+\|\nabla U\|_{L^2(\mathbb{R}^m)}^2\right]^{1/2}
  \le c_{\max}^2\Big(1+\Big|\frac{s}{2c_{\max}}\Big|^2\Big)^{1/2}\sigma^{-1}|s|^{-1}\|F\|_{L^2(\mathbb{R}^m)}\\
  &\le c_{\max} ^2(1+|s|^2)^{1/2}\,\sigma^{-1}\,|s|^{-1}\|F\|_{L^2(\mathbb{R}^m)}
  \le\  \sqrt{2}\, c_{\max} ^2\,\underline{\sigma}^{-1}\,\sigma^{-1}\|F\|_{L^2(\mathbb{R}^m)},
\end{aligned}
\]
which proves \eqref{eq:04a:prop:estimates}--\eqref{eq:04:prop:estimates} and \eqref{eq:05:prop:estimates} for $r=1$. The result for $r\in(1,2)$ follows by interpolation.
\end{proof}

In order to simplify the exposition we will use  in this section 
the usual $s$--dependent norm, equivalent to the $H^1(D)$ norm,
\[
 \|U\|_{s,D}^2=
 |s|^2\|U\|^2_{L^2(D)}+
 \|\nabla U\|^2_{L^2(D)}.
\]
Clearly,
\[
 \mathrm{K}(s) F|_{\Omega_2} =  \mathrm{K^{tr}}(s) (F ,\gamma_{\Sigma} \mathrm{K}(s)F),
\]
where
\begin{equation}\label{eq:def:Ktr}
\mathrm{K^{tr}}(s)(F,g_\Sigma):= U,\quad\text{with}\quad
\left|
 \begin{array}{l}
  U \in H^{1}(\Omega_2),\  \gamma_{\Sigma} U = g_\Sigma,\\
  \dfrac{s^2}{c^2} U    - \Delta U  = F,\quad \text{in }\Omega_2.
 \end{array}
 \right.
\end{equation}

\begin{proposition}\label{prop:sec5.1:1}For any $s$ with  $\sigma =\Re s>0$ the following holds:
\begin{equation}
\begin{aligned}
\label{eq:00:theo:sec5.1:1}
\|\mathrm{\mathrm{K}^{\rm tr}}(s)(F,g_\Sigma)\|_{s,\Omega_2}
&\le
C \sigma^{-1}\,c_{\rm max}^2 \, \|F\|_{L^2(\Omega_2)} \\
&\quad+
C \sigma^{-1}\,(1+|s|)^{1/2}\,|s|
\left(1+  c_{\rm max}^2\, c_{0,\rm R}\,\sigma^{-1}|s|\right)
   \|g_\Sigma\|_{H^{1/2}(\Sigma)} .
\end{aligned}
\end{equation}
where $C$ is a constant depending only on the domain $\Omega_2$, and
\[
 c_{\rm max} = \|c\|_{L^\infty(\mathbb{R}^m)}\ge 1,
 \qquad c_{0,\rm R} = \|c_{\rm R}\|_{L^\infty(\Omega_2)}=\|c^{-2}-1\|_{L^\infty(\Omega_2)}.
\]
\end{proposition}
\begin{proof}
Throughout the proof, $C$ denotes a constant depending only on $\Omega_2$, which may change from line to line.

Let $U_1$ and $U_2$ be the respective solutions of the problems
\[
 \left|
 \begin{array}{l}
  U_1 \in H^{1}(\Omega_2),\\
  -\Delta U_1 +s^2 U_1  = 0,\\
  \gamma_\Sigma U_1= g_\Sigma,
 \end{array}
 \right.\qquad
 \left|
 \begin{array}{l}
  U_2 \in H_0^{1}(\Omega_2),\\
  -\Delta U_2 +(s/c)^2 U_2  = F_2,\\
  \gamma_\Sigma U_2 = 0,
 \end{array}
 \right.
\]
with
\begin{equation}\label{eq:01:theo:sec5.1:1}
 F_2 := F -\bigl[c^{-2}-1\bigr]s^2U_1.
\end{equation}
Notice that $U_1+U_2 = \mathrm{K^{tr}}(s)(F,g_\Sigma)$, so it suffices to bound these two functions to prove the result.

For $U_1$, we can use \eqref{eq:02:sec:03} (see also \cite[Lemma 4.13]{BanjaiSayas2023}) to conclude
\begin{equation}\label{eq:02:theo:sec5.1:1}
|s| \|U_1\|_{L^2(\Omega_2)}+\|\nabla U_1\|_{L^2(\Omega_2)}
\le C \sigma^{-1}(1+|s|)^{1/2} |s| \|g_\Sigma\|_{H^{1/2}(\Sigma)}.
\end{equation}
As a consequence,
\begin{equation}\label{eq:03:theo:sec5.1:1}
 \|F_2\|_{L^2(\Omega_2)}\le  \|F\|_{L^2(\Omega_2)}
+ C\,c_{0,\rm R} \,\sigma^{-1}(1+|s|)^{1/2}|s|^2\|g_\Sigma\|_{H^{1/2}(\Sigma)}.
\end{equation}

On the other hand,  and similarly to \eqref{eq:09:sec:03}, we have that $U_2$ solves
\begin{equation}\label{eq:09:sec:03:bis}
\left|
  \begin{array}{l}
 U_2 \in H_0^1(\Omega_2),\\
 a_{s,c,\Omega_2}(U_2,V) = (F_2,V)_{\Omega_2},\quad \forall V  \in H_0^1(\Omega_2).
 \end{array}
 \right.
\end{equation}
Proceeding as in the proof of Proposition \ref{prop:estimates}, \begin{eqnarray}
\label{eq:01:prop:estimates:v2}
\|U_2\|_{L^2(\Omega_2)} &\le&   c_{\rm max}^{2}\sigma^{-1} |s|^{-1}\|F_2\|_{L^2(\Omega_2)}\\
\label{eq:02:prop:estimates:v2}
\|\nabla U_2\|_{L^2(\Omega_2)} &\le& \tfrac12 c_{\rm max}\sigma^{-1}\|F_2\|_{L^2(\Omega_2)},
\end{eqnarray}
From \eqref{eq:01:prop:estimates:v2}-\eqref{eq:02:prop:estimates:v2}, we derive
\begin{equation}\label{eq:04:theo:sec5.1:1}
\| U_2\|_{s,\Omega_2}  \le      \left(c_{\rm max}+\tfrac12\right)^2  \sigma^{-1}\|F_2\|_{L^2(\Omega_2)} .
\end{equation}
Combining \eqref{eq:01:theo:sec5.1:1}--\eqref{eq:04:theo:sec5.1:1}, the result follows.
\end{proof}

\begin{remark}\label{remark:sec:4}
Notice that we can also derive
\begin{eqnarray}
\| \mathrm{K}^{\rm tr }(s)(\cdot,0)\|_{L^2(\Omega_2)\to H_0^1(\Omega_2)} &\le&     c_{\rm max}^2 \left(1+\left|\frac{s}{2c_{\rm  max}}\right|^2\right)^{1/2}\sigma^{-1}|s|^{-1}.\label{eq:04a:prop:estimates:v2}
\end{eqnarray}
The extension of estimate \eqref{eq:05:prop:estimates} to $\mathrm{K}^{\rm tr }(\cdot,0):L^2(\Omega_2)\to H^r(\Omega_2)\cap H_0^1(\Omega_2)$ for $r\in(1,2]$ is somewhat more delicate, since we cannot guarantee that  $U=\mathrm{K}^{\rm tr}(F,0)\in H^2(\Omega_2)$. However, under the extra regularity assumption that the solution of
\begin{subequations}\label{eq:regularity}
\begin{equation}\label{eq:regularity:01}
\left\{
\begin{array}{l}
V \in H_0^1(\Omega_2),\\[2mm]
- \Delta V = G \quad \text{in } \Omega_2,
\end{array}
\right.
\end{equation}
belongs to $H^2(\Omega_2)$ and satisfies
\begin{equation}\label{eq:regularity:02}
\|V\|_{H^2(\Omega_2)}
\le C \|G\|_{L^2(\Omega_2)} = C \|\Delta V\|_{L^2(\Omega_2)},
\end{equation}
\end{subequations}
with $C$ independent of $G$, the proof can be adapted to show that
\begin{equation}\label{eq:convergenceH2Omega2}
\|U\|_{H^2(\Omega_2)}
\le C_{\Omega_2}
\left(\frac{c_{\max}}{c_{\min}}\right)^{2}
(1+|s|)^2\sigma^{-1}|s|\,
\|F\|_{L^2(\Omega_2)}.
\end{equation}
By interpolation we can extend the result to $H^r(\Omega_2)$, $r\in(1,2)$.

Such regularity holds   for domains with a $C^{1,1}$ boundary \cite[Thm.~2.2.2.3]{Grisvard1985} or convex domains in any dimension \cite[Thm.~3.2.1.2]{Grisvard1985}, in particular for convex polygonal and polyhedral domains, which is not a restrictive assumption in practice. In the general case, when $\Omega_2$ is only Lipschitz, full $H^2$-regularity
 fails: if $G \in L^2(\Omega_2)$, elliptic regularity results imply that $V \in H^{r}(\Omega_2)$ for any $r<3/2$ (see \cite{mclean:2000}), and the exponent $3/2$ cannot be improved in general, not even for $C^1$ domains \cite{Costabel:2019}.
\end{remark}

Notice that, provided $\sigma\sim |s|$, estimate \eqref{eq:00:theo:sec5.1:1} can be rewritten
\[
 \begin{aligned}
\|\mathrm{\mathrm{K}^{\rm tr}}(s)(F,g_\Sigma)\|_{s,\Omega_2}
&\le
C |s|^{-1}\,c_{\rm max}^2 \, \|F\|_{L^2(\Omega_2)}
+
C \,(1+|s|)^{1/2}
\left(1+  c_{\rm max}^2 c_{0,\rm R} \right)
   \|g_\Sigma\|_{H^{1/2}(\Sigma)}, 
\end{aligned}
\]
which shows an extra factor $|s|^{3/2}$ in the Dirichlet datum with respect to $F$.

Observe that, up to this point, we have not yet made use of the fact that $F$ has compact support in our applications. This property, combined with the fact that $\sigma = \Re s > 0$, implies that $\mathrm{K}(s)F$ decays exponentially to zero as the distance from $\Omega$ increases. Deriving a precise estimate for this decay constitutes the objective of the remainder of this section.

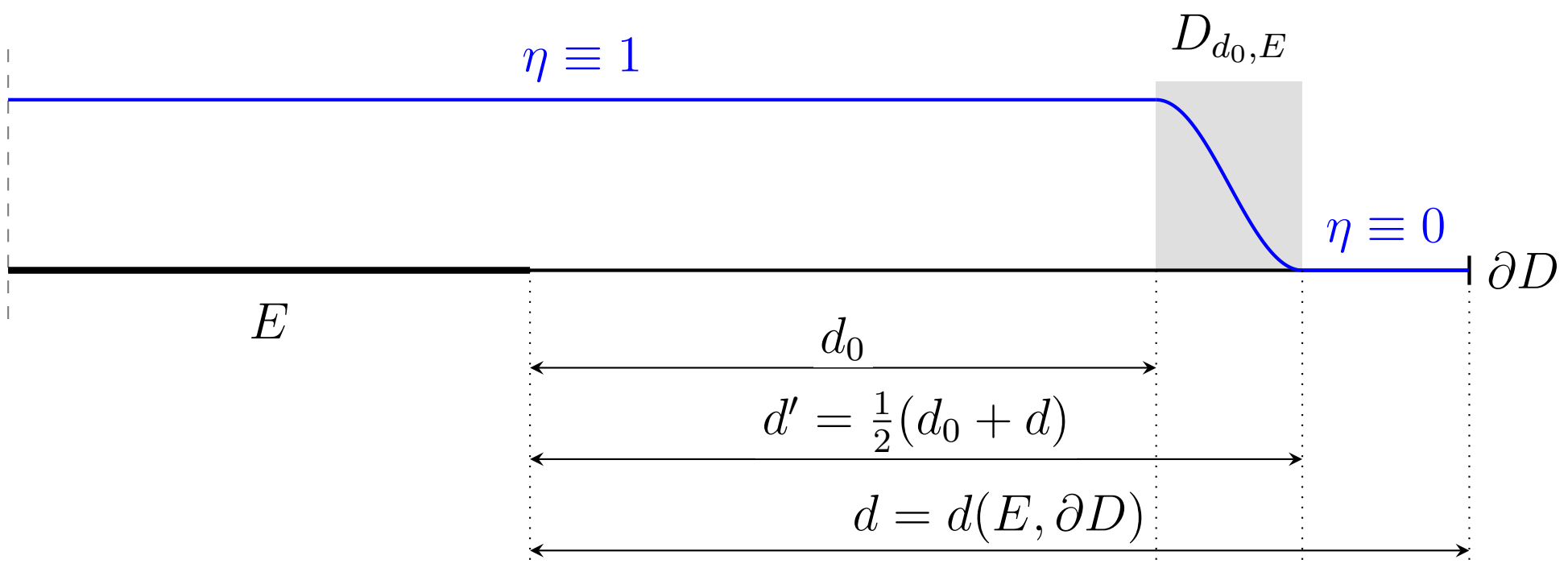
\begin{figure}[htb] 
\centering
\pgfmathsetlengthmacro{\xunit}{(\linewidth-1cm)/7}
\begin{tikzpicture}[>=stealth,x=0.8*\xunit]
 \def\Ea{2.5}\def\dZero{3}\def\dPrime{3.7}\def\dFull{4.5}\def\Hgt{1.4}
 \fill[gray!25] ({\Ea+\dZero},0) rectangle ({\Ea+\dPrime},{\Hgt+0.15});
 \node[above] at ({\Ea+0.5*(\dZero+\dPrime)},{\Hgt+0.15}) {$D_{d_0,E}$};
 \draw[dashed,gray] (0,-0.4) -- (0,{\Hgt+0.5});
 \draw[thick] (0,0) -- ({\Ea+\dFull},0);
 \draw[thick] ({\Ea+\dFull},-0.12) -- ({\Ea+\dFull},0.12);
 \node[right] at ({\Ea+\dFull},0) {$\partial D$};
 \draw[line width=1.6pt] (0,0) -- (\Ea,0) node[midway,below=3pt] {$E$};
 \draw[blue,thick] (0,\Hgt) -- ({\Ea+\dZero},\Hgt) node[midway,above] {$\eta\equiv1$};
 \draw[blue,thick,domain={\Ea+\dZero}:{\Ea+\dPrime},samples=50,smooth]
   plot (\x,{\Hgt*(0.5+0.5*cos(180*(\x-\Ea-\dZero)/(\dPrime-\dZero)))});
 \draw[blue,thick] ({\Ea+\dPrime},0) -- ({\Ea+\dFull},0) node[midway,above] {$\eta\equiv0$};
 \foreach \x in {\Ea,{\Ea+\dZero},{\Ea+\dPrime},{\Ea+\dFull}} \draw[dotted] (\x,0) -- (\x,-2.4);
 \draw[<->] (\Ea,-0.8)  -- ({\Ea+\dZero},-0.8)
   node[midway,above,fill=white,inner sep=1.5pt] {$d_0$};
 \draw[<->] (\Ea,-1.55) -- ({\Ea+\dPrime},-1.55)
   node[midway,above,fill=white,inner sep=1.5pt] {$d'=\tfrac12(d_0+d)$};
 \draw[<->] (\Ea,-2.3)  -- ({\Ea+\dFull},-2.3)
   node[midway,above,fill=white,inner sep=1.5pt] {$d=d(E,\partial D)$};
\end{tikzpicture}
\caption{\label{fig:cutoff}One--dimensional sketch  of the construction in the proof of Proposition~\ref{prop:estimates:0}: the cut--off function $\eta$ equals $1$ on $E_{d_0}$ and $0$ outside $E_{d'}$, so that $\nabla\eta$ is supported in the shell $D_{d_0,E}=\{d_0<d(\mathbf x,E)<d'\}$ (shaded), which lies at distance at least $d_0$ from $E$ and inside $D$.}
\end{figure}

\subsection{Exponential decay away from the inhomogeneity}
\begin{proposition}\label{prop:estimates:0}
Let $D\subset\mathbb{R}^m$, $m=2,3$, be open, let $E\subset D$ be bounded, and let $0<d_0<d:=d(E,\partial D)$ and $\sigma_0>0$. Set
\[
 d':=\tfrac12\big(d_0+d\big),\qquad
 D_{d_0,E}:=\big\{{\bf x}\in D:\ d_0<d({\bf x},E)<d'\big\}.
\]
Then for every multi-index $\bm{\alpha}\in\mathbb{N}_0^m$ there exists $C>0$, depending only on $m$, $|\bm\alpha|$, $\sigma_0$, $d_0$, $d-d_0$ and $|D_{d_0,E}|$, such that for every $s$ with $\sigma=\Re s\ge\sigma_0$ and every $W\in H^1(D)$ satisfying
\[
 s^2W-\Delta W=0\quad\text{in } D,
\]
the following holds:
\begin{equation}\label{eq:00:prop:estimates:0}
 \|\partial_{\bm{\alpha}}W\|_{L^\infty(E)}\ \le\ C\,e^{-\sigma d_0}\,|s|^{|\bm\alpha|+(m-1)/2} \|W\|_{L^2(D_{d_0,E})}.
\end{equation}
\end{proposition}
\begin{proof}
Let $\eta\in C^\infty(\mathbb{R}^m)$ be a cut-off function with $0\le\eta\le1$, $\eta\equiv1$ on $E_{d_0}$ and $\eta\equiv0$ outside $\overline{E}_{d'}$, with 
\[
E_{\varepsilon}:=\{{\bf x}\in\mathbb{R}^m \ : \ d({\bf x},E)<\varepsilon\}
\]
so that $\operatorname{supp}\nabla\eta\subset\overline{D}_{d_0,E}$. Constructing $\eta$ by mollification of the characteristic function of $E_{(d_0+d')/2}$ with a kernel of radius $\tfrac12(d'-d_0)$, we have in addition
\begin{equation}\label{eq:cutoff:prop:estimates:0}
 \|\partial_{\bm\beta}\eta\|_{L^\infty(\mathbb{R}^m)}\le C_{\bm\beta}\,\big(d-d_0\big)^{-|\bm\beta|},
\end{equation}
with $C_{\bm\beta}$ depending only on $m$ and $\bm\beta$, in particular neither on $E$ nor on $D$ (see Figure \ref{fig:cutoff} for 1D sketch of this construction).

Let $V:=\eta W$, extended by zero to $\mathbb{R}^m$. Then $V\in H^1(\mathbb{R}^m)$, $\operatorname{supp}V\subset\overline{E}_{d'}\subset D$, $V|_E=W|_E$, and
\[
 s^2V-\Delta V=f:=-2\nabla\eta\cdot\nabla W-(\Delta\eta)\,W\in L^2(\mathbb{R}^m),\qquad \operatorname{supp}f\subset\overline{D_{d_0,E}} .
\]
Since $\sigma=\Re s>0$, $V$ is the unique $H^1(\mathbb{R}^m)$ solution of this equation, and therefore
\[
V({\bf x}) = \int_{\mathbb{R}^m}\phi({\bf x}-{\bf y}) f({\bf y})\,{\rm d}{\bf y},\quad \text{where }
 \phi({\bf z};s)=\left\{\begin{array}{ll}
 \dfrac{{\rm i}}{4}H_0^{(1)}({\rm i}s|{\bf z}|), & m=2,\\[8pt]
 \dfrac{e^{-s|{\bf z}|}}{4\pi|{\bf z}|}, & m=3,
 \end{array}\right.
\]
is the fundamental solution of $s^2-\Delta$. For ${\bf x}\in E$ and ${\bf y}\in D_{d_0,E}$ we have $|{\bf x}-{\bf y}|\ge d_0$, so that $\phi({\bf x}-\cdot\,;s)$ is smooth on $\overline{D}_{d_0,E}$ and, since $\phi({\bf x}-\cdot\,;s)\nabla\eta$ has compact support in $D$, we may integrate by parts to obtain
\begin{equation}\label{eq:01:prop:estimates:0}
 W({\bf x})=V({\bf x})=\int_{D_{d_0,E}}\Big[\,2\nabla_{\bf y}\phi({\bf x}-{\bf y};s)\cdot\nabla\eta({\bf y})+\phi({\bf x}-{\bf y};s)\,\Delta\eta({\bf y})\Big]\,W({\bf y})\,{\rm d}{\bf y},\qquad {\bf x}\in E .
\end{equation}
Next, for $|{\bf z}|\ge d_0$, $\Re s\ge\sigma_0$ and any multi-index $\bm\beta$,
\begin{equation}\label{eq:02:prop:estimates:0}
 |\partial_{\bm\beta}\phi({\bf z};s)|\le C\,e^{-\sigma|{\bf z}|}\,|s|^{|\bm\beta|+(m-3)/2},
\end{equation}
with $C$ depending only on $m$, $|\bm\beta|$, $d_0$ and $\sigma_0$. For $m=3$ this follows by direct differentiation, each derivative producing a factor $|s|$ or $|{\bf z}|^{-1}\le d_0^{-1}\le|s|/(\sigma_0d_0)$. For $m=2$ it follows from $|H_\nu^{(1)}(w)|\le C_\nu|w|^{-1/2}e^{-\Im w}$, valid for $\Im w\ge0$ and $|w|\ge\sigma_0d_0$, applied with $w={\rm i}s|{\bf z}|$, so that $\Im w=\sigma|{\bf z}|$, together with $\big(H^{(1)}_0\big)'=-H^{(1)}_1$ and the recurrences for higher derivatives.

Differentiating \eqref{eq:01:prop:estimates:0} under the integral sign, applying \eqref{eq:02:prop:estimates:0} with $|\bm\beta|\le|\bm\alpha|+1$ and the Cauchy--Schwarz inequality, we obtain for ${\bf x}\in E$
\[
 |\partial_{\bm\alpha}W({\bf x})|\le C\,e^{-\sigma d_0}\,|s|^{|\bm\alpha|+(m-3)/2}\Big(|s|\,\|\nabla\eta\|_{L^2(D_{d_0,E})}+\|\Delta\eta\|_{L^2(D_{d_0,E})}\Big)\|W\|_{L^2(D_{d_0,E})}.
\]
By \eqref{eq:cutoff:prop:estimates:0},
\[
\begin{aligned}
\|\nabla\eta\|_{L^2(D_{d_0,E})}\ & \le C\big(d-d_0\big)^{-1}|D_{d_0,E}|^{1/2}\\
\|\Delta\eta\|_{L^2(D_{d_0,E})}\ &\le C\big(d-d_0\big)^{-2}|D_{d_0,E}|^{1/2}\le C\big(d-d_0\big)^{-2}\sigma_0^{-1}|s|\,|D_{d_0,E}|^{1/2},
\end{aligned}
\]
 which proves  \eqref{eq:00:prop:estimates:0}.
\end{proof}

The next result analyzes the decay of $\mathrm{K}(s) F$ far away from $\Omega$, the support of $F$.

\begin{proposition}\label{prop:estimates:02}
Let $\sigma_0>0$, let $E\subset\Omega^{\rm c}$ be bounded, let $0<d<d_0:=d(E,\Omega)$ and $\bm{\alpha}\in\mathbb{N}_0^m$. Then there exists $C>0$, depending only on $m$, $|\bm\alpha|$, $\sigma_0$, $d$, $d_0$ and $E$, such that for any $s$ with $\sigma=\Re s\ge\sigma_0$
\begin{equation}\label{eq:01:prop:estimates:02}
\|\partial_{\bm \alpha}\mathrm{K}(s)\|_{L^2(\Omega)\to L^\infty(E)}\le C c_{\rm max}^2\sigma^{-1} e^{-\sigma d}\, |s|^{(m-3)/2+|\bm{\alpha}|},
\qquad c_{\rm max} = \|c\|_{L^\infty(\mathbb{R}^m)}\ge 1.
\end{equation}
\end{proposition}
\begin{proof}
For $F\in L^2(\Omega)$, $W_1:=\mathrm{K}(s)F|_{\Omega^{\rm c}}\in H^1(\Omega^{\rm c})$ satisfies $s^2W_1-\Delta W_1=0$ in $\Omega^{\rm c}$, since $c\equiv1$ and $F\equiv 0$ there. Proposition~\ref{prop:estimates:0}, applied with $D=\Omega^{\rm c}$ and $d_0$ replaced by $d<d(E,\partial D)$, together with \eqref{eq:01:prop:estimates}, gives
\[
\begin{aligned}
 \|\partial_{\bm\alpha}W_1\|_{L^\infty(E)} \ &\le\   C e^{-\sigma d}|s|^{|\bm\alpha|+(m-3)/2}\,|s|\,\|\mathrm{K}(s)F\|_{L^2(\mathbb{R}^m)}\\
&\le\  C c_{\rm max}^2\sigma^{-1}e^{-\sigma d}|s|^{|\bm\alpha|+(m-3)/2}\|F\|_{L^2(\Omega)}.
\end{aligned}
\]
\end{proof}

\subsection { Proof of  Theorems \ref{th:main}, \ref{theo:the:theo:01} and \ref{theo:the:theo:02}}\label{proofs}

\subsubsection*{Proof of Theorem \ref{th:main}}
Following \cite{BanjaiSayas2023}, let $\mathcal{A}(\mu;\mathcal{B}(X,Y))$ denote the class of
analytic functions $\mathrm{F}:\mathbb{C}_+\to\mathcal{B}(X,Y)$ such that
\[
\|\mathrm{F}(s)\|\le C_{\mathrm F}(\sigma)\,|s|^{\mu},\quad  s\in\mathbb{C}_+
\]
 with
$C_{\mathrm F}:(0,\infty)\to(0,\infty)$ non-increasing and polynomially bounded as $\sigma\to0$.
By Proposition~\ref{prop:estimates}, estimate \eqref{eq:05:prop:estimates},
\[
 \mathrm{K}\in\mathcal{A}\big(r-1;\mathcal{B}(L^2(\Omega),H^r(\mathbb{R}^m))\big)
 \]
 where
 \[
 C_{\mathrm K}(\sigma)=C_r\,\underline{\sigma}^{-r}\sigma^{-1},
 \qquad
 C_r:=\sqrt2\,c_{\rm max}^2\,(2\sqrt2)^{\,r-1}\,c_{\rm min}^{-2(r-1)}, 
\]
with
\[
\underline{\sigma}:=\min\{1,\sigma\}.
\]
Since the left-hand side at time $t$ only depends on $f|_{[0,t]}$, we may extend $f$ by a
$C^{\widetilde{p}}$ function vanishing for $t\ge T+1$, so that $f\in W^{\widetilde{p}}_0(L^2(\Omega))$
in the notation of \cite{BanjaiSayas2023}. The result is then \cite[Thm.~2.21]{BanjaiSayas2023} for
BDF1 and BDF2, and \cite[Thm.~2.23]{BanjaiSayas2023} for the trapezoidal rule, applied with
$\mu=r-1\in[0,1]$. In both cases the regularity requirement reduces to
$\widetilde{p}>\mu+p+2=r+p+1$ (for the trapezoidal rule, $\max\{2\mu+3,\mu+4\}=\mu+4$ because $\mu\le1$),
and the constant is proportional to $\sigma^{p+r-\widetilde{p}}\,C_{\mathrm K}(\widetilde{\sigma})$ with
\[
 \widetilde{\sigma}=\tfrac12\sigma\min\{1,C_2\}<1
 \]
where $C_2$ depends only on $\delta$ and $\kappa_0$.
Hence, $C_{\mathrm K}(\widetilde{\sigma})=C_r\,\widetilde{\sigma}^{-(r+1)}$, which yields the stated
dependence on $\sigma$, $c_{\rm max}$ and $c_{\rm min}$.

\subsubsection*{Proof of Theorem \ref{theo:the:theo:01}}
Notice that by the trace theorem, if $\varXi$ is a bounded Lipschitz curve/surface contained in  a bounded Lipschitz domain $E$
\begin{equation}\label{eq:varXi}
\begin{aligned}
 \|\gamma_\varXi V\|_{H^{1/2}(\varXi)}\  &\le\  C(\varXi,E) \|V\|_{H^1(E)}\le C(\varXi,E) |E|^{1/2}\left( \|V\|_{L^\infty(E)}+\|\nabla V\|_{L^\infty(E)}\right)\\
 &\le\
 \underbrace{C(\Xi,E)\max\{1,\sigma_0^{-1}\}|E|^{1/2}}_{C' = C'(\sigma_0,\varXi,E)}\left( |s|\|V\|_{L^\infty(E)}+\|\nabla V\|_{L^\infty(E)}\right)
 \end{aligned}
\end{equation}
for all $s$ with $\sigma\ge \sigma_0$.

First, with the notation of this section
\[
 U|_{\Omega_2} = \mathrm{K^{tr}}(s)(F,\gamma_{\Sigma}\mathrm{K}(s)F),\quad
W^{\rm tr}_F = \mathrm{K}^{\rm tr}(s)(F,0).
\]
By Proposition \ref{prop:sec5.1:1},
\[
\|W^{\rm tr}_F-U\|_{s,\Omega_2}\le
C \sigma^{-1}\,(1+|s|)^{1/2}\,|s|
\left(1+  c_{\rm max}^2\, c_{0,\rm R}\,\sigma^{-1}|s|\right)
   \|\gamma_{\Sigma}\mathrm{K}(s)F\|_{H^{1/2}(\Sigma)} .
\]
Take now $E$ a bounded open neighborhood of $\Sigma$ with $d<d(E,\Omega)$. Proposition~\ref{prop:estimates:02} with $|\bm\alpha|\le1$ and \eqref{eq:varXi} with $\varXi=\Sigma$ give
\[
\|\gamma_\Sigma\mathrm{K}(s)F\|_{H^{1/2}(\Sigma)}\le C c_{\rm max}^2\sigma^{-1}e^{-\sigma d}|s|^{(m-1)/2}\|F\|_{L^2(\Omega)},
\]
and \eqref{eq:01:theo:the:theo:01} follows after absorbing the factor $(1+|s|)^{1/2}\le(1+\sigma_0^{-1})^{1/2}|s|^{1/2}$ into the constant.

The proof of \eqref{eq:02:theo:the:theo:01} is very similar, noticing first that
 \[
W^{\rm tr}_F = \mathrm{K}^{\rm tr}(s)(F,0)
 \]
 which yields a first bound by Proposition \ref{prop:sec5.1:1}:
 \[
  \|W^{\rm tr}_F\|_{s,\Omega_2}\le C c_{\rm max}^2\sigma^{-1}\|F\|_{L^2(\Omega_2)}.
 \]
Take now $\varXi= \Gamma$, $E$ a bounded open neighborhood of $\Gamma$ with $\tilde{d}= \min\{d_1,d_2\}<d(E,\Omega\cup\Omega_2^{\rm c})$,  $D= \Omega_2\cap \Omega^{\rm c}$ and notice that $W^{\rm tr}_F$
solves the constant--coefficient Helmholtz equation with complex wave number $\mathrm{i}\, s$ in $D$.
Then,  by \eqref{eq:varXi} and Proposition \ref{prop:estimates:0} (with $d_0=\tilde d$)  we conclude
\[
 \begin{aligned}
\| \gamma_\Gamma W^{\rm tr}_F \|_{H^{1/2}(\Gamma)}  \ &\le C  \left(|s|\|W^{\rm tr}_F\|_{L^\infty(E)}+\|\nabla W^{\rm tr}_F\|_{L^\infty(E)}\right) \le C' e^{-\sigma   \tilde{d}}  |s|^{(m-1)/2 } \| W^{\rm tr}_F\|_{s,D}\\
&\le\
C''  c_{\rm max}^2  e^{-\sigma   \tilde{d}}   \,|s|^{(m-1)/2 }\sigma^{-1}\|F\|_{L^2(\Omega_2)},
\end{aligned}
\]
which is \eqref{eq:02:theo:the:theo:01}.

\subsubsection*{Proof of Theorem \ref{theo:the:theo:02} } 
Bound \eqref{eq:01a:theo:the:theo:02} follows by the argument used to prove \eqref{eq:02:theo:the:theo:01}, with $W^{\rm tr}_F$ replaced by
\[
 V := \mathrm{K}^{\rm tr}_{\Omega_2\Sigma}(s)g_\Sigma= \mathrm{K}^{\rm tr}(s)(0,g_\Sigma),
\]
which also solves the constant--coefficient Helmholtz equation in  $D= \Omega_2\cap \Omega^{\rm c}$ and satisfies
\[
\| V\|_{s,D}\le
\| V\|_{s,\Omega_2}\le C \sigma^{-1}(1+|s|)^{1/2}|s|
\left(1+  c_{\rm max}^2\, c_{0,\rm R}\,\sigma^{-1} |s|\right)\|g_\Sigma\|_{H^{1/2}(\Sigma)}
\]
as a consequence of Proposition \ref{prop:sec5.1:1}. Now, with $E$ a bounded open neighborhood of $\Gamma$ with $\tilde{d}= \min\{d_1,d_2\}<d(E,\Omega\cup\Omega_2^{\rm c})$ we get after applying Proposition \ref{prop:estimates:0} (with $d_0=\tilde d$) and \eqref{eq:varXi} (with $\varXi = \Gamma$)
\[
\begin{aligned}
\| \gamma_\Gamma V\|_{H^{1/2}(\Gamma)}  \ &\le C\left( |s|\|V\|_{L^\infty(E)}+\|\nabla V\|_{L^\infty(E)}\right) \\
&\le\ C' e^{-\sigma   \tilde{d}}  |s|^{(m-1)/2}\| V\|_{s,\Omega_2}\\
&\le\
C'' e^{-\sigma   \tilde{d}}  \,(1+|s|)^{1/2}\,|s|^{(m+1)/2}\sigma^{-1}
\left(1+  c_{\rm max}^2\, c_{0,\rm R}\,\sigma^{-1} |s|\right)\|g_\Sigma\|_{H^{1/2}(\Sigma)}.
\end{aligned}
\]
Absorbing again the factor $(1+|s|)^{1/2}$ into $|s|^{1/2}$, this proves \eqref{eq:01a:theo:the:theo:02}.

Finally, to prove (\ref{eq:01b:theo:the:theo:02}) we simply take $\varXi=\Sigma\subset E\subset D:=\Omega_1^{\rm c}$ in \eqref{eq:varXi}, $E$ being a bounded neighborhood of  $\Sigma$ with $d(E,\Gamma)>d_1$ so that for
\[
 W =\mathrm{K}_{\Omega_1^{\rm c}\Gamma}(s)g_\Gamma
\]
we can conclude
\[
\begin{aligned}
\|\gamma_\Sigma W\|_{H^{1/2}(\Sigma)}\ &\le\
C' e^{-\sigma   d_1}  |s|^{(m-1)/2}\| W\|_{s,\Omega_1^{\rm c}}\\
 &\le\
 C''e^{-\sigma   d_1}   |s|^{(m+1)/2}(1+|s|)^{1/2} \sigma^{-1}\|g_\Gamma\|_{H^{1/2}(\Gamma)},
  \end{aligned}
\]
where we have applied estimate \eqref{eq:02:sec:03} in the last step. Since $(1+|s|)^{1/2}\le (1+\sigma_0^{-1})^{1/2}|s|^{1/2}$, this proves \eqref{eq:01b:theo:the:theo:02}.

\section{Numerical experiments} \label{sec:numerics}

All computational results in this section were obtained with fully discrete counterparts of Algorithms~\ref{alg:allsteps} and~\ref{alg:coupling}: the all-steps-at-once Algorithm~\ref{alg:allsteps} of Subsection~\ref{subsec:cq:allsteps} organizes the computation, and at each quadrature frequency $s_\ell$ (by \eqref{eq:sym:01}, only half of them need to be computed) the evaluation of $\mathrm{K}(s_\ell)F$, i.e., the solution of problem~\eqref{eq:theproblem:BEMFEM:01}, is carried out independently and in parallel by the FEM--BEM coupling of Algorithm~\ref{alg:coupling} of Section~\ref{sec:coupling}. We next briefly describe how these components are discretized. (A rigorous convergence analysis of the spatial discretization is beyond the scope of this article, and is left for future work.)

\subsection{FEM--BEM discretization and coupling}
\emph{Finite elements.} For the heterogeneous bounded domain problem \eqref{eq:alg02:01} on $\Omega_2$ we use standard conforming $\mathbb{P}_2$ Lagrange elements. This simple choice is adequate here: for frequencies with $\Re s$ proportional to $|s|$, the sesquilinear form is coercive with constant $\Re s/|s|$ and hence the Galerkin solution is quasi-optimal with constants independent of $s$. Since $f$ is smooth in time,  the weakly damped frequencies case with $|s|$ large, which the mesh cannot resolve, carry a right-hand side $F^\kappa_{\ell,r}$ at round-off level. The mesh resolves $\partial\Omega$ and the interfaces between the pieces of $\Omega$, which is why $\Omega$ is chosen polygonal; $c$ is then smooth on each element and $c^{-2}$ is integrated by element--wise Gaussian quadrature of its analytic expression. The Dirichlet datum on $\Sigma$ is imposed by nodal interpolation at the vertices and edge midpoints lying on $\Sigma$. The resulting matrices, $s^2M_c+S$ with $M_c$ the $c^{-2}$--weighted mass matrix and $S$ the stiffness matrix, are complex symmetric but not Hermitian, and are solved either with a sparse direct solver or with the multigrid iteration described at the end of this section; in both cases the factorization (or the hierarchy) is computed once per frequency and reused in all the solves at that frequency.

\noindent\emph{Boundary elements.} The exterior problem defining $\mathrm{K}_{\Omega_1^{\rm c}\Gamma}(s)$ is solved by an indirect single layer potential representation, $W^{\rm c}=\mathrm{SL}({\rm i}s)\lambda$, with $\lambda$ solving the first--kind equation $\mathrm{V}({\rm i}s)\lambda=g_\Gamma$ on $\Gamma$, where ${\rm i}s$ is the wave number of the associated Helmholtz equation, see \cite{Kress:2014,mclean:2000}; the strong damping makes more elaborate (combined) formulations unnecessary, and we observed no difference when using them. For numerically solving the first--kind equation,  we use the deltaBEM method \cite{DoLuSa:2014a,DoLuSa:2014b}, a Nystr\"om--type fully discrete Calder\'on calculus of order three on a uniform grid of $n$ points in the parameter domain of $\Gamma$, which is local in nature, very simple to implement and requires no splitting of the kernel into regular and singular parts. (For the complex frequencies arising in our CQ formulation, such a splitting causes both terms to grow exponentially with $|\Re s|\,|\mathbf{z}|$, even though their sum decays exponentially, which leads to severe cancellation.)
 An alternative in the same spirit is a Galerkin discretization with piecewise polynomials. In either case it is worth noticing that, since the fundamental solution kernels decay like $e^{-\sigma|\mathbf{x}-\mathbf{y}|}$, the matrices become numerically banded --- entries corresponding to well--separated points fall below any prescribed tolerance --- precisely for strongly damped frequencies, that is, in the regime where, by Theorem~\ref{theo:the:theo:01}, the correction is least needed: the cost of the boundary term fades at the same rate as its influence on the solution. Conversely, for the few weakly damped frequencies the matrices are dense, but there the correction is essential. Finally, the trace $\gamma_\Sigma W^{\rm c}$ is obtained by evaluating the discrete potential at the Dirichlet nodes of the finite element mesh on $\Sigma$; since these lie at distance at least $d_1$ from $\Gamma$, this is a plain quadrature with a smooth kernel, encoded in a dense matrix assembled once per frequency.

\noindent\emph{Traces and coupling.} The trace $\gamma_\Gamma$ of a finite element function, needed for $f_\Gamma=\gamma_\Gamma W^{\rm tr}_F$ and for the operator $\mathrm{K}_\Gamma(s)$ of Remark~\ref{remark:considerations}, is computed by evaluating the $\mathbb{P}_2$ function at the $n$ deltaBEM points of $\Gamma$, which lie in the interior of the elements. With these two trace operators, the Schur complement equation \eqref{eq:the:equation} becomes a dense $n\times n$ system, which is never assembled: it is solved by GMRES, each iteration requiring one finite element solve, one boundary element solve and the two trace evaluations. The GMRES tolerance is $10^{-8}$ on the relative residual, and the whole coupling step is skipped whenever $\|\gamma_\Sigma\mathrm{K}_{\Omega_1^{\rm c}\Gamma}(s)f_\Gamma\|_\infty\le10^{-7}\|W^{\rm tr}_F\|_\infty$, in which case $U|_{\Omega_2}$ is taken to be $W^{\rm tr}_F$ (see \eqref{eq:alg02:02} in Algorithm \ref{alg:coupling}).

\subsection{Numerical example and simulation}
 \begin{figure}
\includegraphics[width=0.52\textwidth]{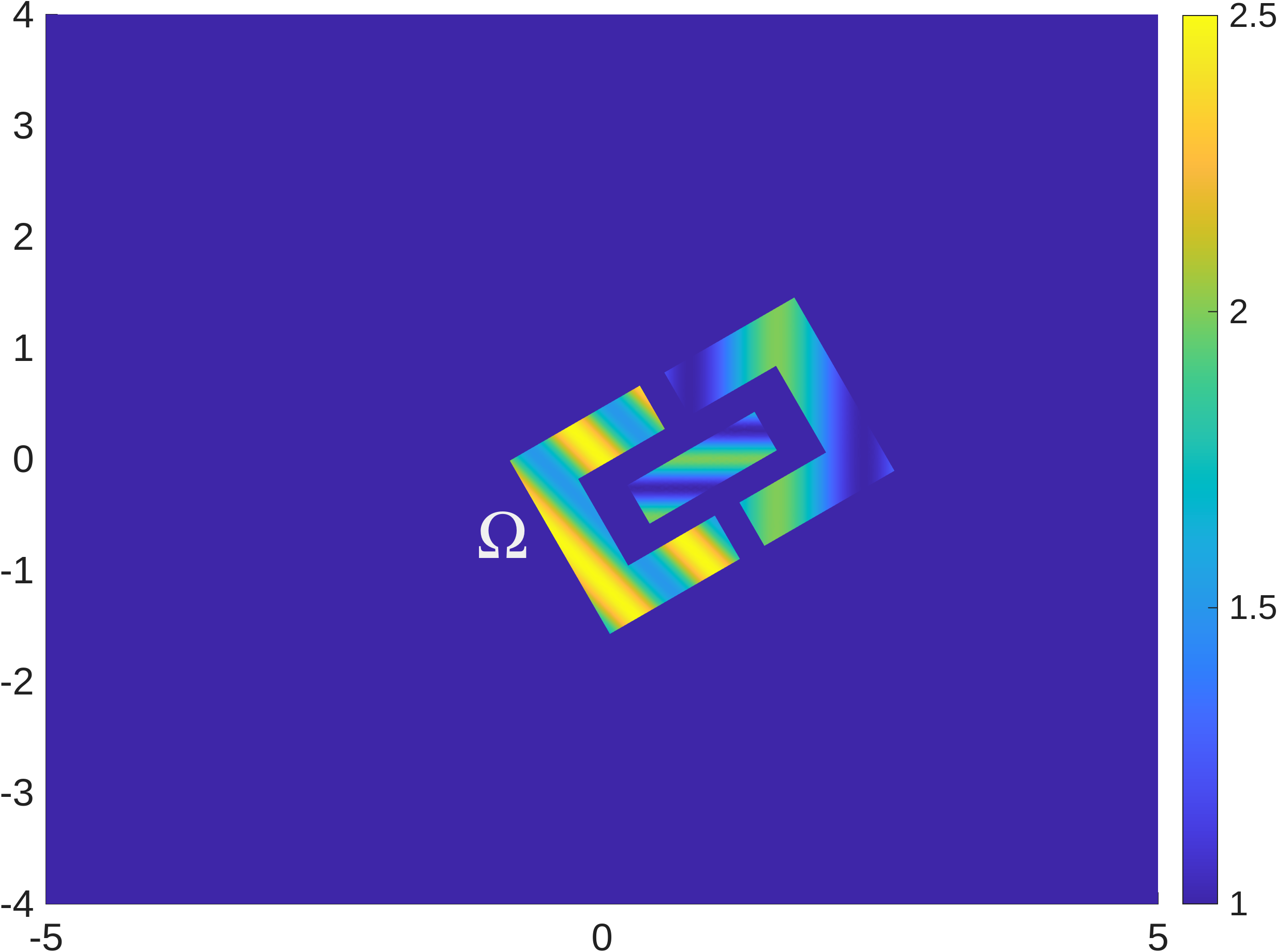}
\hfill
\includegraphics[width=0.47\textwidth]{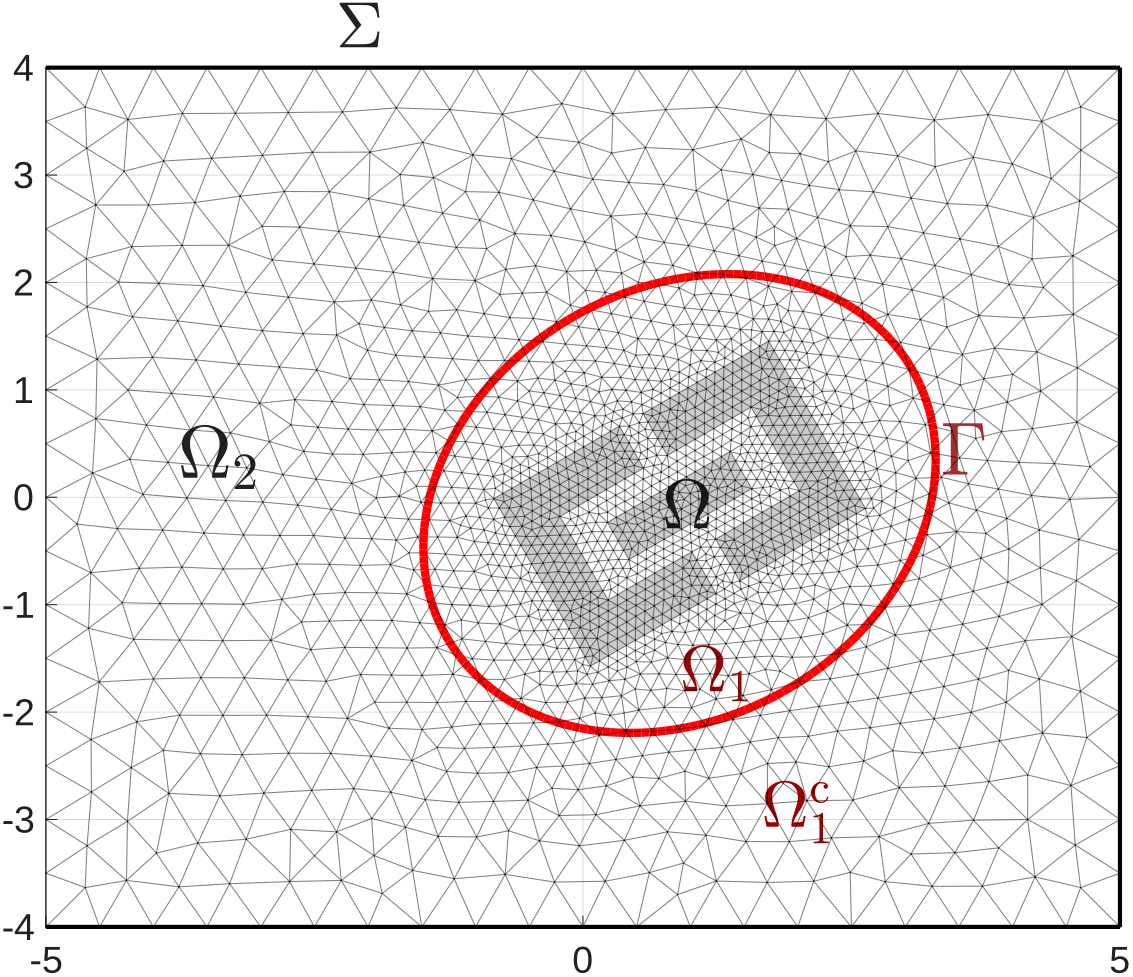}
\caption{\label{fig:exp01:config}Configuration of the experiment. Left: the wave speed $c$ in $\Omega_2=(-5,5)\times(-4,4)$, with the chain-link inhomogeneity $\Omega$ made up of three polygonal pieces. Right: initial finite element mesh, the ellipse $\Gamma$ (red) supporting the boundary element unknowns, and the domains $\overline{\Omega}\subset\Omega_1$ and $\overline{\Omega}_1\subset\Omega_2$.}
 \end{figure}

In this section we consider the example configuration depicted in Figure~\ref{fig:exp01:config}, based on the  heterogeneous region $\Omega$. The FEM domain is the rectangle $\Omega_2 = (-5,5)\times(-4,4)$, so that $\Sigma=\partial\Omega_2$, and the BEM curve $\Gamma$ is the ellipse with semi-axes $2.5$ and $2$, rotated by $\pi/6$ and centered at approximately $(0.90,-0.06)$, whose parameterization is given in \eqref{eq:ellipse}.

The heterogeneous domain $\Omega$ is the chain-link shaped region shown in Figure~\ref{fig:exp01:config}: a non-convex domain with re-entrant corners and holes, made up of three polygonal pieces $\Omega^{(1)}$ (left link), $\Omega^{(2)}$ (central overlap) and $\Omega^{(3)}$ (right link). The wave speed is $c\equiv 1$ outside $\overline{\Omega}$, whereas inside each piece it is a smooth, oscillatory perturbation of the background,
\begin{equation}\label{eq:exp01:c}
 c(x,y) = 1 + \sum_{j=1}^{3} f_j(x,y)\,\chi_{\Omega^{(j)}}(x,y),
 \qquad
 \begin{aligned}
  f_1(x,y) &= \tfrac12\big(2+\cos(4\sqrt2\,(x+y))\big),\\
  f_2(x,y) &= \tfrac12\big(1+\cos(12y)\big),\\
  f_3(x,y) &= \tfrac12\big(1+\cos(4x)\big),\\
 \end{aligned}
\end{equation}
where $\chi_{\Omega^{(j)}}$ denotes the characteristic function of
$\Omega^{(j)}$. Hence $c$ takes values in $[1,2]$ on $\Omega^{(1)}$ and
$\Omega^{(3)}$, and in $[3/2,5/2]$ on $\Omega^{(2)}$, so that
$c_{\min}=1$ and $c_{\max}=5/2$. Notice that $c$ is discontinuous across
$\partial\Omega$ (except at isolated points) and across the interfaces
between the three pieces, which produces a right-hand side $F$ with jumps,
as anticipated in Section~\ref{sec:estimates}. The three oscillation directions
$(1,0)$, $(1,1)/\sqrt2$ and $(0,1)$ and the increasing frequencies $4$,
$4\sqrt{2}$ and $12$ are chosen to avoid any accidental alignment with
the mesh.

The problem is driven by an incident plane pulse propagating in the
background medium, where $c\equiv 1$, with unit speed:
\begin{equation}\label{eq:exp01:uinc}
\begin{aligned}
 u^{\rm inc}(\mathbf{x},t)
 = \exp\Big(-\frac{\big(t-(\mathbf{x}-\mathbf{x}_0)\cdot\mathbf{d}\big)^2}{4\varsigma^2}\Big)
   \cos\big(\omega\,(t-(\mathbf{x}-\mathbf{x}_0)\cdot\mathbf{d})\big),\\
 \mathbf{x}_0 = (-3,0),\quad
 \mathbf{d} = \frac{(1,-0.2)}{\sqrt{1.04}},\quad
 \varsigma= 0.2,\quad \omega = \pi,
 \end{aligned}
\end{equation}
that is, a Gaussian-modulated wave front of width of order $\varsigma$ and
carrier frequency $\omega$, traveling in the direction $\mathbf{d}$. At
$t=0$ the pulse is centered on the line through $\mathbf{x}_0$ orthogonal
to $\mathbf{d}$, at distance about $2$ from $\overline{\Omega}$, so that
its trace on the inhomogeneity is of order $e^{-25}$ and the homogeneous
initial conditions in \eqref{eq:theproblem}
hold up to round-off. We choose the final time for our simulation to be $T=8$. 

The initial mesh, generated with Gmsh \cite{Gmsh:2009}, consists of \num{1146} triangles and \num{2365} nodes. Finer meshes are obtained by uniform
 refinements, so that the meshes form a nested hierarchy; this is what
makes multigrid solvers applicable, a point to which we return below. This initial mesh is fitted to the material interfaces, which comes at no cost here since the interfaces are polygonal, and it preserves the full order of the finite element approximation under refinement. For curved interfaces some accuracy is lost, less so if the mesh follows the shape of the inhomogeneity with straight-sided elements, as in~\cite{DoGaLe:2025}, than if the interface is ignored by the mesh. 
The BEM unknowns live on the ellipse with parameterization $\mathbf{p}$:
\begin{equation}\label{eq:ellipse}
 \Gamma:\quad
 \mathbf{p}(x) =
 \begin{pmatrix} \cos\tfrac{\pi}{6} & -\sin\tfrac{\pi}{6}\\[2pt]
                 \sin\tfrac{\pi}{6} & \phantom{-}\cos\tfrac{\pi}{6}\end{pmatrix}
 \begin{pmatrix} 0.75 + 2.5\cos x\\[2pt] -0.5 + 2\sin x\end{pmatrix},
 \qquad  x\in[0,2\pi),
 \end{equation}
which is discretized with $n=300$ points.

\paragraph{Error measurement metrics.}
Since no closed-form solution of the above example is available, the role of the exact solution for error measurements is
played by a reference solution $u^{\rm ref}$ computed with the same scheme and a
finer discretization of the parameter under study, the others being kept fixed:
the time step $\kappa=8/3072$ ($n_{\rm CQ}=3072$), half that of the finest
discretization considered, $n_{\rm CQ}=192,\,384,\,768,\,1536$, when the
convergence in time is examined; and one further refinement of the finite
element mesh or of the boundary element discretization when the convergence
in space is studied. For each time level $t_j=j\kappa$ of the numerical
solution under study we compute the relative error
\[
\big(E^{\kappa}(t_j)\big)^2 :=
\frac{\displaystyle\sum_{\ell}
      \big|u^{\kappa}({\bf x}_\ell,t_j)-u^{\rm ref}({\bf x}_\ell,t_j)\big|^2}
     {\displaystyle\sum_{\ell}\big|u^{\rm ref}({\bf x}_\ell,t_j)\big|^2}
\;\approx\;
\frac{\displaystyle\int_{\Omega_2}
      \big|u^{\kappa}({\bf x},t_j)-u^{\rm ref}({\bf x},t_j)\big|^2\,{\rm d}{\bf x}}
     {\displaystyle\int_{\Omega_2}\big|u^{\rm ref}({\bf x},t_j)\big|^2\,{\rm d}{\bf x}},
\]
where the sums run over the vertices and edge midpoints of the mesh of
$u^{\rm ref}$, on which the coarser finite element solution, being piecewise
polynomial, is evaluated exactly; this is a discrete relative $L^2(\Omega_2)$
error at time $t_j$. Here $u$ denotes the total wave; note that the error itself
is the same for the total and for the scattered field, since the incident wave
is known exactly, and the total wave is used only for normalization purposes
because its norm stays away from zero for all $t\in[0,T]$.

In addition to $E^{\kappa}(t)$ we report its maximum in time
\begin{equation}\label{def:Ekappa:infty}
E^{\kappa}_{\infty} := \max_j E^{\kappa}(t_j).
\end{equation}

\begin{figure}[htb]
\centering
\includegraphics[width=0.85\textwidth]{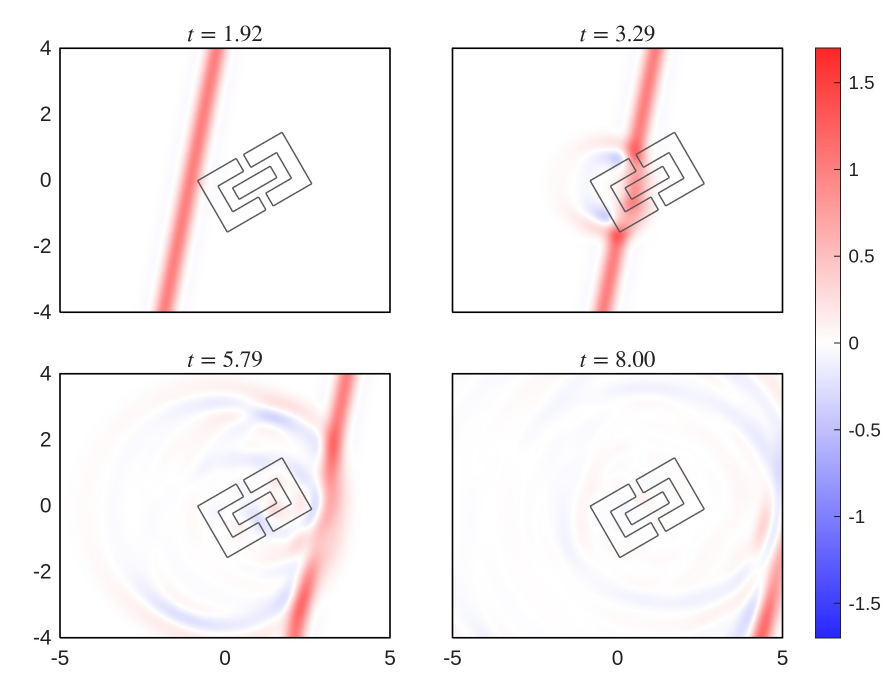}
\caption{Total field $u = u^{\rm inc}+u^{\rm sc}$ at
$t = 1.92$, $3.29$, $5.79$ and $8.00$, computed on the FEM domain
$\Omega_2=(-5,5)\times(-4,4)$, whose boundary $\Sigma$ is the frame of each
panel; the boundary of the inhomogeneity $\Omega$ is drawn in gray. The
color scale is the same in all panels and is saturated at $\pm1.7$. Since
$c>1$ inside $\Omega$, the pulse travels faster through the inhomogeneity
than in the background, and the wave front is visibly advanced after
crossing it ($t=5.79$). Note the absence of spurious reflections at
$\Sigma$, where the exact exterior condition is enforced through the BEM.}
\label{fig:exp01:total}
\end{figure}

\begin{figure}[htb]
\centering
\includegraphics[width=0.85\textwidth]{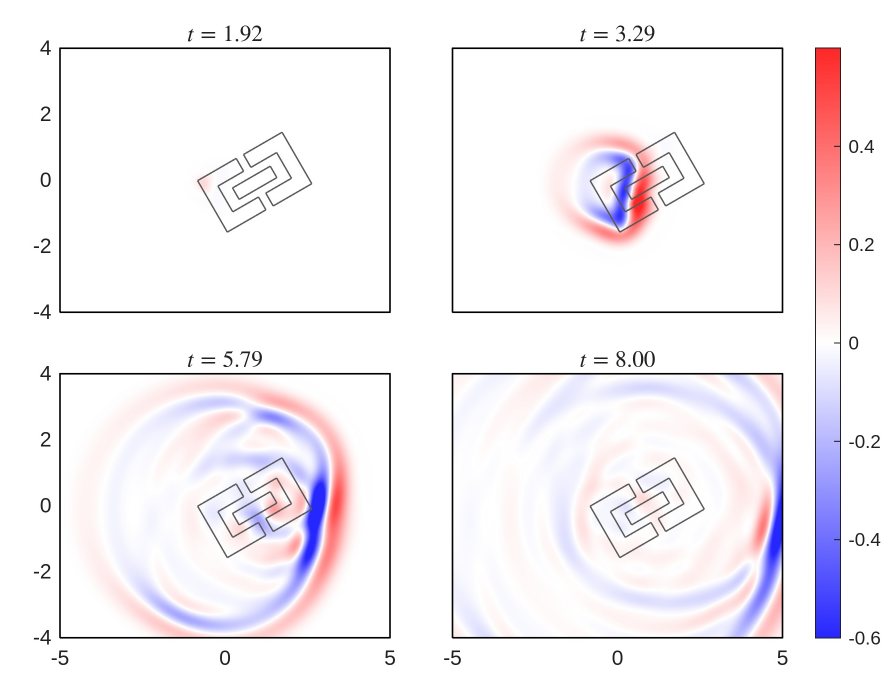}
\caption{Scattered field $u^{\rm sc}$ at the same times
as in Figure~\ref{fig:exp01:total}. The color scale is again common to all
panels, now saturated at $\pm0.6$ (the amplitude of $u^{\rm sc}$ is
considerably smaller than that of the incident pulse); saturated regions
mark the concentration of the scattered energy inside $\Omega$ and at its
right tip. The waves generated by the multiple internal reflections in
the chain-link geometry leave the computational domain through $\Sigma$
without any visible artifact.}
\label{fig:exp01:scatt}
\end{figure}

\begin{figure}[htb]\centering
\includegraphics[width=0.47\textwidth]{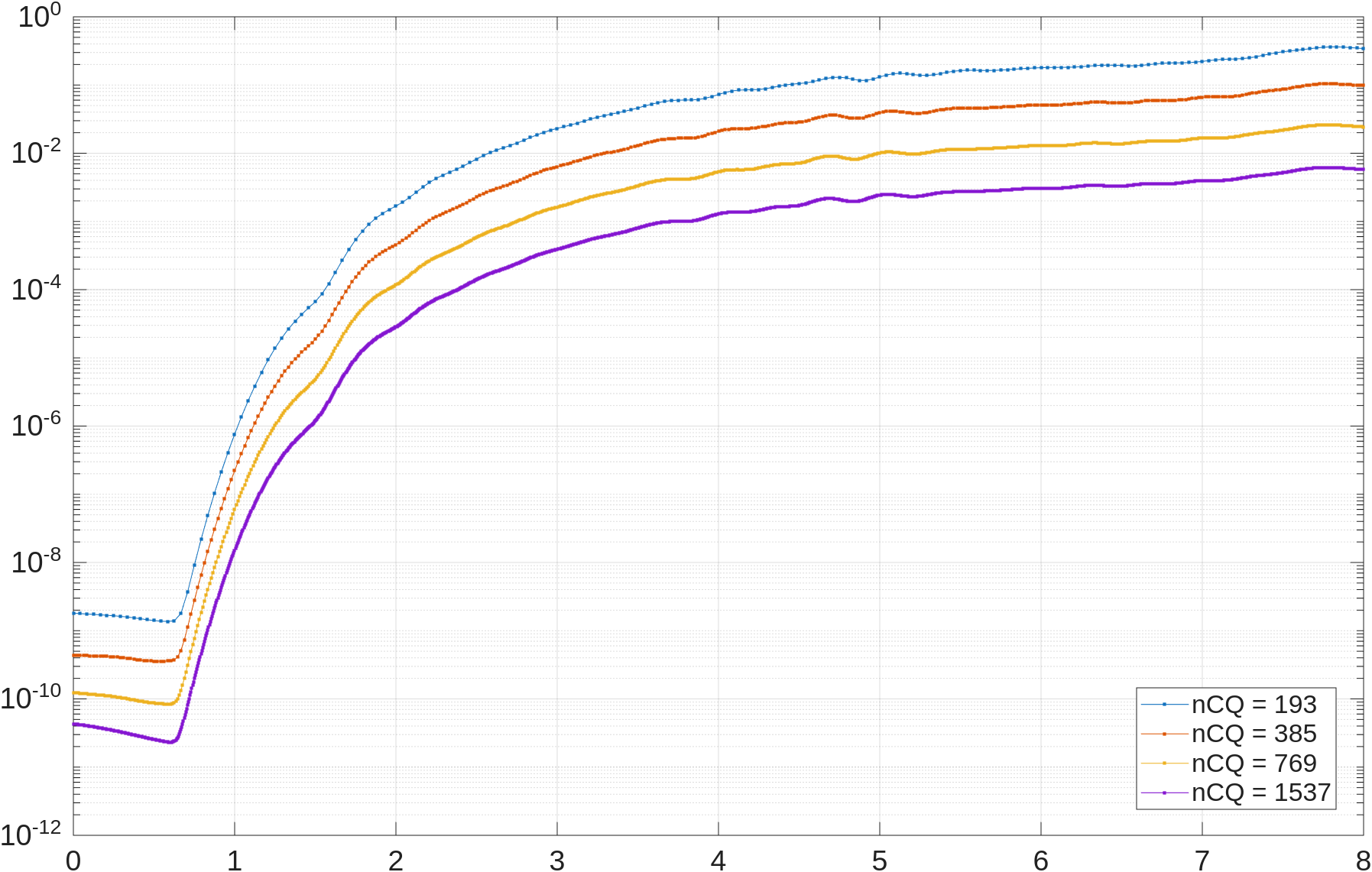}\quad
\includegraphics[width=0.47\textwidth]{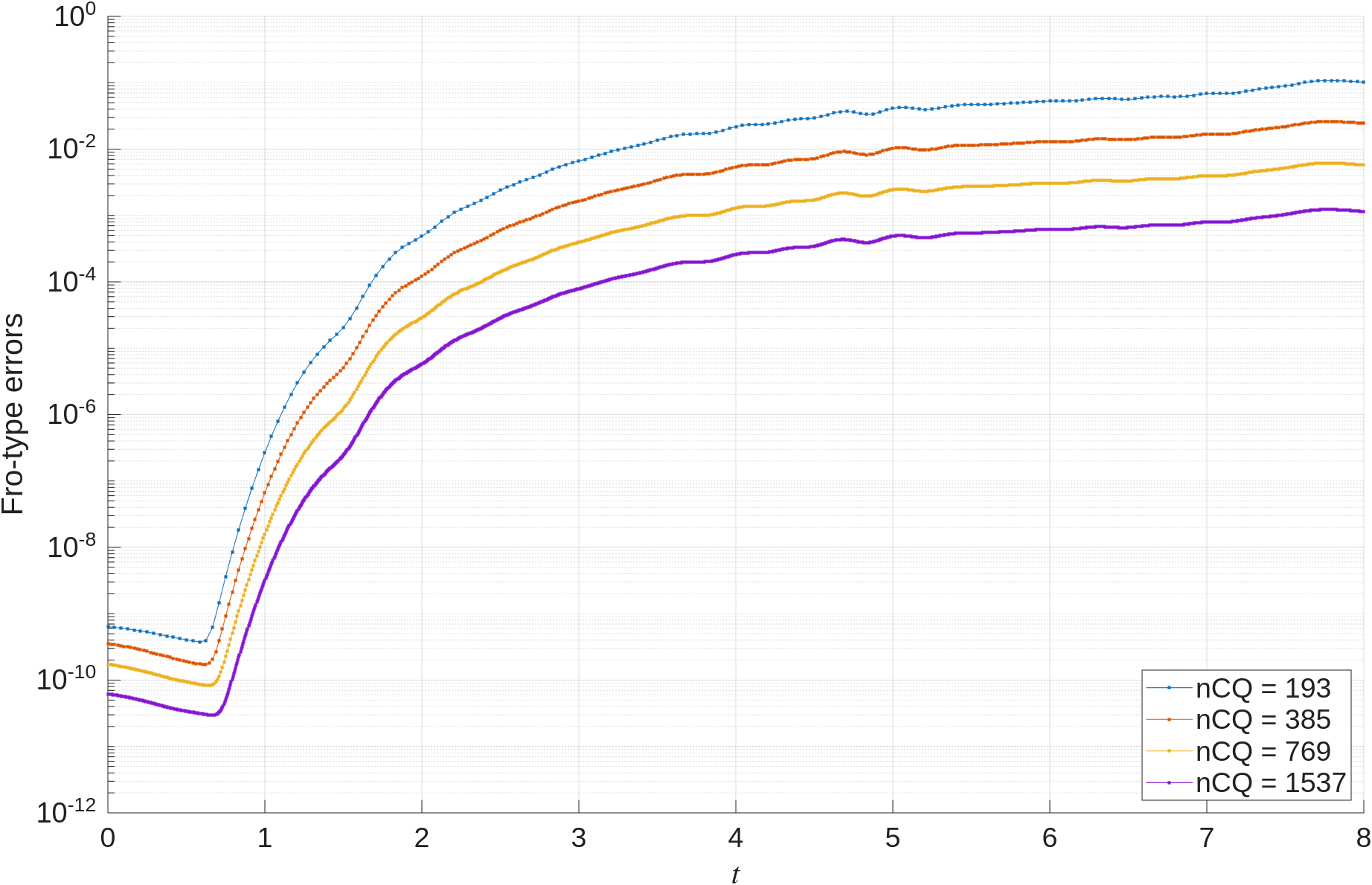}\\
\includegraphics[width=0.47\textwidth]{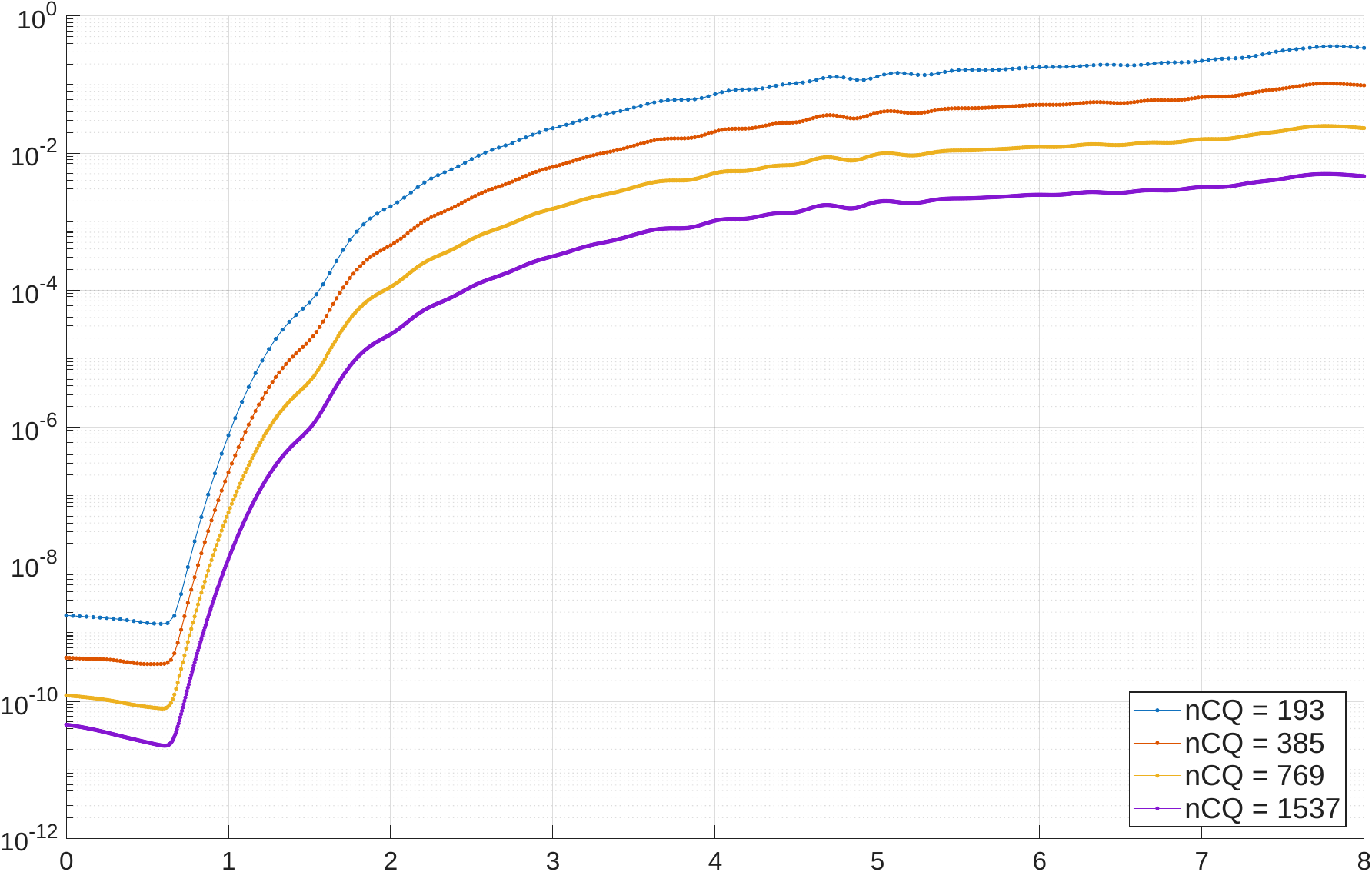}
\quad
\includegraphics[width=0.47\textwidth]{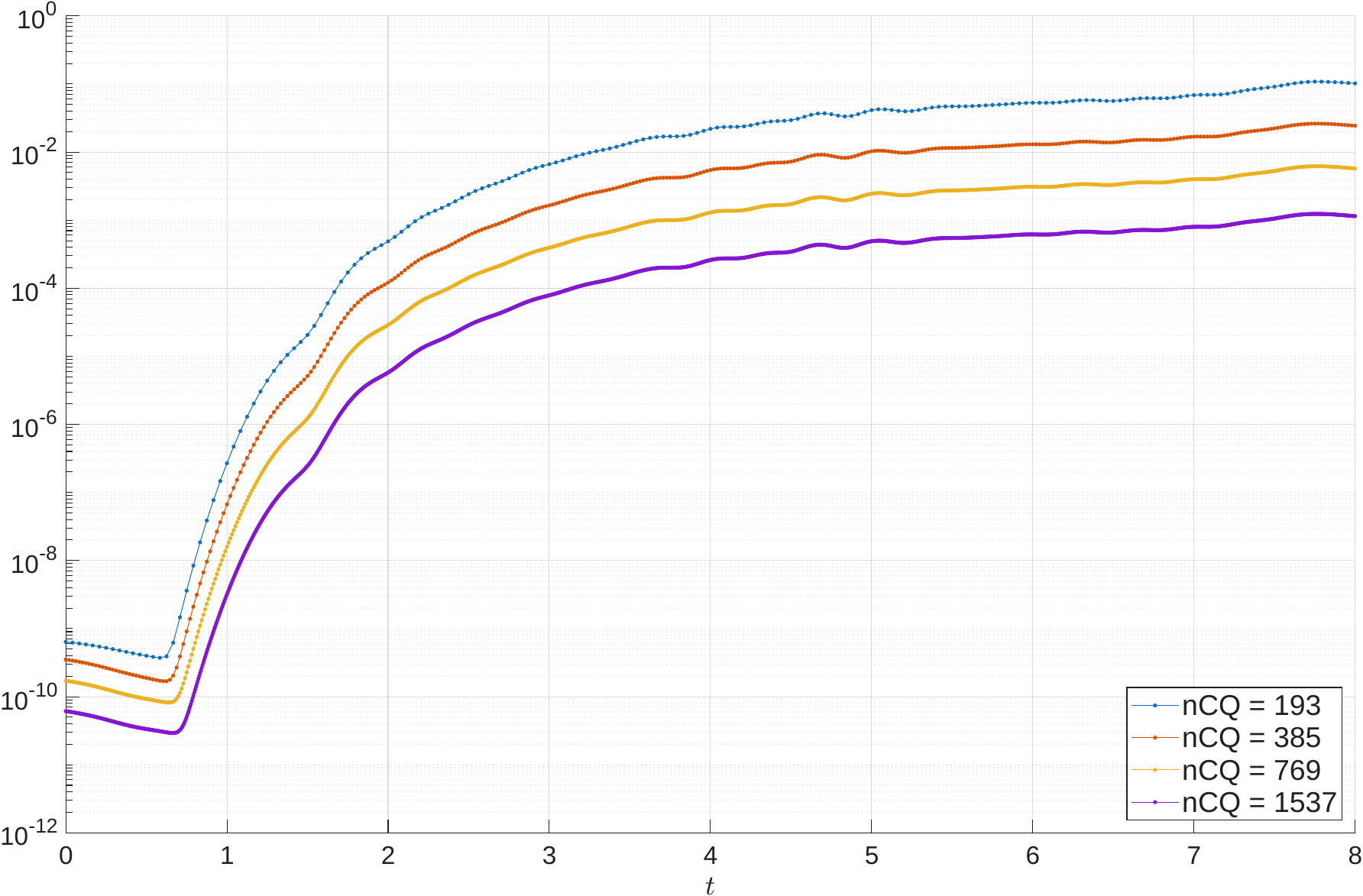}
\caption{
\label{fig:exp01:error}Relative error $E^{\kappa}(t)$ at the time steps $t_j=j\kappa$ for the
BDF2- (left) and trapezoidal-rule-based (right) CQ schemes,
$n_{\rm CQ}=192,\,384,\,768,\,1536$; the reference solution uses the same scheme
with $n_{\rm CQ}=3072$.
Top: a continuous $\mathbb{P}_2$ finite element space on the second uniform
refinement of the initial mesh, \num{36961} nodes (vertices and edge midpoints)
on \num{18336} triangles. Bottom:
two further uniform refinements, \num{587905} nodes on \num{293376} triangles.}
\end{figure}

\begin{figure}[htb]
\centering
\begin{subfigure}[t]{0.85\textwidth}
  \centering
  \includegraphics[width=\textwidth]{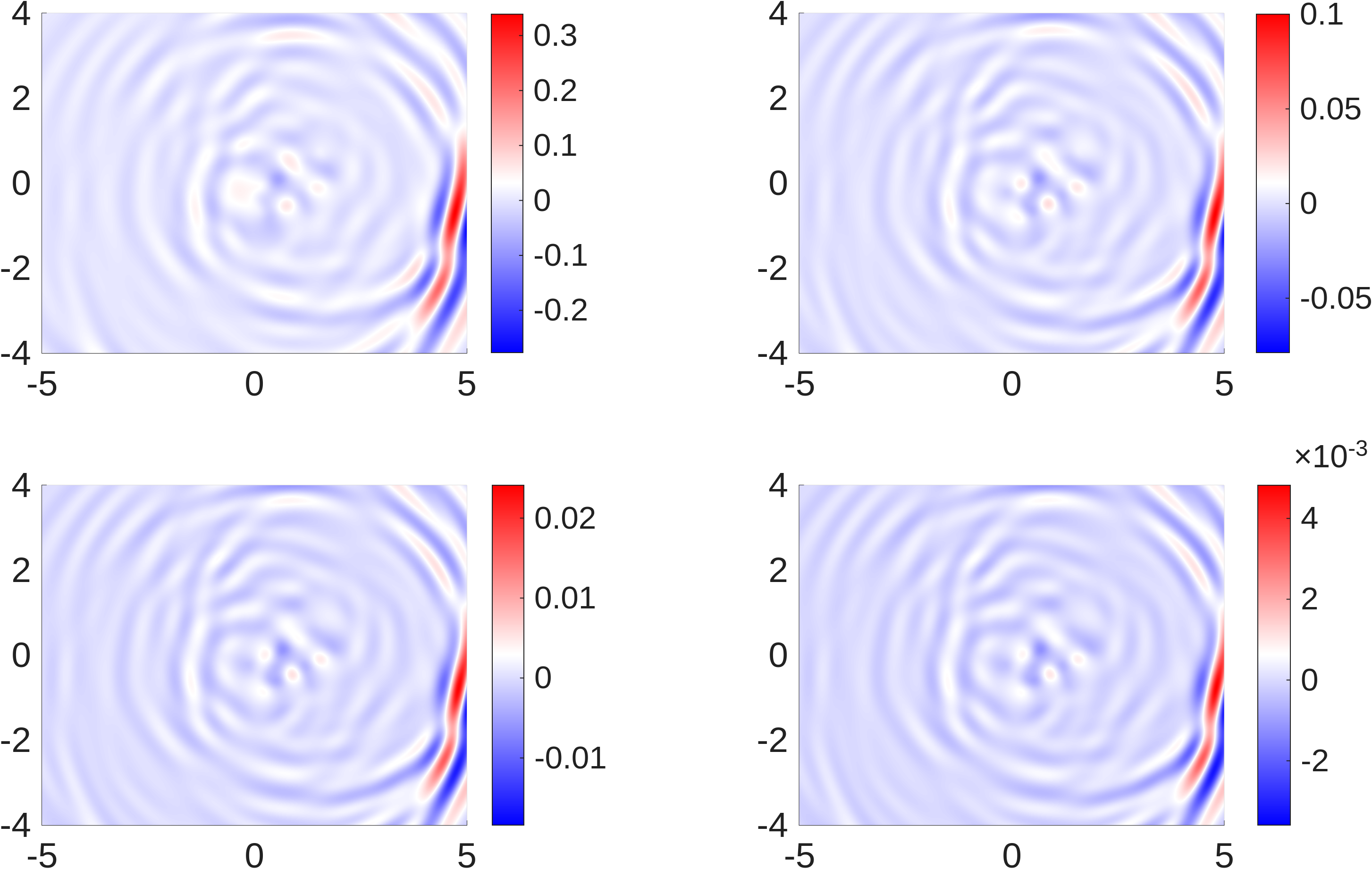}
  \caption{BDF2.}
  \label{fig:exp01:err8:bdf2}
\end{subfigure}
\\[1ex]
\begin{subfigure}[t]{0.85\textwidth}
  \centering
  \includegraphics[width=\textwidth]{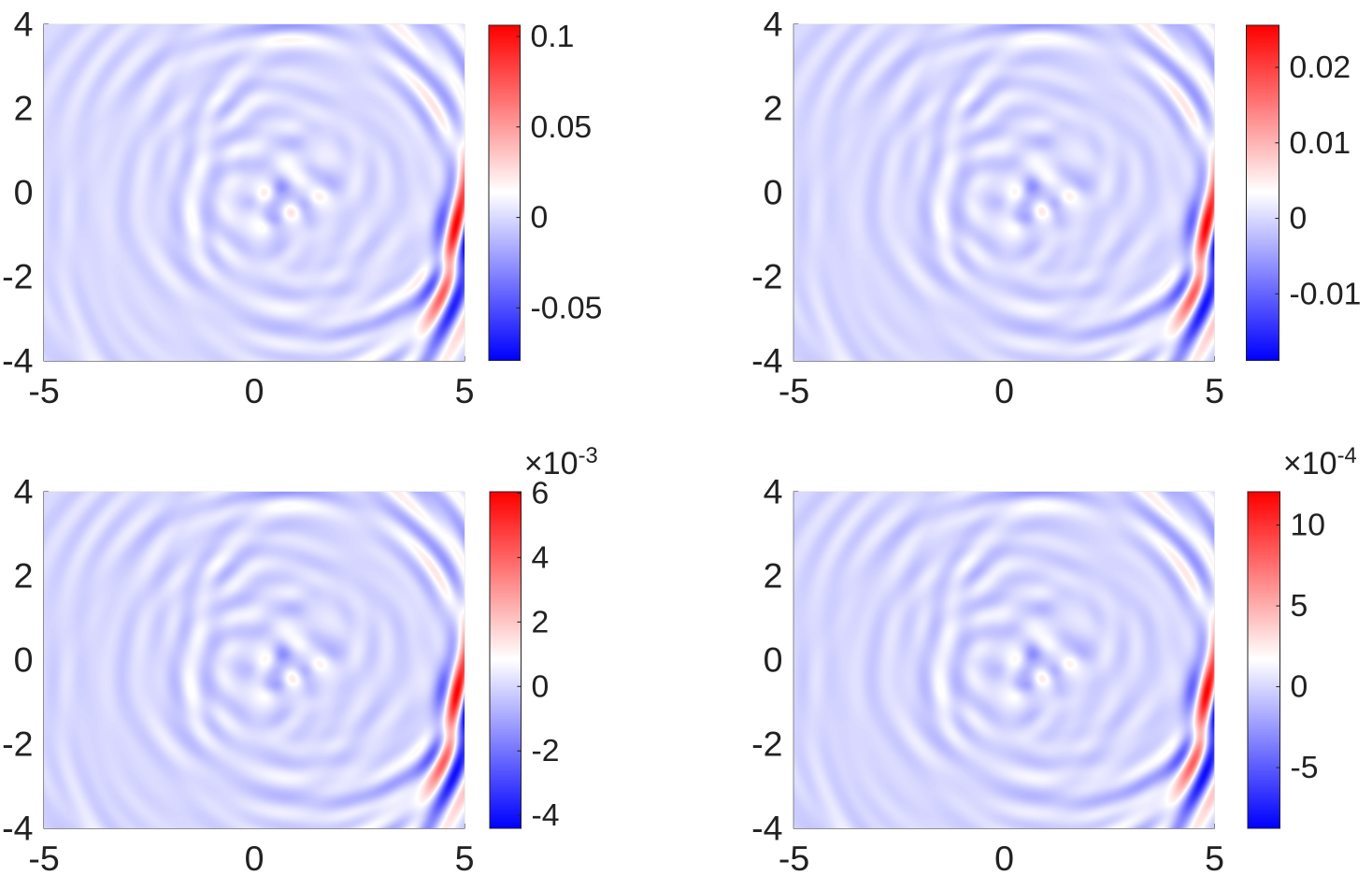}
  \caption{Trapezoidal rule.}
  \label{fig:exp01:err8:tr}
\end{subfigure}
\caption{Error $u^{\kappa}-u^{\rm ref}$ at $t=8$ for
$n_{\rm CQ}=192,\,384,\,768,\,1536$ (left to right), each against a reference
solution of the same scheme with $n_{\rm CQ}=3072$. Note the different color
scales: the error is roughly divided by four at each refinement while its spatial
pattern is unchanged. It concentrates at the incident front, which is leaving
$\Omega_2$ through $\Sigma$ at this time.}
\label{fig:exp01:err8}
\end{figure}

\begin{figure}[htb]
\centering
\includegraphics[width=0.47\textwidth]{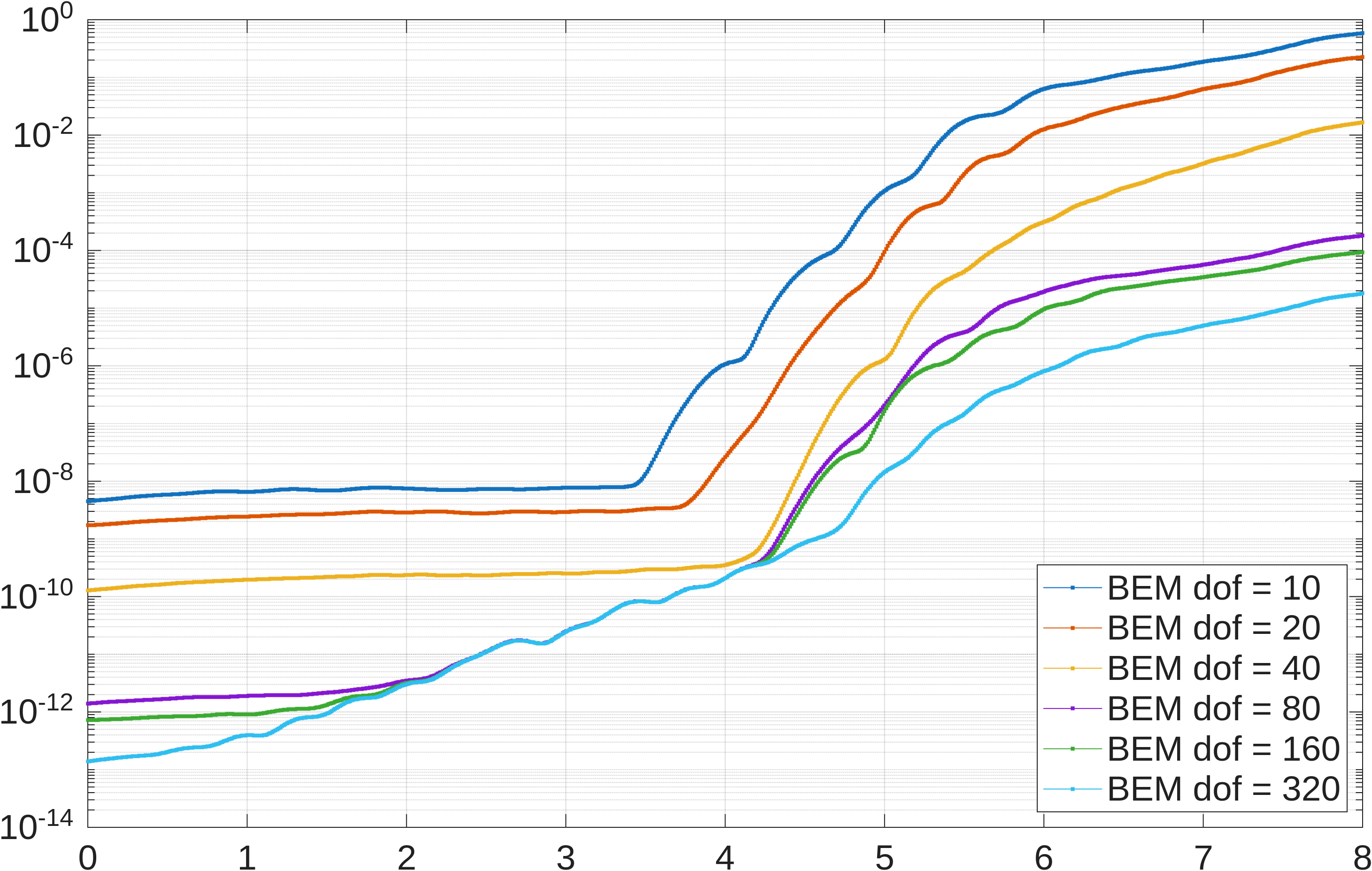}\quad
\includegraphics[width=0.47\textwidth]{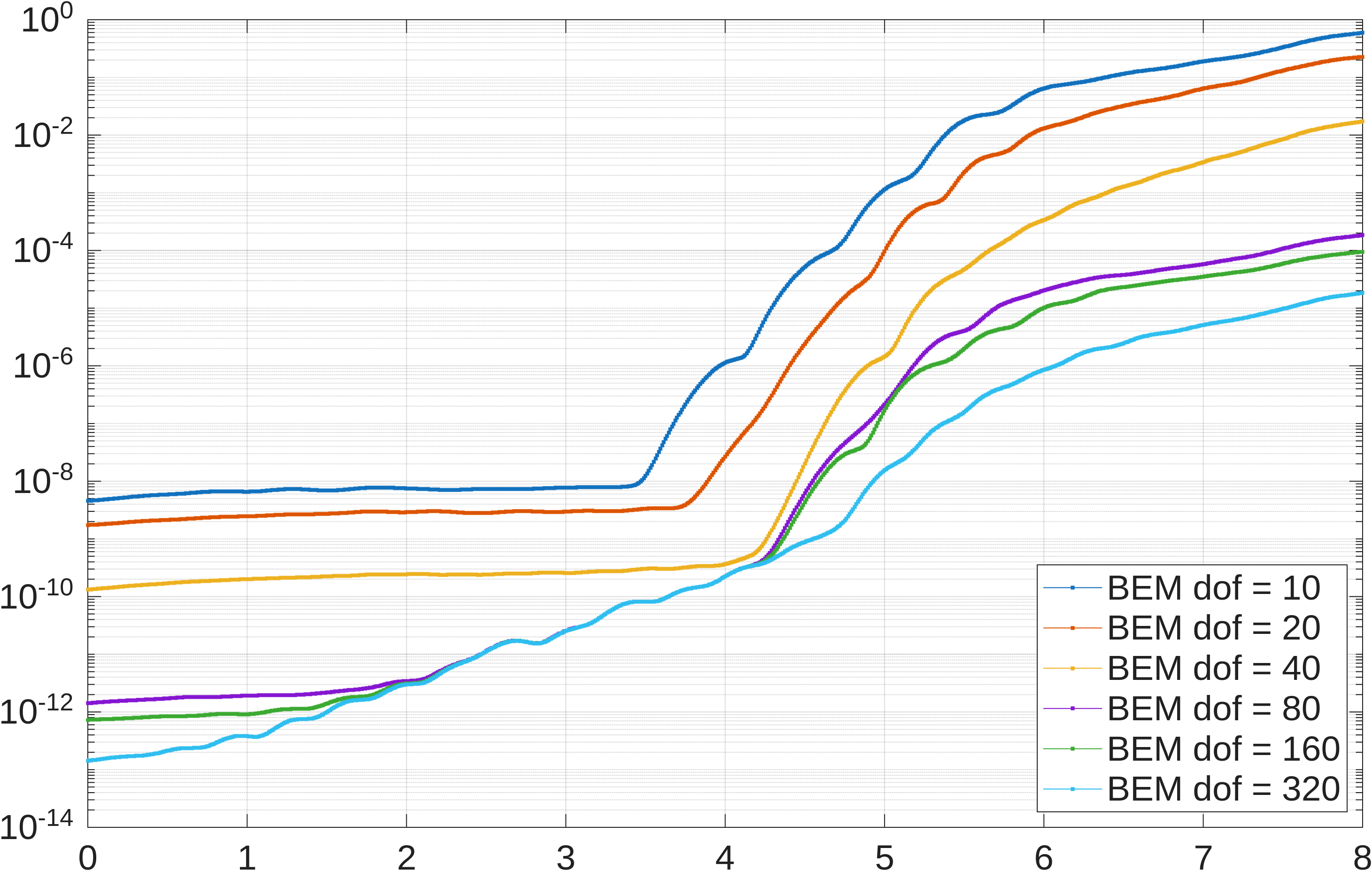}
\caption{Relative error $E^{\kappa}(t)$ at the time steps $t_j=j\kappa$ for the BDF2- (left) and trapezoidal-rule-based (right) CQ schemes with $n_{\rm CQ}=768$ and a $\mathbb{P}_2$ finite element space with \num{36961} nodes on \num{18336} triangles (second uniform refinement). The deltaBEM solver uses 10, 20, 40, 80, 160 and 320 degrees of freedom, the reference solution being computed with 640.}
\label{fig:exp01:BEMerror}
\end{figure}

\begin{figure}[htb]
\centering
\includegraphics[width=0.47\textwidth]{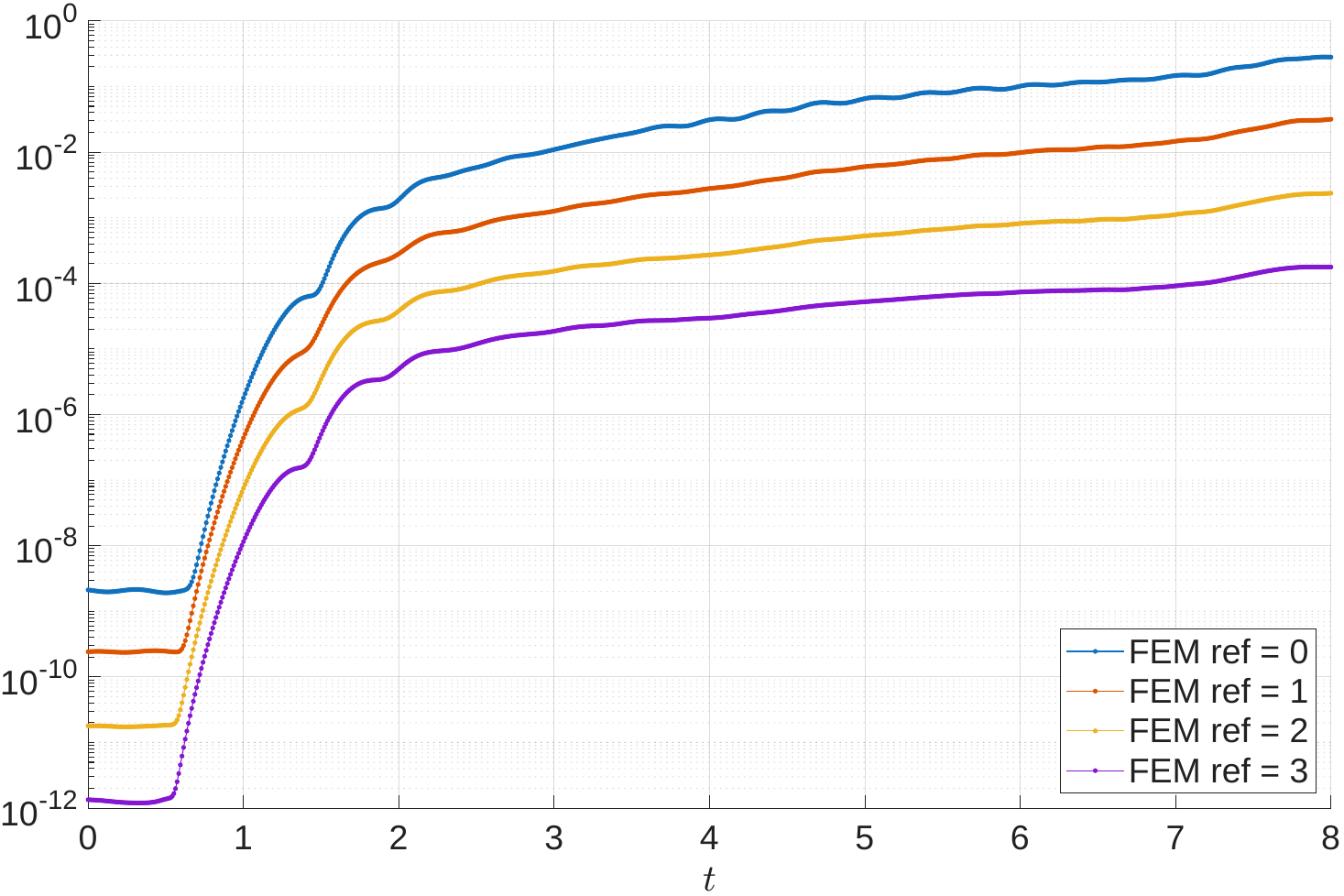}\quad
\includegraphics[width=0.47\textwidth]{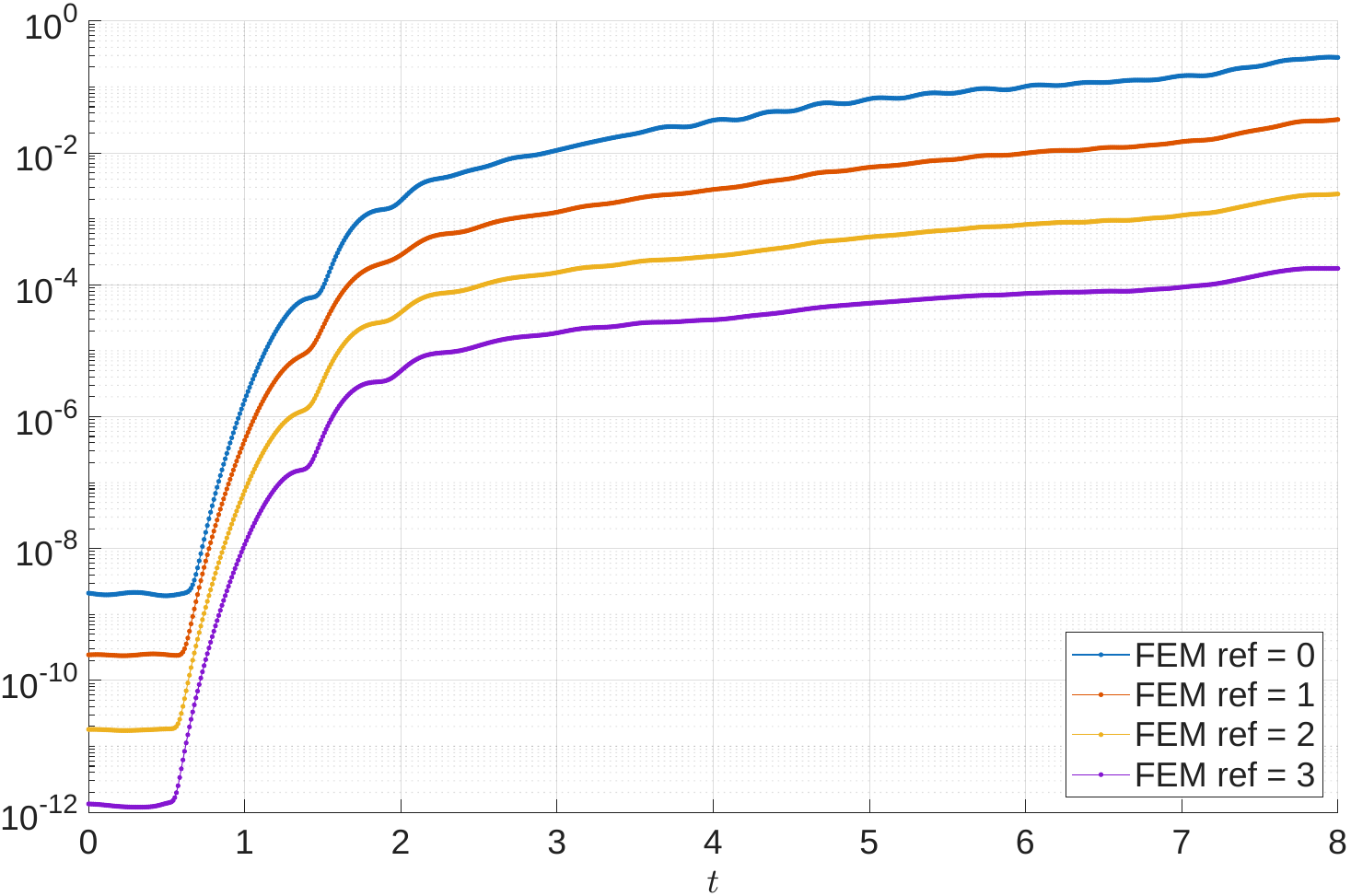}
\caption{Relative error $E^{\kappa}(t)$ at the time steps $t_j=j\kappa$ for the
BDF2- (left) and trapezoidal-rule-based (right) CQ schemes with
$n_{\rm CQ}=768$ and \num{320} degrees of freedom in the deltaBEM solver. The
error is computed from the initial mesh (\num{2365} nodes on \num{1146}
triangles) up to its third uniform  refinement, the reference solution
being computed with the same discretization parameters on the next refinement
(\num{587905} nodes on \num{293 376} triangles).}
\label{fig:exp01:FEMerror}
\end{figure}

\begin{figure}[htb]
\begin{subfigure}[t]{0.85\textwidth}
\centering
  \includegraphics[width=0.9\textwidth]{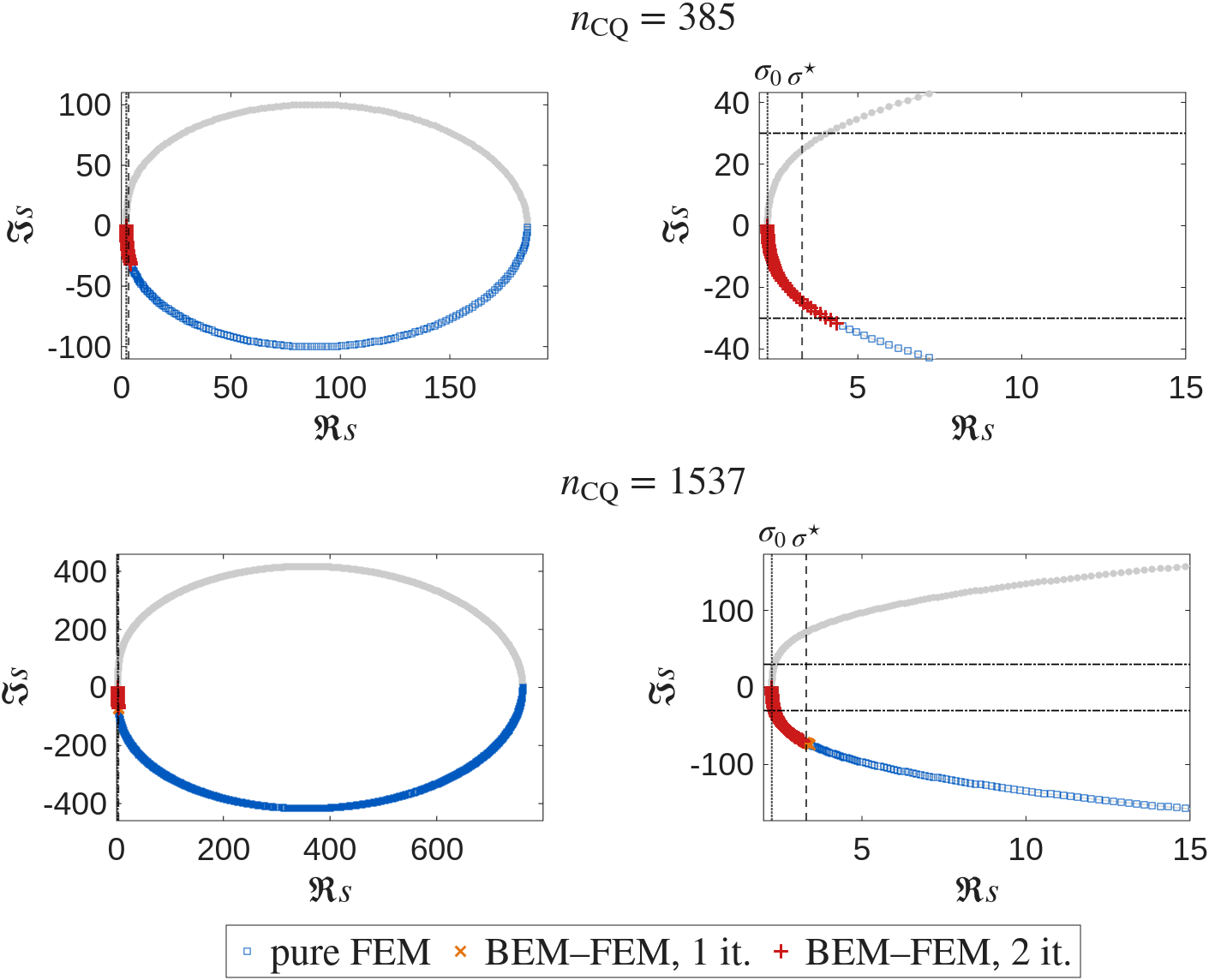}
  \caption{BDF2.}
\end{subfigure}
\begin{subfigure}[t]{0.85\textwidth}
\centering
  \includegraphics[width=0.9\textwidth]{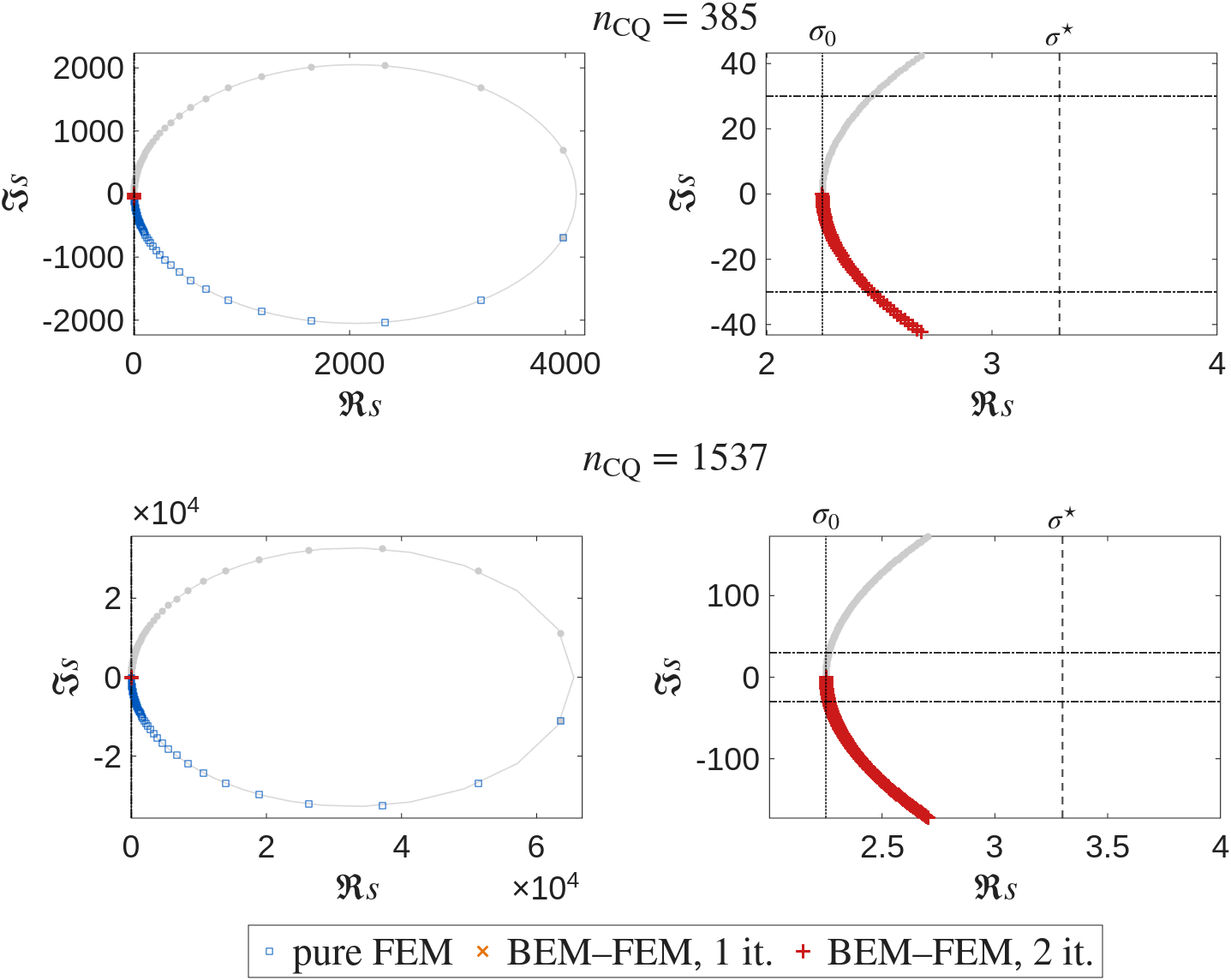}
  \caption{Trapezoidal rule.}
\end{subfigure}
\caption{\label{fig:BEMFEM:iterations}Frequencies $s_\ell = \delta(re_{\ell,n_Q})/\kappa$ used by the
all-steps-at-once solver for $n_{\rm CQ}=384$ and $n_{\rm CQ}=1536$, marked
according to the solver applied: pure FEM (squares) or FEM--BEM coupling with
one or two GMRES iterations (crosses). Right: detail of the region
$\sigma_0\le\Re s\le\sigma^\star$, where the coupling is active. Fourth uniform
refinement of the initial mesh, quadratic continuous finite elements with
\num{587905} nodes (vertices and midpoints) on \num{293376} triangles, and 300
dof on the BEM solver set on $\Gamma$.}
\end{figure}

\begin{figure}[htb]
\begin{subfigure}[t]{0.85\textwidth}
\centering
  \includegraphics[width=0.9\textwidth]{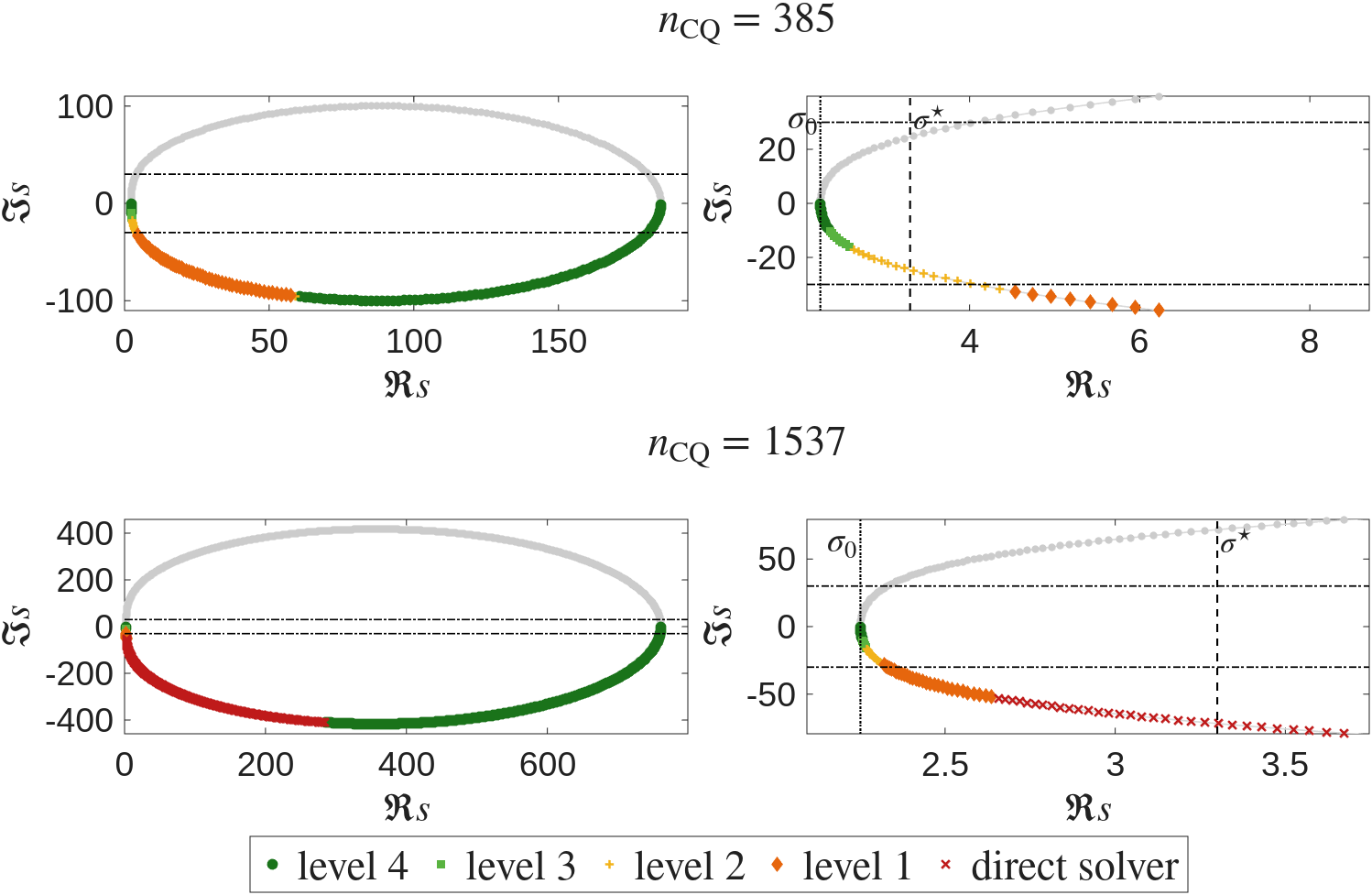}
  \caption{BDF2.}
\end{subfigure}
\begin{subfigure}[t]{0.85\textwidth}
\centering
  \includegraphics[width=0.9\textwidth]{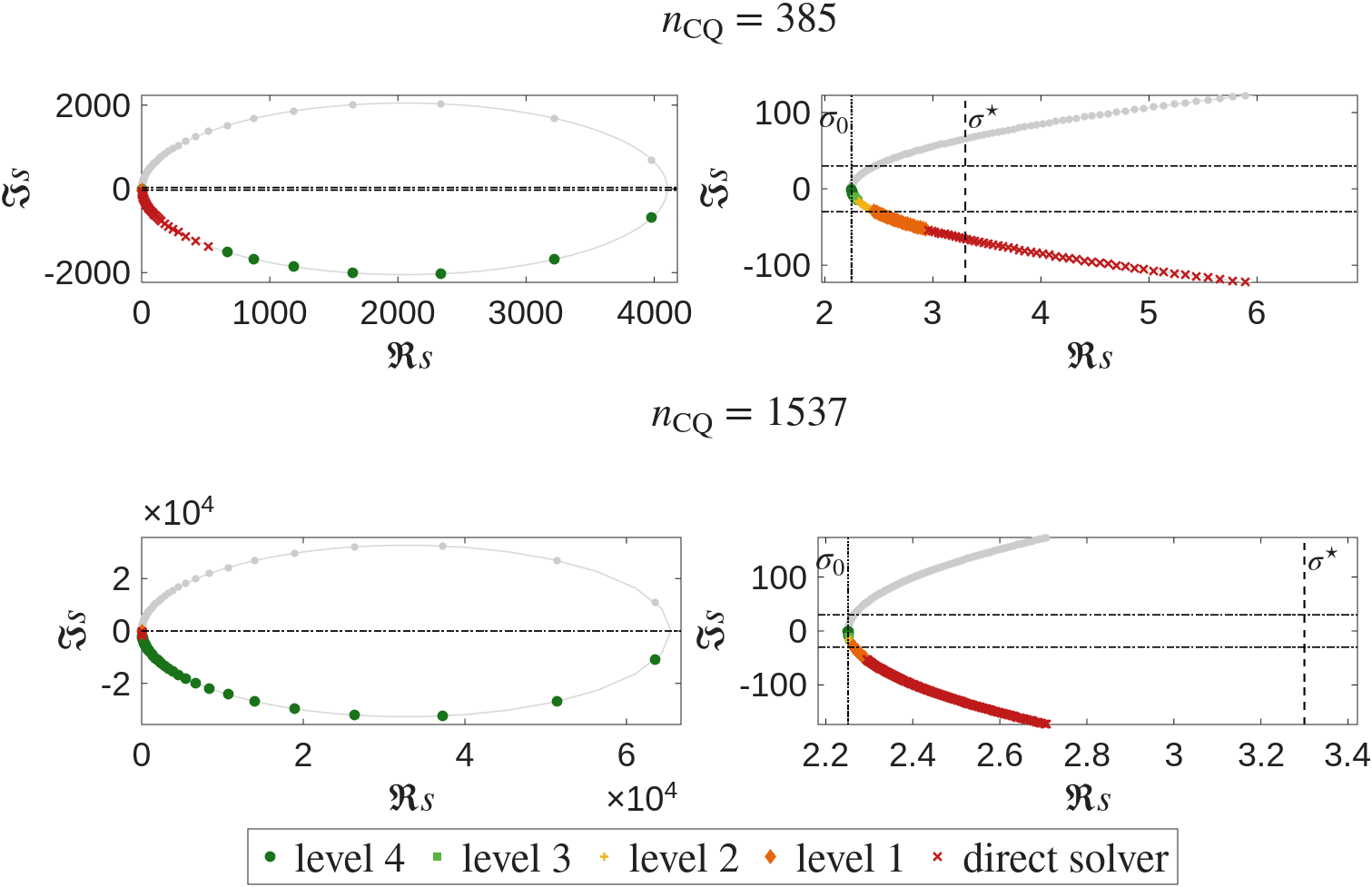}
  \caption{Trapezoidal rule.}
\end{subfigure}
\caption{\label{fig:BEMFEM:multigrid}Frequencies $s_\ell = \delta(re_{\ell,n_Q})/\kappa$ used by the all-steps-at-once solver for $n_{\rm CQ}=384$ and $n_{\rm CQ}=1536$, marked according to the number of levels of the V-cycle that converged, the cross denoting the frequencies for which the multigrid solver failed and was replaced by a sparse direct solver. Right: detail of the region of small $\sigma=\Re s$. Same discretization as in Figure~\ref{fig:BEMFEM:iterations}.}
\end{figure}

\paragraph{Results.}
We have run the experiment with the BDF2- and the trapezoidal-rule-based CQ
schemes for $n_{\rm CQ}=192,\,384,\,768,\,1536$,\footnote{In the legends and
titles of Figures~\ref{fig:exp01:error}, \ref{fig:BEMFEM:iterations} and
\ref{fig:BEMFEM:multigrid}, $n_{\rm CQ}$ counts the time levels including
$t=0$, that is, $n_{\rm CQ}+1$ in the notation of Section~\ref{sec:FULLCQ}.}
i.e.\ $\kappa=1/24,\,1/48,\,1/96,\,1/192$, with continuous $\mathbb{P}_2$
finite elements on the second and fourth uniform refinements of the initial
mesh (\num{36961} and \num{587905} degrees of freedom, respectively) and 300
boundary element unknowns on $\Gamma$. Figure~\ref{fig:exp01:total} shows
the total field computed with the trapezoidal rule and $n_{\rm CQ}=3072$, the
discretization used below as reference solution, at four instants, and
Figure~\ref{fig:exp01:scatt} the corresponding scattered field.

Figure~\ref{fig:exp01:error} shows $E^{\kappa}(t)$ for both schemes, with the
same BEM discretization (single layer formulation with 300 degrees of freedom)
and two levels of refinement of the finite element space. The curves are
practically indistinguishable, 
which indicates that the time error does not depend on the spatial resolution. For illustrative purposes we also depict the error
at the final time $T=8$ in Figure~\ref{fig:exp01:err8}.

The maximum in time of the relative error $E^\kappa_\infty$ (see \eqref{def:Ekappa:infty}) decays from \num{0.31} to \num{4.4e-3} for BDF2 and from \num{0.1} to \num{1.2e-3} for the trapezoidal rule, with estimated orders of convergence $\log_2(E^{\kappa}_\infty/E^{\kappa/2}_\infty)$ of $1.6,\ 2.3,\ 2.1$ and $2.0,\ 2.1,\ 2.3$, respectively. The last value of each sequence is slightly biased upwards, since the reference solution has only half the time step of the finest discretization.

Both schemes exhibit the second order rate of convergence predicted by Theorem~\ref{th:main}, but they do so in different ways. The trapezoidal rule is in the asymptotic regime from the coarsest step onwards, whereas BDF2 shows a reduced order between the two coarsest levels (EOC $\approx 1.6$) and reaches second order only for $\kappa\le 1/48$. The latter is a pre-asymptotic effect: a relative error of $0.31$ is close to saturation, and the error can no longer scale with $\kappa^2$ at that level.

The comparison between the two schemes is otherwise governed by the phase error, which for a mode of frequency $\omega$ accumulates at the rate $\kappa^2\omega^3/12$ per unit time for the trapezoidal rule and $\kappa^2\omega^3/3$ for BDF2, a factor of four \cite[Sect.~4.13]{BanjaiSayas2023}; accordingly, the errors of the trapezoidal rule are 3--4 times smaller than those of BDF2 at equal $\kappa$ throughout the whole range, or, equivalently, the trapezoidal rule with step $\kappa$ matches BDF2 with step $\kappa/2$, that is, with twice as many frequency problems. This is consistent with the recommendation of \cite[Sect.~4.13]{BanjaiSayas2023} in favor of the trapezoidal rule on accuracy grounds. 

We also consider the influence of the spatial discretization, keeping $n_{\rm CQ}=768$ fixed. Figure~\ref{fig:exp01:BEMerror} shows $E^{\kappa}(t)$ for an increasing number of degrees of freedom on $\Gamma$, the finite element space being fixed. Up to $t\approx4$ the scattered field has not yet reached $\Gamma$, the coupling step is skipped at almost all frequencies, and the curves merely reflect the round-off level of the scaling by $r$; the boundary element discretization only comes into play afterward. From then on the error at the final time decreases rapidly, by one to two orders of magnitude each time the number of degrees of freedom is doubled, from about $0.5$ with 10 to about $10^{-4}$ with 80, and then levels off around $10^{-5}$ for 160 and 320. This plateau is three orders of magnitude below the time discretization error at this $n_{\rm CQ}$ and may be caused by the tolerances of the coupling step rather than by the boundary element resolution. With 320 degrees of freedom the boundary element contribution therefore plays no role in the results reported above.

Keeping instead the boundary discretization fixed and refining the finite element space (Figure~\ref{fig:exp01:FEMerror}), the estimated orders of convergence are about $3$. This is the $L^2$ rate of $\mathbb{P}_2$ elements for piecewise $H^3$ solutions, and it is the rate one expects here: the mesh resolves $\partial\Omega$ and the interfaces between the pieces of $\Omega$, so that the jumps of $c$ and $F$ across them do not affect the elementwise interpolation error, while the corners of these curves, where $\Delta U$ is discontinuous, and the right-angle corners of $\Sigma$ \cite[Section 5.1]{Grisvard1985} only reduce the regularity to $H^{3-\varepsilon}(\Omega_2)$, piecewise, for every $\varepsilon>0$.

\paragraph{Iterative solvers for the Schur complement equation.}
The number of GMRES iterations for the Schur complement equation
\eqref{eq:the:equation} is usually one or two, as Figure~\ref{fig:BEMFEM:iterations}
shows for two representative runs, and we never observed more than four in any
of our experiments; at most frequencies the coupling step is skipped altogether.
This is better than what Theorem~\ref{theo:the:theo:02} guarantees, since the
contraction estimate \eqref{eq:contraction:factor} does not hold in the weakly
damped window where the coupling is active, see Remark~\ref{remark:considerations}.
As mentioned above, the coupling is activated whenever $\gamma_\Sigma\mathrm K_{\Omega_1^{\rm c}\Gamma}(s)f_\Gamma$,
the first iterate of \eqref{eq:the:equation}, exceeds $10^{-7}\|W^{\rm tr}_F\|_\infty$.
By \eqref{eq:02:theo:the:theo:01} and \eqref{eq:01b:theo:the:theo:02} this quantity
decays like $e^{-\sigma(d_1+\min\{d_1,d_2\})}$ up to polynomial factors, so that
the coupling is needed only for $\Re s\le\sigma^\star$, with $\sigma^\star\approx3.3$
fixed by the geometry and the tolerance. On the other hand, the scaling radius $r={\tt eps}^{1/(2n_{\rm CQ})}$ keeps every frequency at a distance from the imaginary axis that is independent of $\kappa$. Indeed, for the three methods considered $\Re\delta(re^{{\rm i}\theta})$ attains its minimum over $\theta$ at $\theta=0$, so that $\Re s_\ell\ge\delta(r)/\kappa=:\sigma_0$ for all $\ell$, and since $\delta(e^{-x})=x+O(x^2)$,
\[
 \sigma_0=\frac{\delta(r)}{\kappa}\approx\frac{|\ln{\tt eps}|}{2n_{\rm CQ}}\cdot\frac{n_{\rm CQ}}{T}=\frac{|\ln{\tt eps}|}{2T}\approx2.25 .
\]

The BEM is therefore needed only in
the fixed window $\sigma_0\le\Re s\le\sigma^\star$, marked by the dashed lines in
the right panels of Figure~\ref{fig:BEMFEM:iterations}, and the proportion of pure
FEM problems is $80\%$--$88\%$ for BDF2 and, much more steadily, $57\%$--$60\%$
for the trapezoidal rule.


The frequencies that actually contribute to the solution can be identified from the
transformed data. Writing the quadrature node as $re_{\ell,n_{\rm CQ}}=re^{{\rm i}\theta_\ell}$,
$\theta_\ell=2\pi\ell/(n_{\rm CQ}+1)$, and using $r\approx e^{-\kappa\sigma_0}$, step~\ref{step:2}
of Algorithm~\ref{alg:allsteps} reads
\begin{equation}\label{eq:which:F:matters}
  \begin{aligned}
 F^\kappa_{\ell,r}\ & =\ \frac1{n_{\rm CQ}+1}\sum_{j=0}^{n_{\rm CQ}}f(j\kappa)\,e^{-(\sigma_0-{\rm i}\tau_\ell)\,j\kappa}
 \\
 & \ \approx\ \frac1T\int_0^T f(t)\,e^{-\sigma_0t}\,e^{{\rm i}\tau_\ell t}\,{\rm d}t,\quad \tau_\ell:=\frac{\theta_\ell}{\kappa}=\frac{2\pi\ell}{T+\kappa}\approx \frac{2\pi\ell}{T},
  \end{aligned}
\end{equation}
that is, heuristically and up to the factor $1/T$, the Fourier transform of the damped
datum at the frequencies $\tau_\ell$ of a Fourier series on $[0,T]$, whatever the multistep
method. For the pulse \eqref{eq:exp01:uinc} this transform is, up to smooth factors, the
sum of two Gaussians $e^{-\varsigma^2(\tau\mp\omega)^2}$ of width $1/\varsigma=5$
centered at $\pm\omega$, which falls below $10^{-13}$ of its maximum for
$|\tau|\gtrsim30$. Hence only the nodes with $\tau_\ell\lesssim30$, about $30T/\pi\approx75$
of them regardless of $n_{\rm Q}$, carry a nonnegligible right--hand side; refining the
time step only adds nodes with data at round--off level. For these nodes
$\delta(re_{\ell,n_{\rm CQ}})\approx\kappa\sigma_0-{\rm i}\theta_\ell$, so that
$s_\ell\approx\sigma_0-{\rm i}\tau_\ell$ (their conjugates give the other half of the contour):
they are the least damped frequencies of the contour, the ones for which the coupling is
needed and, as discussed next, the ones for which the interior solver is hardest. All
the others, however many, are cheap.

\paragraph{Iterative solvers for the FEM problems.}
Most of the FEM problems are in fact solved for completeness of the FFT step only: for the pulse \eqref{eq:exp01:uinc}, $F^\kappa_{\ell,r}$ is smaller than $10^{-13}$ relative to its maximum for $|\Im s|\gtrsim 30$, as
discussed in \eqref{eq:which:F:matters}, and these frequencies could be skipped altogether by setting $Y^\kappa_{\ell,r}=0$. Moreover, for such strongly damped frequencies the solution is close to $c^2F/s^2$, a rescaled copy of the datum, up to boundary layers of width $O(1/|s|)$ along the interfaces where $F$ jumps; the mesh does not resolve these layers, but they carry a negligible part of an already negligible solution and do not propagate.

Regarding the FEM solver, and motivated by the damping present in the complex wave number, we have implemented a multigrid solver based on V-cycles with Gauss--Seidel smoothing: one pre-smoothing step, no post-smoothing, undamped coarse-grid correction. The
depth of the V-cycle is therefore limited by the number of levels available in
each mesh hierarchy. Let us stress that designing a multigrid solver for this problem is not among the aims of this paper: a naive implementation suffices for our purposes, and what follows is merely an exploration of the terrain, intended to show that iterative solvers are a realistic option here and to identify the regimes in which they can be expected to work. Its interest is not speed --- on the present two-dimensional problem it is comparable to a sparse direct solver --- but memory: no fill-in is stored, only the matrices of the mesh hierarchy. In three dimensions, where sparse direct solvers cease to be affordable, some such alternative is not merely convenient but unavoidable. It should also be emphasized that we are in a far more comfortable situation than that of the genuine Helmholtz equation: the wave number $k={\rm i}s$ has a positive imaginary part bounded away from zero, uniformly in $\kappa$, and it is this damping that makes a plain V-cycle work at all for a sizable part of the frequencies.


Figure~\ref{fig:BEMFEM:multigrid} shows where the V-cycle converged and with how many levels. Two features are worth pointing out. First, the number of levels that can be used decreases as $|s|$ grows, one level being lost roughly each time $|s|$ doubles, which is consistent with the coarsest mesh of the hierarchy having to resolve the length scale $1/|s|$; the thresholds are the same for both schemes, as one would expect from a mechanism that depends on the spatial mesh alone. Second, for sufficiently large $|s|$ and insufficient relative damping the smoother fails, and a sparse direct solver is used instead. This is consistent with the analysis of the shifted Helmholtz operator introduced in \cite{ErOoVu:2006}: a complex shift of size $O(k^{2})$, that is, of relative size $O(1)$, is required for this operator to be invertible by multigrid \cite{CoGa:2017}, whereas a shift of size at most $O(k)$ is needed for it to be a good preconditioner for the unshifted problem \cite{GanderGrahamSpence2015}. In our setting the second requirement plays no role, since the shifted problem is the one to be solved. In terms of the Laplace parameter, $k^{2}=-s^{2}$ and the relative size of the shift is $|\Im k^{2}|/|\Re k^{2}|\approx 2\,\Re s/|s|$, so that the threshold observed here corresponds to a shift of the order predicted there.

The proportion of frequencies for which the smoother fails grows as the time step decreases --- from $0\%$ to $48\%$ for BDF2 between $n_{\rm CQ}=384$ and $1536$ --- and is consistently larger for the trapezoidal rule, $63\%$ and $78\%$ for the same two values, whose contour accumulates a large number of poorly damped frequencies. 
The two contours differ where it matters for the smoother. The relative damping $\sigma/|s|=\Re\delta(z)/|\delta(z)|$ depends on the quadrature node $z=re^{{\rm i}\theta}$ alone, $\theta=2\pi\ell/(n_{\rm Q}+1)$, and near $\theta=0$ it is $\approx(1-r)/\theta$ for both schemes, with $1-r\approx\kappa\sigma_0$ forced by round--off: in the band of frequencies that carry the data the two schemes are equally hard, and both resort there to the direct solver or to few levels. Beyond that band they part ways. For BDF2, 
\[
\frac{\sigma}{|s|}\approx \frac1{\theta}\big((1-r)+\frac{1}4\theta^4\big),
\] which reaches a minimum of order $(1-r)^{3/4}$ and then grows to $O(1)$: the weakly damped frequencies form a narrow band, and the remaining failures, at large $|s|$, are due to the coarsest meshes not resolving the length scale $1/|s|$, an \emph{inverse} CFL condition $h\lesssim\kappa$ that constrains the robustness of the iterative solver but not the stability of the scheme. For the trapezoidal rule, instead, 
\[ 
\frac{\sigma}{|s|}\approx\frac{1-r}{\sin\theta},
\] 
which stays below $2(1-r)$ on two thirds of the nodes. A small relative damping can still be compensated by resolution, since the coarse--grid correction of the near--resonant modes requires the coarsest mesh to resolve the wavelength to a degree that grows as the damping decreases \cite{CoGa:2017},  but for the trapezoidal rule this would have to hold on most of the contour, where $|s|$ is already of order $\kappa^{-1}$ while the relative damping is only of order $\kappa$, which is not a practical option; see Figure~\ref{fig:BEMFEM:multigrid}(b). An iterative interior solver for this scheme would rather require smoothers adapted to the weakly damped Helmholtz regime. This is the situation of the classical Helmholtz equation and remains a challenge.

Table~\ref{tab:counts} gives a detailed account of the frequency problems behind the time-convergence experiment of Figure~\ref{fig:exp01:error}. The second column is the number of independent frequencies generated by the quadrature once conjugate symmetry is taken into account, $n_{\rm CQ}/2+1$, which grows linearly with the number of time steps. The third column shows that only 40 of them were actually solved, in all runs and for both schemes: these are the frequencies with the largest data, all with $|\Im s|\lesssim30$, and the remaining ones, whose data lie below $10^{-11}$ of the maximum, were set to zero. Since the convergence rates of Figure~\ref{fig:exp01:error} are unaffected, the cost of the solver phase does not depend on $n_{\rm CQ}$; what grows with the number of time steps is only the FFT and the bookkeeping of the discarded frequencies. The fourth column counts the problems in which the BEM correction was activated. For the trapezoidal rule this happens in all 40 problems, which lie within the window $\sigma_0\le\Re s\le\sigma^\star$ of Figure~\ref{fig:BEMFEM:iterations}; for BDF2 and the two coarsest time steps some of them fall outside it, since for large $\kappa$ the nodes carrying the data spread over a wider arc of the contour, where the BDF2 damping is larger, and are solved as pure FEM problems. In every case the coupled system was solved with two GMRES iterations. The last columns describe the interior solver. All 40 problems were solved by the multigrid V-cycle, at full depth or after reducing the number of coarse levels, and the direct solver on the finest mesh was never needed. The distribution of depths is the same for all time steps except the coarsest one, for the reason discussed below.

\begin{table}[htb]
\centering
{
\begin{tabular}{lrrccrrrrr}
\hline
 & $n_{\rm CQ}$ & frequencies & solved & BEM corr. & \multicolumn{4}{c}{multigrid depth} & LU\\
 & & & & & 4 & 3 & 2 & 1 & \\
\hline
BDF2        &  192 &   97 & 40 & 24 & 15 & 12 & 13 & 0 & 0\\
            &  384 &  193 & 40 & 38 & 14 & 10 & 16 & 0 & 0\\
            &  768 &  385 & 40 & 40 & 14 & 10 & 14 & 2 & 0\\
            & 1536 &  769 & 40 & 40 & 14 &  9 & 15 & 2 & 0\\
            & 3072 & 1537 & 40 & 40 & 14 &  9 & 14 & 3 & 0\\
\hline
Trapezoidal &  192 &   97 & 40 & 40 & 14 & 10 & 12 & 4 & 0\\
            &  384 &  193 & 40 & 40 & 14 &  9 & 14 & 3 & 0\\
            &  768 &  385 & 40 & 40 & 14 &  9 & 14 & 3 & 0\\
            & 1536 &  769 & 40 & 40 & 14 &  9 & 14 & 3 & 0\\
            & 3072 & 1537 & 40 & 40 & 14 &  9 & 14 & 3 & 0\\
\hline
\end{tabular}}
\caption{\label{tab:counts}Frequency problems in the time-convergence experiment, for each scheme and time step: independent frequencies after conjugate symmetry ($n_{\rm CQ}/2+1$); problems actually solved, the remaining ones carrying data below $10^{-11}$ of the maximum and being set to zero; problems in which the BEM correction was activated; problems in which the multigrid V-cycle converged with four, three, two or one coarse levels (Gauss--Seidel smoothing, LU on the coarsest level after a first failure at full depth); and problems that required a direct solver on the finest mesh.}
\end{table}

\section{Conclusions}
\label{sec:conclusions}
We have proposed and analyzed a fully CQ method for efficient simulation of transient acoustic waves in a heterogeneous and unbounded medium. Multistep CQ applied to the solution operator reduces the problem to independent Helmholtz--type problems with complex frequencies, each of them solved by an overlapping FEM--BEM coupling in which the finite element domain and the boundary element curve are chosen independently of the geometry of the inhomogeneity. The analysis of the time discretization rests on bounds for the solution operators of the transformed problems that are explicit in the Laplace parameter and in the wave speed; the exponential decay of the transformed solutions away from the inhomogeneity shows, in turn, that the boundary element correction is needed only in a fixed window of weakly damped frequencies, independent of the time step, and that the coupled operator is a contraction outside it.

The numerical experiments confirm the second order in time of the BDF2 and trapezoidal schemes, with the errors in time and in space decoupled from each other, and show that the coupled equation is solved by GMRES in very few iterations, usually one or two, at every frequency where the coupling is active, a behavior that our contraction estimate does not yet explain. They also show that a plain multigrid V--cycle solves the interior problems for a large part of the frequencies, failing only where the relative damping is small, in agreement with the known theory for the shifted Helmholtz operator.

Several questions remain open. The analysis of the spatial discretization, with constants explicit in the Laplace parameter, is the natural next step; the regularity of the transformed solutions on polygonal domains and the $s$--dependence of the lifting operators on $\Sigma$ are the main difficulties. Higher order in time calls for Runge--Kutta convolution quadrature, which fits the same all--steps--at--once framework. In three dimensions, where the method is most attractive, the interior solver becomes the bottleneck, and the design  of smoothers adapted to the weakly damped frequencies deserves further study. Finally, the extension to the acoustic model with variable density, where the interior operator takes divergence form while the coupling is unaffected, is left for future work.

\section*{Acknowledgements}

vD is deeply grateful to Francisco~J.~Gaspar for several enlightening discussions
on multigrid methods for wave problems, which greatly improved his understanding
of this iterative method in this context and led to a sound implementation in the present
setting, and for his generous and selfless help.

\section*{Declarations}

\paragraph{Funding.}
The first author (vD) is partially supported by the project
``T\'ecnicas innovadoras para la resoluci\'on de problemas evolutivos'',
grant PID2022-136441NB-I00, funded by MCIN/AEI/10.13039/501100011033
(Ministerio de Ciencia e Innovaci\'on, Government of Spain), and by the project
``Adquisici\'on de conocimiento y miner\'{\i}a de datos, funciones especiales y
m\'etodos num\'ericos avanzados'' of Universidad P\'ublica de Navarra.
The second author (MG) is supported by the Simons Foundation through Grant No.~518882.

\paragraph{Competing interests.}
The authors have no competing interests to declare.

\paragraph{Data availability.}
The MATLAB code and the data generated in the numerical experiments are
available from the corresponding author on reasonable request.

\paragraph{AI use declaration.}
During the preparation of this work the first author (vD) used Claude (Anthropic),
mainly in the debugging and documentation of the MATLAB code used in the numerical
experiments, for checking some of the computations and, to a lesser extent, revising
parts of the text. The authors reviewed and edited the content as needed and take
full responsibility for the content of the publication.

\bibliography{biblio}

\end{document}